\documentclass[a4paper,10pt,notitlepage,reqno]{article}
\usepackage[english]{babel}
\usepackage{amssymb,amsthm,bm} 
\usepackage{mathtools}
\usepackage{mathrsfs}
\usepackage{enumitem}
\usepackage{authblk}
\usepackage{color}
\usepackage{stackrel}
\usepackage{geometry}
\usepackage{hyperref}
\theoremstyle{plain} 
\newtheorem{theorem}{Theorem}[section]
\newtheorem{lemma}{Lemma}[section]

\theoremstyle{definition}

\newtheorem{remark}{Remark}[section]

\numberwithin{equation}{section}

\DeclarePairedDelimiter{\abs}{\lvert}{\rvert} 
\DeclarePairedDelimiter{\norm}{\lVert}{\rVert}

\DeclareMathOperator{\realpart}{Re}
\renewcommand{\Re}{\realpart}

\DeclareMathOperator{\imaginarypart}{Im}
\renewcommand{\Im}{\imaginarypart}

\DeclareMathOperator{\sgn}{sgn}
\DeclareMathOperator{\curl}{curl}
\DeclareMathOperator{\supp}{supp}
\DeclareMathOperator{\dist}{dist}

\newcommand{\R}{\mathbb{R}}
\newcommand{\C}{\mathbb{C}}

\newcommand{\half}{\mbox{$\frac{1}{2}$}}
\newcommand{\quarter}{\mbox{$\frac{1}{4}$}}
\newcommand{\Dom}{\mathcal{D}}

\newcommand{\normeq}[1]{{\left\vert\kern-0.25ex\left\vert\kern-0.25ex\left\vert #1 
    \right\vert\kern-0.25ex\right\vert\kern-0.25ex\right\vert}}

\newenvironment{system}%
{\left\lbrace\begin{array}{r@{\hspace{1mm}}ll}}%
{\end{array}\right.}

\usepackage{color}
\usepackage[normalem]{ulem}
\definecolor{DarkGreen}{rgb}{0,0.5,0.1} 

\newcommand\soutD{\bgroup\markoverwith
{\textcolor{DarkGreen}{\rule[.5ex]{2pt}{1pt}}}\ULon}

\definecolor{Darkgblue}{rgb}{0.3,0.3,0.5}

\title{\textbf{Absence of eigenvalues of Dirac and Pauli Hamiltonians via the method of multipliers}}

\author[1]{Lucrezia Cossetti} 
\author[2]{Luca Fanelli}
\author[3]{David Krej\v{c}i\v{r}\'ik}

\affil[1]{Fakult\"{a}t f\"{u}r Mathematik, Institut f\"{u}r Analysis, Karlsruher Institut f\"{u}r Technologie, Englerstra{\ss}e 2, 76131 Karlsruhe, Germany; lucrezia.cossetti@kit.edu}

\affil[2]{Dipartimento di Matematica, SAPIENZA Universit\`{a} di Roma, P. le Aldo Moro 5, 00185 Rome, Italy; fanelli@mat.uniroma1.it}

\affil[3]{Department of Mathematics, Faculty of Nuclear Sciences and Physical Engineering, Czech Technical University in Prague, Trojanova 13, 12000 Prague 2, Czechia; david.krejcirik@fjfi.cvut.cz}

\begin{document}

\date{\small 5 December 2019}

\maketitle

\begin{abstract}
\noindent
	By developing the method of multipliers, we establish sufficient conditions on the magnetic field and the complex, matrix-valued electric potential, which guarantee that the corresponding system of Schr\"odinger operators has no point spectrum. In particular, this allows us to prove analogous results for Pauli operators under the same electromagnetic conditions and, in turn, as a consequence of the supersymmetric structure, also for  magnetic Dirac operators.  
\end{abstract}

							%
							%

\section{Introduction}

\subsection{Objectives and state of the art}
Understanding electromagnetic phenomena has played 
a fundamental role in quantum mechanics.
The simplest mathematical model for the Hamiltonian of an electron,
subject to an external electric field described 
by a scalar potential $V:\R^3\to\R$
and an external magnetic field $B=\curl A$ 
with a vector potential $A\colon \R^3 \to \R^3$,
is given by the Schr\"odinger operator
\begin{equation}\label{Schrodinger}
  -\nabla_{\!A}^2 + V
  \qquad\mbox{in}\qquad
  L^2(\R^3;\C)
  \,,
\end{equation}
where $\nabla_{\!A} := \nabla + i A$ is the magnetic gradient.

Unfortunately, the mathematically elegant model~\eqref{Schrodinger}
is not sufficient to explain finer electromagnetic effects,
for it disregards an inner structure of electrons, namely their \emph{spin}.
A partially successful attempt to take the spin into account 
is to enrich the algebraic structure
of the Hilbert space and consider the Pauli operator
\begin{equation}\label{Pauli_perturbed}
  H_\textup{P}(A,\bm{V})
  :=-\nabla_{\!A}^2 \, \bm{I_{\C^2}} + \sigma \cdot B + \bm{V}
  \qquad\mbox{in}\qquad
  L^2(\R^3;\C^2)
  \,,
\end{equation}
where $\sigma := (\bm{\sigma_1},\bm{\sigma_2},\bm{\sigma_3})$ 
are Pauli matrices.
Here the term $\sigma \cdot B$ describes the interaction
of the spin with the magnetic field
and $\bm{V} := V \bm{I_{\C^2}}$
stands for the electric interaction as above. 

To get a more realistic description of the electron, 
subject to an external electromagnetic field,
one has to take relativistic effects into account.
A highly successful model is given by the Dirac operator
\begin{equation}\label{magnetic_Dirac}
  H_\textup{D}(A,\bm{V})
  :=-i  \alpha \cdot \nabla_{\!A} + \half \bm\beta + \bm{V}
  \qquad\mbox{in}\qquad
  L^2(\R^3;\C^4)
  \,,
\end{equation}
where $\alpha := (\bm{\alpha_1},\bm{\alpha_2},\bm{\alpha_3})$
and~$\beta$ are Dirac matrices and $\bm{V} := V \bm{I_{\C^4}}$.
 
The principal objective of this paper is to develop 
the so-called \emph{method of multipliers}
in order to establish spectral properties of the Pauli and Dirac operators. 
This technique comes from partial differential equations,
but it seems to be much less known in spectral theory.
We are primarily interested in physically relevant sufficient conditions, 
which guarantee the absence of point spectra
(including possibly embedded eigenvalues).
We proceed in greater generality by allowing
$V:\R^3 \to \C$ to be \emph{complex-valued} in~\eqref{Schrodinger}
and $\bm{V}:\R^3 \to \C^{2 \times 2}$ to be 
a general matrix-valued potential, possibly \emph{non-Hermitian},
in~\eqref{Pauli_perturbed}. 
However, some of our results are new even in the self-adjoint setting.
Since the spin-magnetic term $\sigma \cdot B$ can be included in~$\bm{V}$,
we simultaneously consider matrix 
electromagnetic Schr\"odinger operators
\begin{equation}\label{eq:Schro_int-unit}
  H_\textup{S}(A, \bm{V})
  :=- \nabla_{\!A}^2 \, \bm{I_{\C^2}} + \bm{V}
  \qquad\mbox{in}\qquad
  L^2(\R^3;\C^2)
  \,.
\end{equation}
Since the operator acts on spinors, 
we occasionally call the corresponding spectral problem 
the \emph{spinor Schr\"o\-dinger equation}.

As the last but not least generalisation to mention, 
in the main body of the paper,
we shall consider the Pauli and Dirac operators 
in the Euclidean space~$\R^d$ of \emph{arbitrary dimension} $d \geq 1$.

The study of spectral properties of 
scalar Schr\"odinger operators~\eqref{Schrodinger}
constitutes a traditional domain of mathematical physics 
and the literature on the subject is enormous.
Much less is known in the mathematically challenging 
and still physically relevant situations
where~$V$ is allowed to be complex-valued,
see \cite{FKV,FKV2} and references therein.
Works concerning non-self-adjoint Pauli operators 
are much more sparse in the literature, 
see~\cite{Sambou} and references therein.
More results are available in the case of non-self-adjoint Dirac operators,
see~\cite{Cuenin-Laptev-Tretter_2014,Dubuisson_2014,
Cuenin_2014,Sambou_2016,
Cuenin_2017,Enblom_2018,Cuenin-Siegl_2018,FK9}.

The paper~\cite{FKV} represents a first application of 
the method of multipliers to spectral theory:
the authors established sufficient conditions,
which guarantee the \emph{total} absence of eigenvalues of~\eqref{Schrodinger}.
It is remarkable that the conditions are physically relevant 
in the sense that they involve the magnetic field~$B$
rather than the vector potential~$A$.
The two-dimensional situation was covered later in~\cite{FKV2}.
The robustness of the method of multipliers 
has been demonstrated in its successful application 
to the half-space instead of the whole Euclidean space
in~\cite{Cossetti-Krejcirik_2018}
and to Lam\'e instead of Schr\"odinger operators 
in~\cite{Cossetti_2017}.
In the present paper,  
we push the analysis forward by investigating how 
the unconventional method provides meaningful and interesting results 
in the same direction also in the less explored setting 
of the spinorial Hamiltonians.

\subsection{The strategy}   
The main ingredient in our proofs is the method of multipliers 
as developed in~\cite{FKV} 
for scalar Schr\"odinger operators~\eqref{Schrodinger}. 
In the present paper, however, we carefully revisit the technique
and provide all the painful details,
which were missing in the previous works.
We identify various technical hypothesis about the electromagnetic 
potentials to justify the otherwise formal manipulations. 
We believe that this part of the paper will be of independent
interest for communities interested in spectral theory 
as well as partial differential equations. 

The next, completely new contribution is the adaptation of the method 
to the matrix electromagnetic Schr\"odinger operators~\eqref{eq:Schro_int-unit}.
The Pauli Hamiltonians~\eqref{Pauli_perturbed} 
are then covered as a particular case.

The method of multipliers does not seem to apply directly 
to Dirac operators, 
because of the lack of positivity of certain commutators.
Our strategy is to employ the supersymmetric structure
of Dirac operators (\emph{cf}.~\cite[Ch.~5]{Thaller}). 
More specifically, using the standard representation
\begin{equation}\label{std-repr3d}
	\bm{\alpha_\mu}=
	\begin{pmatrix}
		\bm{0} & \bm{\sigma_\mu}\\
		\bm{\sigma_\mu} & \bm{0}
	\end{pmatrix}, 
	\qquad 
	\bm{\beta}=
	\begin{pmatrix}
		\bm{I_{\C^2}} & \bm{0}\\
		\bm{0} & -\bm{I_{\C^2}}
	\end{pmatrix},
	\qquad \mu=1,2,3,
\end{equation}
and the commutation properties of the Pauli matrices,
it is easy to see that the square of 
the purely magnetic Dirac operator  
$H_\textup{D}(A,\bm{0}) =: H_\textup{D}(A)$
satisfies
\begin{equation}\label{favorable_form}
	H_\textup{D}(A)^2=
	\begin{pmatrix}
		H_\textup{P}(A) + \quarter \bm{I_{\C^2}} & \bm{0}\\
		\bm{0} & H_\textup{P}(A) + \quarter\bm{I_{\C^2}}
	\end{pmatrix},
\end{equation}
where $H_\textup{P}(A) := H_\textup{P}(A,\bm{0})$
is just the purely magnetic Pauli operator~\eqref{Pauli_perturbed}. 
This allows us to ensure the absence of the point spectrum of 
the Dirac operator $H_\textup{D}(A)$,
once the corresponding result for the Pauli operator $H_\textup{P}(A)$ 
is available, 
which, in turn, follows as a consequence of the corresponding result 
for the general Schr\"odinger operators $H_\textup{S}(A, \bm{V})$
with matrix-valued potentials~$\bm{V}$. 
Notice that, in this way, we are not able 
to treat magnetic Dirac operators
with electric perturbations.

\subsection{The results in three dimensions}   
As usual, the sums on the right-hand sides of~\eqref{Schrodinger},
\eqref{Pauli_perturbed} and~\eqref{eq:Schro_int-unit}
should be interpreted in a form sense (\emph{cf}.~\cite[Ch.~VI]{Kato}).
More specifically, the operators are introduced as the Friedrichs extension
of the operators initially defined on smooth functions of compact support. 
The regularity hypotheses and the functional inequalities
stated in the theorems below ensure that the operators
are well defined as m-sectorial operators.
The Dirac operator~\eqref{magnetic_Dirac} with $\bm{V}=\bm{0}$
is a closed symmetric operator under the stated assumptions.

Henceforth, we use the notation $r(x) := |x|$
for the distance function from the origin of~$\R^d$
and $\partial_r f(x) := \frac{x}{|x|}\cdot\nabla f(x)$
for the radial derivative of a function $f:\R^d\to\C$. 
We also set $f_\pm(x) := \max\{\pm f(x),0\}$ if~$f$ is real-valued.

For matrix Schr\"odinger operators~\eqref{eq:Schro_int-unit}, 
we prove the following result.
\begin{theorem}[Spinor Schr\"odinger equation]\label{thm:Schro_simplified}
Let $A\in L^2_\textup{loc}(\R^3; \R^3)$ 
be such that $B\in L^2_\textup{loc}(\R^3;\R^3).$ 
Suppose that $\bm{V}\in L^1_\textup{loc}(\R^3; \C^{2\times 2})$ 
admits the decomposition 
$\bm{V} = \bm{V^{(1)}} + \bm{V^{(2)}}$ with components 
$\bm{V^{(1)}}\in L^1_\textup{loc}(\R^3)$ 
and $\bm{V^{(2)}}=V^{(2)} \bm{I_{\C^2}}$,
where $V^{(2)} \in L^1_\textup{loc}(\R^3)$ 
is such that $[\partial_r(r \Re V^{(2)})]_+\in L^1_\textup{loc}(\R^3)$ 
and $r \bm{V^{(1)}}, r(\Re V^{(2)})_-, r\Im V^{(2)}\in L^2_\textup{loc}(\R^3).$ Assume that there exist numbers $a,b, \beta, \mathfrak{b}, c\in [0,1)$ satisfying
	\begin{equation}\label{Schro_cond_numbers_simplified}
		2(b + \beta + 2a)<1 
		\qquad \text{and} \qquad 
		2c + 2\beta + 6a + \mathfrak{b}^2 + \sqrt{2}(b + a)(\sqrt{\beta} + \sqrt{a})<1
	\end{equation}
	such that, for all two-vector $u$ with components in $C^\infty_0(\R^3),$ the inequalities
\begin{equation}\label{eq:conditions_B_V-3d}
\int_{\R^3} r^2 \abs{\bm{V^{(1)}}}^2 \abs{u}^2\ 
\leq a^2\int_{\R^3} \abs{\nabla_{\!A} u}^2 , 
			\qquad 
\int_{\R^3} r^2 \abs{B}^2 \abs{u}^2
\leq c^2\int_{\R^3} \abs{\nabla_{\!A} u}^2,
\end{equation}
and
\begin{gather}\label{Pauli-cond-V^2-3d}
		\int_{\R^3} r^2 (\Re V^{(2)})_-^2 \abs{u}^2
\leq b^2 \int_{\R^3} \abs{\nabla_{\!A} u}^2,
		\qquad 
		\int_{\R^3} r^2\abs{\Im V^{(2)}}^2 \abs{u}^2
\leq \beta^2 \int_{\R^3} \abs{\nabla_{\!A} u}^2,
		\\
		\label{Pauli-cond-radV^2-3d}
		\int_{\R^3} [\partial_r (r \Re V^{(2)})]_+ \abs{u}^2
\leq \mathfrak{b}^2 \int_{\R^3} \abs{\nabla_{\!A} u}^2,
\end{gather}	
	hold true. 
If in addition $A\in W^{1,3}_\textup{loc}(\R^3)$ 
and $V^{(2)}\in W^{1, 3/2}_\textup{loc}(\R^3),$
then $H_\textup{S}(A, \bm{V})$ has no eigenvalues, 
\emph{i.e.}\ $\sigma_\textup{p}(H_\textup{S}(A, \bm{V}))= \varnothing.$
\end{theorem}

As a consequence of the previous result, 
one has the corresponding theorem for Pauli operators.
\begin{theorem}[Pauli equation]\label{thm:main_result}
Under the hypotheses of Theorem~\ref{thm:Schro_simplified},
with~\eqref{Schro_cond_numbers_simplified} being replaced by
\begin{equation}\label{cond_numbers}
				2\big(b + \beta + 2(a + \sqrt{3}c)\big)<1 
		\qquad \text{and} \qquad 
		2c + 2\beta + 6(a+ \sqrt{3}c) + \mathfrak{b}^2 + \sqrt{2}\big(b + (a + \sqrt{3}c)\big)(\sqrt{\beta} + \sqrt{a + \sqrt{3}c})<1,					
\end{equation}
the operator $H_\textup{P}(A,\bm{V})$ has no eigenvalues, 
\emph{i.e.}\ $\sigma_\textup{p}(H_\textup{P}(A,\bm{V}))=\varnothing.$
\end{theorem}  

Due to the supersymmetric structure~\eqref{favorable_form} 
of the Dirac operator, the spectra of the Dirac and Pauli operators 
are intimately related. 
In particular, we deduce the following result from the previous theorem.
\begin{theorem}[Dirac equation]\label{thm:magnetic_Dirac3d}
	Let $A\in L^2_\textup{loc}(\R^3;\R^3)$ be such that $B\in L^2_\textup{loc}(\R^3; \R^3).$ Assume that there exists a number $c\in [0,1)$ satisfying
	\begin{equation}
	\label{cond_number_Dirac}
		4\sqrt{3}c<1
		\qquad \text{and} \qquad		
		2c + 6\sqrt{3}c + \sqrt{2}(\sqrt{3}c)^{3/2}<1
	\end{equation}
	such that, for all four-vector $u$ with components in $C^\infty_0(\R^3),$ the inequality
	\begin{equation}\label{cond_Dirac_B}
		\int_{\R^3} r^2\abs{B}^2\abs{u}^2 
\leq c^2 \int_{\R^3} \abs{\nabla_{\!A} u}^2 
	\end{equation}
	holds true.
If in addition $A\in W^{1,3}_\textup{loc}(\R^3),$
then $H_\textup{D}(A)$ has no eigenvalues, 
\emph{i.e.}\ $\sigma_\textup{p}(H_\textup{D}(A))=\varnothing.$
\end{theorem}
\begin{remark}

Notice that the conditions in~\eqref{cond_number_Dirac} are overabundant, in the sense that if $c$ is such that the second inequality 
of~\eqref{cond_number_Dirac} holds true, 
then $4\sqrt{3}c<1$ is automatically satisfied.
Indeed, the second inequality of~\eqref{cond_number_Dirac} 
requires $c < c_{1}^*$ where $c_1^*\approx 0.075,$ whereas the first requires $c<c_2^*$ where $c_2^*\approx 0.14.$
We decided to keep both conditions anyway 
in order to have a faster comparison 
with the corresponding results concerning the other theorems. 
\end{remark}  
%


\subsection{Organisation of the paper}
Even though so far we have considered only the three-dimensional framework, in this work we shall actually provide variants of the results presented above in \emph{any dimension}. (We anticipate already now that the two-dimensional framework will be excluded in the settings of Pauli and Dirac operators because of the well-known Aharonov--Casher effect.) In order to state our results in any dimension, however, an auxiliary material will be needed in order to introduce the general framework for the Pauli and Dirac Hamiltonians.
We therefore postpone the presentation of the general results
to Section~\ref{sec:main-higher},
while Section~\ref{Sec.any} is devoted to the definition of Dirac and Pauli operators
to any dimension (this section can be skipped by an experienced reader).
The method of multipliers for scalar Schr\"odinger operators
is revisited with all the necessary details in 
Section~\ref{Sec:scalar-Schro:1}.
The development of the method for Schr\"odinger operators
with matrix-valued potentials is performed in Section~\ref{Sec:Schro}.
The application of this general result to Pauli and Dirac operators
is given in Section~\ref{sec:Pauli-Dirac}.

\subsection{Notations}
Here we summarise specific notations and conventions 
that we use in this paper.

\begin{itemize}
\item We adopt the convention to write matrices in boldface.
\item For any dimension $d\geq 2,$ the physically relevant quantity associated to a given magnetic vector potential $A\colon \R^d\to \R^d$ is the $d\times d$ matrix-valued quantity
\begin{equation*}
	\bm{B} := (\nabla A) - (\nabla A)^t . 
\end{equation*}
Here, as usual, $(\nabla A)_{jk}= \partial_j A_k$ 
and $(\nabla A)^t_{jk}=(\nabla A)_{kj}$
with $j,k=1,2\dots, d$.
In $d=2$ and $d=3$ the magnetic tensor $\bm{B}$ 
can be identified with the scalar field 
$B_{12} = \partial_1 A_2 - \partial_2 A_1$ 
or the vector field $B=\curl A,$  respectively. 
More specifically, one has
\begin{equation*}
	\bm{B} w=
	\begin{system}
		&B_{12}\, w^\perp &\text{if}\quad d=2,\, \quad w\in \R^2\\ 
		&-B \times w &\text{if}\quad d=3,\, \quad w\in \R^3,
	\end{system}
\end{equation*}
where for any $w=(w_1,w_2)\in \R^2,$ $w^\perp := (w_2,-w_1)$ 
and the symbol $\times$ denotes the cross product in $\R^3.$

Notice that we did not comment on the case $d=1.$ In one dimension, in fact, the addition of a magnetic potential is trivial, in the sense that it is always possible to remove it by a suitable gauge transformation.
We refer to~\cite{Cazacu-Krejcirik_2016} for a complete survey on the concept of magnetic field in any dimensions and its definition in terms of differential forms and tensor fields.
\item  
We adopt the standard notation $|\cdot|$ for the Euclidean norm on $\C^d.$ 
We use the same symbol $|\cdot|$ for the operator norm: 
if $\bm{M}$ is a $d\times d$ matrix, we set
\begin{equation*}
\abs{\bm{M}}:=\sup_{\stackrel[v \not= 0]{}{v\in \C^d}} 
\frac{\abs{\bm{M} v}}{\abs{v}} .
\end{equation*}
\item
Let $v, w\in \R^d,$ the centered dot operation $v\cdot w$ designates the scalar product of the two vectors $v,w$ in $\R^d.$
\item
Given two vectors $v, w\in \R^d$ and a $d\times d$ matrix $\bm{M},$ the \emph{double}-centered dot operation $v \cdot \bm{M}\cdot w$ stands for the vector-matrix-vector product which returns the following scalar number
\begin{equation*}
	v\cdot \bm{M} \cdot w:=\sum_{j,k=1}^d v_k M_{k j} w_j.  
\end{equation*}
\item
We use the following definition for the $L^2$-norm of a vector-valued function 
$u=(u_1,u_2, \dots, u_n)$ on $\R^d$: 
\begin{equation*}
	\norm{u}_{[L^2(\R^d)]^n}:=\Bigg(\sum_{j=1}^n \norm{u_j}_{L^2(\R^d)}^2 \Bigg)^{1/2}.
\end{equation*}
\end{itemize}

\section{Definition of Dirac and Pauli Hamiltonians in any dimension}\label{Sec.any}
%
As already mentioned, our results will be stated in all dimensions $d\geq 1.$ In particular, this requires a more careful analysis on the Dirac and Pauli operators as their explicit form changes according to the underlying dimension.
Since here we are just interested in identifying the correct action of the operators,
we disregard issues with the operator domains for a moment.

\subsection{The Dirac operator}
%
Generalising the expression~\eqref{magnetic_Dirac} to arbitrary dimensions
requires ensuring existence of $d+1$ Hermitian matrices 
$\alpha:=(\bm{\alpha_1},\bm{\alpha_2}, \dots, \bm{\alpha_d})$ and $\bm{\beta}$ satisfying the anticommutation relations
\begin{equation}\label{eq:anticommutation}
\begin{gathered}
	\bm{\alpha_\mu} \bm{\alpha_\nu} + \bm{\alpha_\nu} \bm{\alpha_\mu}= 2\delta_{\mu\nu} \bm{I_{\C^{n(d)}}},\\
	\bm{\alpha_\mu} \bm{\beta} + \bm{\beta} \bm{\alpha_\mu}
  =\bm{0_{\C^{n(d)}}},\\
	\bm{\beta}^2=\bm{I_{\C^{n(d)}}},
\end{gathered}
\end{equation}
for $\mu, \nu\in \{1,2,\dots,d\}$,
where $\delta_{\mu \nu}$ represents the Kronecker symbol. 
The possibility to find such matrices clearly depends on the dimension $n(d)$ of the matrices themselves. In this regard one can verify that the following distinction is needed:
\begin{equation}\label{eq:n(d)}
 	n(d) :=
 	\begin{system}
 		&2^\frac{d+1}{2} &\text{if}\quad d \; \text{is odd},\\
 		&2^\frac{d}{2} &\text{if}\quad d \; \text{is even} . 
 	\end{system}
\end{equation}

Even though all that really cares are the anticommutation relations that the Dirac matrices satisfy, for the purpose of visualisation of the supersymmetric structure of the Dirac  operator, we shall rely on a particular representation of these matrices, that is the so-called \emph{standard representation}. According to the standard representation one defines the $d+1$ matrices $\alpha=(\bm{\alpha_1},\bm{\alpha_2}, \dots, \bm{\alpha_d})$ and $\bm{\beta}$ iteratively (with respect to the dimension) distinguishing between odd and even dimensions. For sake of clearness in the following the Dirac matrices are written with a superscript $^{(d)}$ to stress that these are constructed at the step corresponding to working in $d$ dimensions, 
\emph{e.g.}, $\alpha=(\bm{\alpha_1^{(d)}}, \bm{\alpha_2^{(d)}}, \dots, \bm{\alpha_d^{(d)}})$ and $\bm{\beta^{(d)}}$ are the $d+1$ Dirac matrices constructed in $d$ dimensions. Moreover, for notation convenience, we denote the matrix $\bm{\beta^{(d)}}$ as the $(d+1)$-th $\alpha$-matrix, namely $\bm{\beta^{(d)}}:=\bm{\alpha_{d+1}^{(d)}}.$

 \subsubsection*{Odd dimensions}
 If $d$ is odd, let us assume to know the $n(d-1) \times n(d-1)$ matrices $\bm{\alpha_1^{(d-1)}}, \bm{\alpha_2^{(d-1)}}, \dots, \bm{\alpha_{d}^{(d-1)}}$ corresponding to a previous step in the iteration. We then define $n(d)\times n(d)$ matrices (where, according to~\eqref{eq:n(d)}, $n(d)=2 n(d-1)$) in the following way:
 \begin{equation*}
 	\bm{\alpha_\mu^{(d)}}= 
 	\begin{pmatrix}
 		\bm{0} & \bm{\alpha_\mu^{(d-1)}}\vspace{0.2cm}\\
 		\bm{\alpha_\mu^{(d-1)}} & \bm{0}
 	\end{pmatrix},
 	\qquad
 	\bm{\beta^{(d)}}:=\bm{\alpha_{d+1}^{(d)}}=
 	\begin{pmatrix}
 		\bm{I_{\C^{n(d-1)}}} & \bm{0}\vspace{0.2cm}\\
 		\bm{0} & -\bm{I_{\C^{n(d-1)}}}
 	\end{pmatrix},
 	\quad \mu =1,2,\dots, d.	
 \end{equation*} 
 \subsubsection*{Even dimensions}
 If $d$ is even, we define $n(d)\times n(d)$ matrices (where, according to~\eqref{eq:n(d)}, $n(d)=n(d-1)=2n(d-2)$) as follows:
 \begin{equation*}
 	\bm{\alpha_1^{(d)}}=
 	\begin{pmatrix}
 		\bm{0} & \bm{I_{\C^{n(d-2)}}} \vspace{0.2cm}\\
 		\bm{I_{\C^{n(d-2)}}} & \bm{0}
	\end{pmatrix},
	\qquad
	\bm{\alpha_{\mu+1}^{(d)}}= 	 
 	\begin{pmatrix}
 		\bm{0} & -i\bm{\alpha_\mu^{(d-2)}}\vspace{0.2cm}\\
 		i\bm{\alpha_\mu^{(d-2)}} & \bm{0}
 	\end{pmatrix},
 	\quad \mu =1,2,\dots, d-1,
 \end{equation*}
 and 
 \begin{equation*}
 	\bm{\beta^{(d)}}:=\bm{\alpha_{d+1}^{(d)}}=
 	\begin{pmatrix}
 		\bm{I_{\C^{n(d-1)}}} & \bm{0}\vspace{0.2cm}\\
 		\bm{0} & -\bm{I_{\C^{n(d-1)}}}
 	\end{pmatrix}.
 \end{equation*} 

  Notice that we are also using the convention that $n(0)=1$ and that the $1\times 1$ matrix $\alpha_1^{(0)}=(1).$   This allows us to use the previous rule to construct the Dirac matrices corresponding to the standard representation also in $d=1$ and $d=2.$
  
  According to the construction above, one recognises that the Dirac matrices, regardless of the dimension, have all the following structure
 \begin{equation}\label{eq:higher-std-rep}
 \bm{\alpha_\mu}=
 	\begin{pmatrix}
 		\bm{0} & \bm{a_\mu^\ast}\vspace{0.2cm}\\
 		\bm{a_\mu} & \bm{0}
 	\end{pmatrix}, 
 	\qquad
 	\bm{\beta}=
 	\begin{pmatrix}
 		\bm{I_{\C^{n(d)/2}}} & \bm{0}\vspace{0.2cm}\\
 		\bm{0} & -\bm{I_{\C^{n(d)/2}}}
 	\end{pmatrix},
 	\quad 
 	\mu=1,2,\dots, d,
 \end{equation}
 where $\bm{a_\mu}$ are  $n(d)/2 \times n(d)/2$ matrices (Hermitian if $d$ is odd) such that
 \begin{equation}\label{a_anticomm}
 	\begin{gathered}
 		\bm{a_\mu} \bm{a_\nu^\ast} + \bm{a_\nu} \bm{a_\mu^\ast}=2\delta_{\mu\nu}\bm{I_{\C^{n(d)/2}}},\\
 		\bm{a_\mu^\ast} \bm{a_\nu} +\bm{a_\nu^\ast} \bm{a_\mu}=2\delta_{\mu\nu}\bm{I_{\C^{n(d)/2}}},
 	\end{gathered}
 \end{equation}
 for $\mu, \nu\in \{1,2,\dots, d\}.$ 
Here, as usual, $\bm{a_\mu^\ast}$ denotes the adjoint to $\bm{a_\mu},$ 
that is the conjugate transpose of $\bm{a_\mu}.$
We set $a:=(\bm{a_1},\dots,\bm{a_d})$.

 \begin{remark}\label{rmk:norm=1} 
 	Notice that, as a consequence of the fact that $\bm{\alpha_\mu}$ are Hermitian (in any dimension) and that $\bm{\alpha_\mu}^2=\bm{I_{\C^{n(d)}}},$ one has  $\abs{\bm{\alpha_\mu}}=1,$  $\mu=1,2,\dots, d.$  Therefore, due to the iterative construction above, one has that also the submatrices $\bm{a_\mu}$ and $\bm{a_\mu^\ast}$ have norm one, \emph{i.e.} $\abs{\bm{a_\mu}}=\abs{\bm{a_\mu^\ast}}=1.$ 
 \end{remark}

In the standard representation, that is using expression~\eqref{eq:higher-std-rep} for the Dirac matrices, the purely magnetic Dirac operator can be defined
through the following block-matrix differential expression
\begin{equation}\label{std_representation}
	H_\textup{D}(A) := 
	\begin{pmatrix}
		\half \bm{I_{\C^{n(d)/2}}} & \bm{D^\ast} \vspace{0.2cm}\\
		\bm{D} & -\half \bm{I_{\C^{{n(d)/2}}}}
	\end{pmatrix},
\end{equation}
where 
\begin{equation*}
	\bm{D} := -i a \cdot \nabla_{\!A},
	\qquad
	\bm{D^\ast} := -i a^\ast \cdot \nabla_{\!A}.
\end{equation*}
Notice that in odd dimension, being the submatrices $\bm{a_\mu}$ Hermitian, 
one has $\bm{D}=\bm{D^\ast}.$

\subsection{The square of the Dirac operator}
From representation~\eqref{std_representation}, 
it can be easily seen that $H_\textup{D}(A)$ can be decomposed as a sum of a $2\times 2$ \emph{diagonal} block and a $2\times 2$ \emph{off-diagonal} block operators. 
More specifically, one has
\begin{equation*}
	H_\textup{D}(A)=H_\textup{diag} + H_\textup{off-diag},
\end{equation*}
where
\begin{equation*}
	H_\textup{diag} :=
	\begin{pmatrix}
		\half  \bm{I_{\C^{n(d)/2}}} & \bm{0}\\
		\bm{0} & -\half \bm{I_{\C^{n(d)/2}}}
	\end{pmatrix},
	\qquad
	H_\textup{off-diag} :=
	\begin{pmatrix}
		\bm{0} & \bm{D^\ast}\\
		\bm{D} & \bm{0}
	\end{pmatrix}
  .
\end{equation*}
As one may readily check, $H_\textup{diag}$ and $H_\textup{off-diag}$ satisfy the anticommutation relation
\begin{equation}\label{supersymmetric_cond}
	H_\textup{diag} H_\textup{off-diag} + H_\textup{off-diag} H_\textup{diag}=\bm{0}. 
\end{equation} 
This distinguishing feature places the Dirac operator  within the class of operators with supersymmetry.  It is consequence of the supersymmetric condition~\eqref{supersymmetric_cond} that squaring out the Dirac operator gives
\begin{equation*}
	H_\textup{D}(A)^2=(H_\textup{diag} + H_\textup{off-diag})^2=H_\textup{diag}^2 + H_\textup{off-diag}^2,
\end{equation*}     
where
\begin{equation*}
	H_\textup{diag}^2=
	\begin{pmatrix}
		\quarter \bm{I_{\C^{n(d)/2}}} & \bm{0}\\
		\bm{0} & \quarter \bm{I_{\C^{n(d)/2}}}
	\end{pmatrix},
	\qquad
	H_\textup{off-diag}^2=					 
		\begin{pmatrix}
		 \bm{D^\ast}\bm{D} & \bm{0}\\
		\bm{0} & \bm{D}\bm{D^\ast}
		\end{pmatrix}.
\end{equation*}
Therefore, $H_\textup{D}(A)^2$ turns out to have the following favorable form
\begin{equation}\label{eq:fav-form-gen}
	H_\textup{D}(A)^2=
	\begin{pmatrix}
		\bm{D^\ast}\bm{D} + \quarter \bm{I_{\C^{n(d)/2}}} & \bm{0}\vspace{0.2cm}\\
		\bm{0} & \bm{D}\bm{D^\ast} + \quarter\bm{I_{\C^{n(d)/2}}}
	\end{pmatrix}.
\end{equation}
From property~\eqref{a_anticomm} of the Dirac submatrices, one can show that 
\begin{equation}\label{eq:DastD-DDast}
	\bm{D^\ast}\bm{D}=- \nabla_{\!A}^2 \bm{I_{\C^{n(d)/2}}} - \frac{i}{2}\, a^\ast\! \cdot \bm{B} \cdot a,
	\qquad
	\bm{D}\bm{D^\ast}=- \nabla_{\!A}^2 \bm{I_{\C^{n(d)/2}}} - \frac{i}{2}\, a\cdot \bm{B} \cdot a^\ast.
\end{equation}

\subsection{Low-dimensional illustrations}
In order to become more confident with the previous construction, we decided to present explicitly the situations of dimensions $d=1$ and $d=2$ in the next two subsections. 
(Dimension $d=3$ was already discussed above.)

  
\subsubsection{Dimension one}\label{rmk:1d}
In the Hilbert space $L^2(\R;\C^2)$,
the 1d Dirac operator reads 
\begin{equation*}
	H_\textup{D}(0):=-i \bm{\alpha} \nabla + \half\bm{\beta} ,
\end{equation*}
where~$\nabla$ is just a weird notation for an ordinary derivative.  
With the notation $H_\textup{D}(0)$ we emphasise that the magnetic potential~$A$ has been chosen to be identically equal to zero, since in one dimension it can be always removed by choosing a suitable gauge.
One can immediately verify that squaring out the operator $H_\textup{D}(0)$ yields
\begin{equation*}
	{H_\textup{D}(0)}^2
  =- \nabla^2\bm{I_{\C^2}} + \quarter \bm{I_{\C^2}} .
\end{equation*}
%
%
According to the rule provided above, in the standard representation, one chooses
$ \bm{\alpha}:= \bm{\sigma_1}$ and $\bm{\beta}:= \bm{\sigma_3},$
where $\bm{\sigma_1}$ and $\bm{\sigma_3}$ are two of the three Pauli matrices. 
Thus, one conveniently writes
\begin{equation*}
	H_\textup{D}(0)=
	\begin{pmatrix}
		\half & D\vspace{0.2cm}\\
		D  &-\half
	\end{pmatrix},
\end{equation*}
where $D :=-i \nabla$
and
\begin{equation}\label{Dirac.1D}
	{H_\textup{D}(0)}^2=
	\begin{pmatrix}
		H_\textup{P}(0) + \quarter& 0\vspace{0.2cm}\\
		0  & H_\textup{P}(0) + \quarter
	\end{pmatrix},
\end{equation}
with the Pauli operator
\begin{equation}\label{Pauli1d}
	H_\textup{P}(0):=- \nabla^2 .
\end{equation} 
Hence, in one dimension, the Pauli operator coincides with 
the free one dimensional Schr\"odinger operator acting in $L^2(\R;\R)$.
%

\subsubsection{Dimension two}\label{rmk:Pauli-2d}
%
In the Hilbert space $L^2(\R^2;\C^2)$,
the 2d Dirac operator reads 
 \begin{equation*}
 	H_\textup{D}(A):=-i \alpha \cdot\nabla_{\!A} + \half\bm{\beta},
 \end{equation*}
 where $\alpha:=(\bm{\alpha_1}, \bm{\alpha_2})$ and $\bm{\beta}$ are $2\times 2$ Hermitian matrices satisfying~\eqref{eq:anticommutation}. Squaring out $H_\textup{D}(A)$ yields
 \begin{equation*}
 	H_\textup{D}(A)^2
  =-\nabla_{\!A}^2 \bm{I_{\C^2}} -\frac{i}{2}[\bm{\alpha_1},\bm{\alpha_2}] B_{12} 
  + \quarter \bm{I_{\C^2}}.
 \end{equation*}
 According to the rule provided above, in the standard representation,
 one chooses $\bm{\alpha_1}:=\bm{\sigma_1},$ $\bm{\alpha_2}:=\bm{\sigma_2}$ and $\bm{\beta}:=\bm{\sigma_3}$. This gives $[\bm{\alpha_1},\bm{\alpha_2}]=2i\bm{\sigma_3}$ and
 \begin{equation*}
 H_\textup{D}(A)=
 	\begin{pmatrix}
 		\half & D^\ast\\
 		D &-\half
 	\end{pmatrix},
 \end{equation*}
 where 
\begin{equation*} 
 D :=-i \partial_{1, A} + \partial_{2, A}
 , \qquad
  D^\ast :=-i \partial_{1, A} - \partial_{2, A},
 \end{equation*}
 and $\partial_{j,A}:= \partial_j + i A_j,$ $j=1,2.$
 Thus
 \begin{equation*}
 	 {H_\textup{D}(A)}^2
         = H_\textup{P}(A) + \quarter\bm{I_{\C^{2}}}
 \end{equation*}
with the Pauli operator
$$
  H_\textup{P}(A) := -\nabla_{\!A}^2 \bm{I_{\C^2}} +  \bm{\sigma_3} B_{12}
$$
%

\subsection{The Pauli operator}
After these illustrations,
let us come back to the general dimension $d \geq 1$.
Recall that the Dirac operator~$H_\textup{D}(A)$ 
has been introduced via~\eqref{std_representation}
and that its square satisfies~\eqref{eq:fav-form-gen}.
The following lemma specifies the form of the square 
according to the parity of the dimension
and offers a natural definition for the Pauli operator in any dimension.
 
\begin{lemma}[Algebraic definition of Pauli operators]\label{lemma:explicit-def}
	Let $d\geq 1$ and let $n(d)$ be as in~\eqref{eq:n(d)}.
	\begin{itemize}
		\item
			If $d$ is odd, then
			\begin{equation}\label{eq:oddDirac}
 				H_\textup{D}^\textup{odd}(A)^2=
				\begin{pmatrix}
					H_\textup{P}^\textup{odd}(A) + \quarter 					\bm{I_{\C^{n(d)/2}}} & \bm{0}\vspace{0.2cm}\\
					\bm{0} & H_\textup{P}^\textup{odd}(A) + \quarter \bm{I_{\C^{n(d)/2}}}
				\end{pmatrix},
			\end{equation} 
			where we define
			\begin{equation}\label{Pauli-odd}
				H_\textup{P}^\textup{odd}(A):= -\nabla_{\!A}^2 \bm{I_{\C^{n(d)/2}}} - \frac{i}{2}\, a \cdot \bm{B} \cdot a	.
			\end{equation}
		\item
			If $d$ is even, then
			\begin{equation}\label{eq:evenDirac}
 				H_\textup{D}^\textup{even}(A)^2=H_\textup{P}^\textup{even}(A) + \quarter \bm{I_{\C^{n(d)}}},
			\end{equation} 
			where we define
			\begin{equation}\label{Pauli-even}
				H_\textup{P}^\textup{even}(A):=
				\begin{pmatrix}
					-\nabla_{\!A}^2 \bm{I_{\C^{n(d)/2}}} - \frac{i}{2}\, a^\ast\! \cdot \bm{B} \cdot a,
					& \bm{0}\\
					\bm{0} &
					-\nabla_{\!A}^2 \bm{I_{\C^{n(d)/2}}} - \frac{i}{2}\, a\cdot \bm{B} \cdot a^\ast
				\end{pmatrix}.	
			\end{equation}	
		\end{itemize}
		\end{lemma} 
		\begin{proof}
			In odd dimensions one has that $\bm{D^\ast}=\bm{D},$ therefore
\begin{equation*}
	\bm{D^\ast}\bm{D}=\bm{D}\bm{D^\ast}=\bm{D}^2=\big[-ia\cdot \nabla_{\!A}\big]^2=-\nabla_{\!A}^2 \bm{I_{\C^{n(d)/2}}} - \frac{i}{2}\, a \cdot \bm{B} \cdot a.
\end{equation*}
Thus, defining 
\begin{equation*}
	H_\textup{P}^\textup{odd}(A):= \bm{D^\ast}\bm{D}
\end{equation*}
and using~\eqref{eq:fav-form-gen} one immediately gets the desired representation in odd dimensions. 
In even dimensions one defines
\begin{equation*}
	H_\textup{P}^\textup{even}(A):=
	\begin{pmatrix}
		 \bm{D^\ast D} & \bm{0}\\
		 \bm{0} & \bm{D^\ast D}
	\end{pmatrix}.	
\end{equation*}
Hence, from~\eqref{eq:fav-form-gen} and~\eqref{eq:DastD-DDast} one readily has the thesis.
\end{proof}

Notice that in even dimensions the Pauli operator is a matrix 
operator with the same dimension as the Dirac Hamiltonian.
In odd dimensions the dimension of the Pauli operator 
is a half of that of the Dirac operator.
Recalling~\eqref{eq:n(d)}, we therefore set 
\begin{equation}\label{eq:n'(d)}
 	n'(d) :=
 	\begin{system}
 		&n(d)/2 &\text{if}\quad d \; \text{is odd},\\
 		&n(d) &\text{if}\quad d \; \text{is even} . 
 	\end{system}
\end{equation}

\subsection{Domains of the operators}
Finally, we specify the domains of the Dirac and Pauli operators.
Notice that the rather formal manipulations of the preceding subsections
can be justified when the action of the operators is considered 
on smooth functions of compact support.
Therefore, we shall define each of the operators as an extension 
of the operator initially defined on such a restricted domain.
We always assume that the vector potential  
$A \in L_\mathrm{loc}^2(\R^d;\R^d)$ is such that 
$\bm{B} \in L_\mathrm{loc}^1(\R^d;\R^{d \times d})$.

We define the Pauli operator $H_\textup{P}(A)$ acting on 
the Hilbert space $L^2(\R^d;\R^{n'(d)})$
as the self-adjoint Friedrichs extension of the operator initially
considered on the domain $C_0^\infty(\R^d;\R^{n'(d)})$;
notice that this initial operator is symmetric.
Disregarding the spin-magnetic term for a moment,
the form domain can be identified with the magnetic Sobolev space
(\emph{cf}.~\cite[Sec.~7.20]{LL})
\begin{equation}\label{Sobolev.magnet}
  H_{\!A}^1(\R^d;\R^{n'(d)}) :=
  \left\{u\in L^2(\R^d;\R^{n'(d)})\colon \
  \partial_{j,A} u\in L^2(\R^d;\R^{n'(d)})
  \mbox{ for every } j\in\{1,\dots,d\} \right\} 
  .
\end{equation}
The operator domain is the subset of $H_{\!A}^1(\R^d;\R^{n'(d)})$
consisting of functions~$\psi$ such that 
$\nabla_{\!A}^2 \psi \in L^2(\R^d;\R^{n'(d)})$. 
To include the spin-magnetic term, 
we make the hypothesis that there exist numbers $a<1$ and $b \in \R$
such that, for every $\psi \in C_0^\infty(\R^d)$, 
\begin{equation}\label{rel.bounded}
  \frac{1}{2} \int_{\R^d} |\bm{B}| |\psi|^2 
  \leq a \int_{\R^d} |\nabla_{\!A}\psi|^2 + b \int_{\R^d} |\psi|^2
  \,.
\end{equation}
Then the spin-magnetic term is a relatively form-bounded perturbation
of the already defined operator with the relative bound less than one
(recall Remark~\ref{rmk:norm=1}),
so the Pauli operator $H_\textup{P}(A)$ 
with the same form domain~\eqref{Sobolev.magnet} 
is indeed self-adjoint.

For the domain of the Dirac operator~\eqref{std_representation}
we take
\begin{equation}\label{domain.Dirac}
  \Dom(H_\textup{D}(A)) := H_{\!A}^1(\R^d;\R^{n(d)}) \,.
\end{equation}
Notice that $H_\textup{D}(A)$ is symmetric.
Using Lemma~\ref{lemma:explicit-def},
for every $\psi \in C_0^\infty(\R^d;\R^{n(d)})$,
which is dense in $\Dom(H_\textup{D}(A))$,
we have the identity (with a slight abuse of notation)
$$
  \|H_\textup{D}(A)\psi\|^2
  = (\psi,H_\textup{D}(A)^2\psi)
  = (\psi,H_\textup{P}(A)\psi) + \quarter \|\psi\|^2
  \,.
$$
Since the quadratic form of the Pauli operator $H_\textup{P}(A)$ 
is closed on the space~\eqref{Sobolev.magnet}, it follows 
that the Dirac operator $H_\textup{D}(A)$ with~\eqref{domain.Dirac}  
is a closed symmetric operator.
Under further assumptions about the vector potential
(see \cite[Sec.~4.3]{Thaller}), 
one can ensure that $H_\textup{D}(A)$ is actually self-adjoint,
but our results hold under the present more general setting.

\section{Statement of the main results in any dimension}\label{sec:main-higher}
%
Now we are in position to state our main results in any dimension. As anticipated, in order to do that, we shall consider separately the three spinorial Hamiltonians. 

\subsection{The spinor Schr\"odinger equation}\label{ss:Schro}

Let us start by considering the matrix Schr\"odinger operator
\begin{equation}\label{Schrodinger.any}
  H_\textup{S}(A, \bm{V})
  :=- \nabla_{\!A}^2 \, \bm{I_{\C^n}} + \bm{V}
  \qquad\mbox{in}\qquad
  L^2(\R^d;\C^n)
  \,,
\end{equation}
which is an extension of~\eqref{eq:Schro_int-unit}
to any dimension $d \geq 1$ and $n \geq 1$.
Here $\bm{V} \in L_\mathrm{loc}^1(\R^d;\C^{n \times n})$
and $A \in L_\mathrm{loc}^2(\R^d;\R^d)$.
The operator is properly introduced as the Friedrichs extension
of the operator initially defined on $C_0^\infty(\R^d;\C^n)$.
The hypotheses in the theorems below ensure that 
$H_\textup{S}(A, \bm{V})$ is well defined as an m-sectorial operator.

\subsubsection{A general result in any dimension}

\begin{theorem}\label{thm:Schro}
Given any $d, n\geq 1$, 
let $A\in L^2_\textup{loc}(\R^d;\R^d)$ be such that $\bm{B}\in L^2_\textup{loc}(\R^d;\R^{d\times d}).$ Suppose that $\bm{V}\in L^1_\textup{loc}(\R^d; \C^{n\times n})$ admits the decomposition $\bm{V}= \bm{V^{(1)}} + \bm{V^{(2)}}$ with components $\bm{V^{(1)}}\in L^1_\textup{loc}(\R^d)$ and $\bm{V^{(2)}}=V^{(2)}\bm{I_{\C^n}}$,
where $V^{(2)}\in L^1_\textup{loc}(\R^d)$ is 
such that $[\partial_r(r \Re V^{(2)})]_+\in L^1_\textup{loc}(\R^d)$ 
and $r \bm{V^{(1)}}, r(\Re V^{(2)})_-, r\Im V^{(2)}\in L^2_\textup{loc}(\R^d).$ Assume that there exist numbers $a_1, a_2, b_1, b_2, \mathfrak{b}, \beta_1, \beta_2, c \in [0,1)$ satisfying
	\begin{equation}\label{Schro_cond_numbers}
		b_1^2 + \beta_1^2 + 2a_1^2<1 
		\qquad \text{and} \qquad 
		2c + 2\beta_2 + 2a_2 + (d-1)a_1^2 + \mathfrak{b}^2 + (b_2 + a_2)(\beta_1 + a_1)<1
	\end{equation}
	such that, for all $n$-vector $u$ with components in $C^\infty_0(\R^d),$
	\begin{equation}\label{hyp:V^1}
		\int_{\R^d} \abs{\bm{V^{(1)}}} \abs{u}^2 
\leq a_1^2 \int_{\R^d} \abs{\nabla_{\!A} u}^2 , 
		\qquad
		\int_{\R^d} r^2\abs{\bm{V^{(1)}}}^2 \abs{u}^2
\leq a_2^2 \int_{\R^d} \abs{\nabla_{\!A} u}^2 ,
	\end{equation}
	\begin{gather}
		\label{hyp:ReV^2-}
		\int_{\R^d} (\Re V^{(2)})_- \abs{u}^2
\leq b_1^2 \int_{\R^d} \abs{\nabla_{\!A} u}^2 , 
		\qquad
		\int_{\R^d} r^2 (\Re V^{(2)})_-^2 \abs{u}^2 
\leq b_2^2 \int_{\R^d} \abs{\nabla_{\!A} u}^2 ,\\
		\label{hyp:radial_der-rReV^2}
		\int_{\R^d} [\partial_r (r \Re V^{(2)})]_+ \abs{u}^2 
\leq \mathfrak{b}^2 \int_{\R^d} \abs{\nabla_{\!A} u}^2 ,
		\end{gather}
		\begin{equation}\label{hyp:ImV^2}
		\int_{\R^d} \abs{\Im V^{(2)}} \abs{u}^2 
\leq \beta_1^2 \int_{\R^d} \abs{\nabla_{\!A} u}^2 , 
		\qquad
		\int_{\R^d} r^2\abs{\Im V^{(2)}}^2 \abs{u}^2 
\leq \beta_2^2 \int_{\R^d} \abs{\nabla_{\!A} u}^2 ,
	\end{equation}
	\begin{equation}\label{hyp:B}
		\int_{\R^d} r^2 \abs{\bm{B}}^2\abs{u}^2
\leq c^2 \int_{\R^d} \abs{\nabla_{\!A} u}^2 .
	\end{equation}
If $d=2$ assume also that the inequality
	\begin{equation}\label{positivity-d=2}
		\frac{1}{2} \int_{\R^2} \frac{~\abs{u}^2}{r}  
\leq \int_{\R^2} r\abs{\nabla_{\!A} u}^2  
+ \int_{\R^2} r (\Re V^{(2)})_+ \abs{u}^2 
	\end{equation}
	holds true.
	If, in addition, one has
	\begin{equation}\label{additional}
		A\in W^{1,2p}_\textup{loc}(\R^d)\quad \text{and}\quad \Re V^{(2)}\in W^{1,p}_\textup{loc}(\R^d) \,,
		\qquad \text{where}\qquad
		\begin{system}
		&p=1  &\text{if}\quad d=1,\\
		&p>1  &\text{if}\quad d=2,\\
		&p=d/2  &\text{if}\quad d\geq 3,
		\end{system}
	\end{equation}
	then $H_\textup{S}(A,\bm{V})$ has no eigenvalues, 
\emph{i.e.}\ $\sigma_\textup{p}(H_\textup{S}(A,\bm{V}))=\varnothing.$ 
\end{theorem}

The theorem is commented on in the following subsections.

\subsubsection{Criticality of low dimensions}\label{Sec.criticality}
Because of the criticality of the Laplacian in $L^2(\R^d)$ with $d=1,2$,
the lower dimensional scenarios are a bit special.

First of all, 
due to the absence of magnetic phenomena in $\R^1,$ the corresponding assumptions \eqref{hyp:V^1}--\eqref{hyp:B} in dimension $d=1$ come with the classical gradient $\nabla$ as a replacement of the magnetic gradient~$\nabla_{\!A}.$ 
Consequently, because of the criticality of the Laplacian in $L^2(\R)$,
necessarily $\bm{V^{(1)}}=0$, $(\Re V^{(2)})_-=0$,
$[\partial_r (r \Re V^{(2)})]_+=0$ and $\Im V^{(2)}=0$.
Moreover, \eqref{hyp:B}~is always satisfied if $d=1$ being $\bm{B}$ equal to zero. 
Hence, if $d=1$, the theorem essentially says that 
the scalar Schr\"odinger operator $-\nabla^2+V$ in $L^2(\R)$
has no eigenvalues, provided that~$V$ is non-negative
and the radial derivative $\partial_r (r V)$ is non-positive.
The requirements respectively exclude non-positive and positive eigenvalues.
The latter is a sort of the classical repulsiveness 
requirement (\emph{cf.}~\cite[Thm.~XIII.58]{RS4}).

Similarly, if $d=2$ and there is no magnetic field 
(\emph{i.e.}\ $\bm{B}=\bm{0}$),
the theorem essentially says that 
the scalar Schr\"odinger operator $-\nabla^2+V$ in $L^2(\R^2)$
has no eigenvalues, provided that~$V$ is non-negative
and the radial derivative $\partial_r (r V)$ is non-positive
(again, the conditions exclude non-positive 
and positive eigenvalues, respectively).
On the other hand, in two dimensions,
the situation becomes interesting if the magnetic field is present.
Indeed, the magnetic Laplacian in $L^2(\R^2)$ is subcritical
due to the existence of magnetic Hardy inequalities
(see~\cite{Laptev-Weidl} for the pioneering work
and~\cite{Cazacu-Krejcirik_2016} for the most recent developments).
The latter guarantee a source of sufficient conditions 
to make the hypotheses \eqref{hyp:V^1}--\eqref{hyp:B} non-trivial
(\emph{cf.}~\cite{FKV2}).

\subsubsection{An alternative statement in dimension two}

	We want to comment more on the additional condition~\eqref{positivity-d=2} in dimension $d=2.$ 
	Using the 2d weighted Hardy inequality
	\begin{equation}\label{eq:weighted_Hardy-2d}
		\int_{\R^2} r \, \abs{\nabla_{\!A} u}^2  
\geq \frac{1}{4} \int_{\R^2} \frac{~\abs{u}^2}{r} ,
	\end{equation}  
	 it is easy to check that requiring ``enough'' positivity to $\Re V^{(2)}$ will guarantee the validity of~\eqref{positivity-d=2}. 
More specifically, the pointwise bound
	\begin{equation*}
		[\Re V^{(2)}(x)]_+\geq \frac{1}{4\abs{x}^2},
	\end{equation*}
valid for almost every $x \in \R^2$
	is sufficient for~\eqref{positivity-d=2} to hold.
	On the other hand, without the positivity of $\Re V^{(2)}$,
 condition~\eqref{positivity-d=2} is quite restrictive. Indeed, if one assumes 
$V^{(2)}=0,$ then ensuring the validity of~\eqref{positivity-d=2}, would require to ensure the existence of vector potentials $A$ for which an improvement of the weighted Hardy inequality~\eqref{eq:weighted_Hardy-2d} holds true
(for \eqref{positivity-d=2} with $V^{(2)}=0$ is nothing but~\eqref{eq:weighted_Hardy-2d} with a better constant). 

	For this reason, following an idea introduced in~\cite[Sec.~3.2]{FKV2}, 
we provide an alternative result, which avoids condition~\eqref{positivity-d=2},
but a stronger hypothesis compared to~\eqref{Schro_cond_numbers} is assumed.
 	
	\begin{theorem} 
\label{thm:alternative2d}
		Let $d=2$ and let $n, A, \bm{B}$ and $\bm{V}$ be as in Theorem~\ref{thm:Schro}. 		
		Assume that there exist numbers $a_1, a_2, b_1, b_2, \mathfrak{b}, \beta_1, \beta_2, c , \epsilon \in [0,1)$ satisfying
	\begin{equation}\label{eq:alt-Schr-cond-2d}
	b_1^2 + \beta_1^2 + 2a_1^2<1
	\qquad \text{and} \qquad 		
	2c + 2\beta_2 + 2a_2 + a_1^2 + \mathfrak{b}^2 + (b_2 + a_2)(\beta_1 + a_1) + 4\epsilon + \frac{17}{\epsilon}(a_1^2 + \beta_1^2)<1,
	\end{equation}
	such that, for all $n$-vector $u$ with components in $C^\infty_0(\R^2),$ inequalities \eqref{hyp:V^1}--\eqref{hyp:B} hold true. If, in addition, one has
	\begin{equation*}
		A \in W^{1,2p}_\textup{loc}(\R^2)
		\qquad \text{and} \qquad
		\Re V^{(2)}\in W^{1,p}_\textup{loc}(\R^2),
		\qquad \mbox{where} \qquad p>1,
	\end{equation*}
	then $H_\textup{S}(A, \bm{V})$ has no eigenvalues, 
\emph{i.e.}\ $\sigma_\textup{p}(H_\textup{S}(A, \bm{V}))=\varnothing.$
	\end{theorem}

\subsubsection{A simplification in higher dimensions}\label{rmk:overabundant}
In dimensions $d\geq 3,$ as a consequence of the diamagnetic inequality (see~\cite{Kato_1972} and~\cite[Thm.~7.21]{LL}) 
\begin{equation}\label{eq:diamagnetic}
	\abs{\nabla \abs{\psi}(x)}\leq \abs{\nabla_{\!A} \psi(x)} \qquad \text{a.e.}\quad x\in \R^d,
\end{equation}    
together with the \emph{classical} Hardy inequality 
	\begin{equation}\label{eq:classical_Hardy}
		 \int_{\R^d} \frac{\abs{\psi(x)}^2}{\abs{x}^2}\, dx
		\leq \frac{4}{(d-2)^2}\int_{\R^d} \abs{\nabla \psi}^2\, dx, \qquad \forall\, \psi \in C^\infty_0(\R^d), \qquad d\geq 3,
	\end{equation}
	applied to $\abs{\psi},$ one can prove the following \emph{magnetic} Hardy inequality
	\begin{equation}\label{H-general}
		\int_{\R^d} \frac{\abs{\psi(x)}^2}{\abs{x}^2}\, dx\leq \frac{4}{(d-2)^2} \int_{\R^d} \abs{\nabla_{\!A} \psi}^2\, dx,
  \qquad \forall\, \psi \in C^\infty_0(\R^d), \qquad d\geq 3.
	\end{equation}
	Using~\eqref{H-general}, it is easy to check that the first inequalities in~\eqref{hyp:V^1},~\eqref{hyp:ReV^2-} and~\eqref{hyp:ImV^2} follow respectively as a consequence of the second inequalities in~\eqref{hyp:V^1},~\eqref{hyp:ReV^2-} and~\eqref{hyp:ImV^2} with
\begin{equation*}
	a_1^2:=\frac{2}{d-2} a_2, \quad b_1^2:=\frac{2}{d-2} b_2, \quad \beta_1^2:=\frac{2}{d-2}\beta_2,
\end{equation*}
and assuming $a_2, b_2, \beta_2<(d-2)/2.$ 
Hence, in the higher dimensions $d\geq 3,$ conditions in~\eqref{Schro_cond_numbers} simplifies to
\begin{equation}\label{Schro_cond_numbers-H}
	\frac{2}{d-2}\Big(b_2 + \beta_2 + 2 a_2\Big)<1
	\qquad \text{and} \qquad
	2c + 2\beta_2 + \frac{2(2d-3)}{d-2} a_2 + \mathfrak{b}^2 + \frac{\sqrt{2}}{\sqrt{d-2}} (b_2 + a_2)(\sqrt{\beta_2} +\sqrt{a_2})<1.
\end{equation}
In particular, this justifies the fact that in Theorem~\ref{thm:Schro_simplified} which is a special case of Theorem~\ref{thm:Schro} for $d=3$ (and $n=2$) we assume \emph{only} the validity of~\eqref{eq:conditions_B_V-3d},~\eqref{Pauli-cond-V^2-3d} and~\eqref{Pauli-cond-radV^2-3d}, moreover~\eqref{Schro_cond_numbers} is replaced by~\eqref{Schro_cond_numbers_simplified} (notice that dropping the subscript $\cdot_2$ in the constants and fixing $d=3$ in~\eqref{Schro_cond_numbers-H} gives~\eqref{Schro_cond_numbers_simplified}).

\subsubsection{The Aharonov--Bohm field}
Let us come back to dimension two 
and consider the \emph{Aharonov--Bohm} magnetic potential
	\begin{equation}\label{eq:A-B}
		A(x,y):= (-\sin \theta, \cos\theta) \, \frac{\alpha(\theta)}{r},
	\end{equation}
	where $(x,y)=(r\cos\theta, r\sin\theta)$ is the parametrisation via polar coordinates, $r\in (0,\infty),$ $\theta\in [0,2\pi),$ and $\alpha\colon [0,2\pi)\to \R$ is an arbitrary bounded function. In this specific case,
there is an explicit magnetic Hardy-type inequality 
(see \cite[Thm.~3]{Laptev-Weidl})
	\begin{equation}\label{eq:magnetic_Hardy}
	\int_{\R^2}\abs{\nabla_{\!A} \psi}^2 
\geq \gamma^ 2  \int_{\R^2} \frac{\abs{\psi}^2}{r^2} ,
\qquad \forall\, \psi\in C_0^\infty(\R^2 \setminus \{0\})
, \qquad 
   \gamma:=\dist\{\bar{\alpha}, \mathbb{Z}\},
	\end{equation}
 where $\bar{\alpha}$ has the physical meaning of the total magnetic flux:
	\begin{equation}\label{eq:mag-flux}
		\bar{\alpha}:= \frac{1}{2\pi} \int_0^{2\pi} \alpha(\theta)\, d\theta.
	\end{equation}
Notice that in this case the magnetic field~$B$ equals zero everywhere except for $x=0$; indeed 
	\begin{equation}\label{eq:B_A-B}
	B=2\pi \bar{\alpha} \delta
	\end{equation}
	in the sense of distribution, where $\delta$ is the Dirac delta function.

	The Aharonov--Bohm potential~\eqref{eq:A-B} is \emph{not} in $L^2_\textup{loc}(\R^2),$ 
so the matrix Schr\"odinger operator is not well defined 
as described below~\eqref{Schrodinger.any} 
and Theorem~\ref{thm:Schro} does not apply to it as such. 
Now the Schr\"odinger operator $H_\textup{S}(A, \bm{V})$
is introduced as the Friedrichs extension of 
the operator~\eqref{eq:Schro_int-unit}
initially defined on $C_0^\infty(\R^2 \setminus \{0\};\C^n)$.
At the same time, it is possible to adapt the method of multipliers 
in such a way that it covers this situation as well.
The following result can be considered as 
an extension of \cite[Thm.~5]{FKV2} in the scalar case
to the spinorial Schr\"odinger equation.
	
	\begin{theorem}\label{thm:A-B}		
		Let $d=2$ and let $A$ be as in~\eqref{eq:A-B} with $\bar{\alpha}\notin \mathbb{Z}$ and $\bm{V}$ as in Theorem~\ref{thm:Schro}.
		Assume that there exist numbers $a, b, \mathfrak{b}, \beta, \epsilon\in [0,1)$ satisfying
		\begin{equation}\label{eq:alt-Schr-cond-AB}
			\frac{1}{\gamma} (b+ \beta + 2a)<1			
			\quad \text{and} \quad
			2\beta + 2a + \frac{a}{\gamma} + \mathfrak{b}^2 + \frac{1}{\sqrt{\gamma}}(b + a)(\sqrt{a} + \sqrt{\beta} )+ \left( \frac{1}{4} - \gamma^2 \right) \left[\epsilon + \frac{(a+\beta)}{\epsilon \gamma^3}\right]<1 ,
		\end{equation}
		with $\gamma:=\dist\{\bar{\alpha}, \mathbb{Z}\},$
		such that, for all $n$-vector $u$ with component in $C^\infty_0(\R^2\setminus \{0\}),$ inequalities
		\begin{equation}\label{eq:V^1-AB}
	\int_{\R^2} r^2 \abs{\bm{V^{(1)}}}^2 \abs{u}^2 
\leq a^2\int_{\R^2} \abs{\nabla_A u}^2 , 
		\end{equation}
		and
		\begin{gather}
		\label{eq:V^2-AB}
		\int_{\R^2} r^2 (\Re V^{(2)})_-^2 \abs{u}^2 
\leq b^2 \int_{\R^2} \abs{\nabla_A u}^2 ,
		\qquad 
		\int_{\R^2} r^2\abs{\Im V^{(2)}}^2 \abs{u}^2 
\leq \beta^2 \int_{\R^2} \abs{\nabla_A u}^2 ,
		\\
		\int_{\R^2} [\partial_r (r \Re V^{(2)})]_+ \abs{u}^2 
\leq \mathfrak{b}^2 \int_{\R^2} \abs{\nabla_A u}^2 
		\end{gather}
		hold true. If, in addition, one has
\begin{equation*}
	\Re V^{(2)}\in W^{1,p}_\textup{loc}(\R^2), \qquad p>1,
\end{equation*}		
		then $H_\textup{S}(A, \bm{V})$ has no eigenvalues, 
\emph{i.e.}\ $\sigma_\textup{p}(H_\textup{S}(A, \bm{V}))=\varnothing.$
	\end{theorem}

\subsubsection{On the regularity condition
\texorpdfstring{\eqref{additional}}{additional}
and their replacement}
	As we will see in more details later on (see Section~\ref{Further hypotheses A,V}), the additional local regularity assumptions~\eqref{additional} on the potentials are needed in order to justify \emph{rigorously} the algebraic manipulations that the method of multipliers introduces. A \emph{formal} proof of Theorem~\ref{thm:Schro} would require just the weaker conditions
$A \in L^2_\textup{loc}(\R^d)$ and $\bm{V}\in L^1_\textup{loc}(\R^d).$ 
	
	The unpleasant conditions~\eqref{additional} can be removed if we consider the situation of potentials $\bm{V}$ and $A$ with \emph{just one singularity} at the origin (see Section~\ref{Sec:one_singularity}). This specific case is worth being investigated as it allows to cover a large class of repulsive potentials, \emph{e.g.}, $\bm{V(x)}= a/ \abs{x}^\alpha \bm{I_{\C^n}}$ 
with $a>0$ and $\alpha>0$, and also the Aharonov--Bohm vector fields~\eqref{eq:A-B} which otherwise would be ruled out by conditions~\eqref{additional}.

\subsubsection{An alternative general result in the self-adjoint setting}

Obviously, Theorem~\ref{thm:Schro} above is valid, with clear simplifications, also in the self-adjoint situation, namely considering \emph{Hermitian} matrix-valued potentials $\bm{V}$.  
In this case, however, we also have an alternative result 
that we have decided to present because the ``repulsivity'' condition~\eqref{hyp:radial_der-rReV^2} is replaced by a ``more classical'' 
assumption in terms of $r\partial_r V^{(2)}.$ 
Furthermore, condition~\eqref{positivity-d=2} is not needed in this context. More precisely we have the following result.
\begin{theorem}\label{thm:self-adjoint}
	Let $d, n\geq 1$ and let $A\in L^2_\textup{loc}(\R^d;\R^d)$ be such that $\bm{B}\in L^2_\textup{loc}(\R^d;\R^{d\times d}).$ Suppose that $\bm{V}\in L^1_\textup{loc}(\R^d; \R^{n\times n})$ admits the decomposition $\bm{V}= \bm{V^{(1)}} + \bm{V^{(2)}}$ with components $\bm{V^{(1)}}\in L^1_\textup{loc}(\R^d)$ and $\bm{V^{(2)}}=V^{(2)}\bm{I_{\C^n}}$, where $V^{(2)}\in L^1_\textup{loc}(\R^d)$ is such that $[r\partial_r V^{(2)}]_+\in L^1_\textup{loc}(\R^d)$ and $r V^{(1)}\in L^2_\textup{loc}(\R^d).$ Assume that there exist numbers $a_1, a_2, b, \mathfrak{b}, c \in [0,1)$ satisfying
	\begin{equation}\label{Schro_cond_numbers:s-a}
	a_1^2 + b^2<1
	\qquad \text{and} \qquad 
	2c+ \mathfrak{b}^2 + d a_1^2 + 2a_2<2
	\end{equation}
	such that, for all $n$-vector $u$ with components in $C^\infty_0(\R^d),$
 \eqref{hyp:V^1} and~\eqref{hyp:B} hold and, moreover,	
	\begin{equation*}
		\int_{\R^d} V^{(2)}_- \abs{u}^2 
\leq b^2 \int_{\R^d} \abs{\nabla_{\!A} u}^2 , 
	\end{equation*}
	\begin{equation}
		\label{hyp:r-radial_der-V^2}
		\int_{\R^d} [r \partial_r V^{(2)})]_+ \abs{u}^2 \leq \mathfrak{b}^2 \int_{\R^d} \abs{\nabla_{\!A} u}^2 .
		\end{equation}	
	If in addition~\eqref{additional} holds true,
	then $H_\textup{S}(A,\bm{V})$ has no eigenvalues, 
\emph{i.e.}\ $\sigma_\textup{p}(H_\textup{S}(A,\bm{V}))=\varnothing.$ 
\end{theorem}  

\begin{remark}
	Here, the first condition in~\eqref{Schro_cond_numbers:s-a} is not explicitly used in the proof of the theorem, but it is needed to give sense to the Hamiltonian $H_\textup{S}(A, \bm{V})$. 
We refer to Section~\ref{preliminary_facts} for details.
\end{remark}

\subsection{The Pauli equation}\label{ss:Pauli}
Recall that the definition of the Pauli operator
depends on the parity of the dimension, 
\emph{cf.}~Lemma~\ref{lemma:explicit-def}.

\begin{theorem}\label{Thm.Pauli}
Let $d\geq 3$ be an integer and let $n'(d)$ be as in~\eqref{eq:n'(d)}. 
Let $A\in L^2_\textup{loc}(\R^d; \R^d)$ 
be such that $\bm{B}\in L^2_\textup{loc}(\R^2;\R^{d\times d}).$ 
Suppose that $\bm{V}\in L^1_\textup{loc}(\R^d; \C^{n'(d)\times n'(d)})$ 
admits the decomposition 
$\bm{V}=\bm{V^{(1)}} + \bm{V^{(2)}}$ with components 
$\bm{V^{(1)}}\in L^1_\textup{loc}(\R^d; \C^{n'(d)\times n'(d)})$ 
and $\bm{V^{(2)}}=V^{(2)}\bm{I_{\C^{n'(d)}}}$,
where $V^{(2)}\in L^1_\textup{loc}(\R^d)$ is such that 
$[\partial_r(r \Re V^{(2)})]_+\in L^1_\textup{loc}(\R^d)$ 
and $r \bm{V^{(1)}}, r(\Re V^{(2)})_-, r\Im V^{(2)}\in L^2_\textup{loc}(\R^d).$
If~$d$ is even, we additionally require 
$\bm{V^{(1)}}=V^{(1)}\bm{I_{\C^{n'(d)}}}$.
Assume that there exist numbers $a, b, \beta, \mathfrak{b}, c\in [0,1)$ 
satisfying
			\begin{equation}\label{cond_numbers-gen}
			\begin{gathered}
				\frac{2}{d-2}\Big(b + \beta + 2\Big(a + \frac{d}{2} c\Big) \Big)<1 ,
		\\ 
		2c + 2\beta + \frac{2(2d-3)}{d-2}\Big(a + \frac{d}{2} c\Big) + \mathfrak{b}^2 + \frac{\sqrt{2}}{\sqrt{d-2}}\Big(b + \Big(a + \frac{d}{2} c\Big)\Big)\Big(\sqrt{\beta} + \sqrt{a + \frac{d}{2} c}\Big)<1 ,
		\end{gathered}		
			\end{equation}
			such that, for all $n'(d)$-vector $u$ 
with components in $C^\infty_0(\R^d),$ the inequalities
			\begin{equation}\label{eq:conditions_B_V}
\int_{\R^d} r^2 \abs{\bm{V^{(1)}}}^2 \abs{u}^2 
\leq a^2\int_{\R^d} \abs{\nabla_{\!A} u}^2 , 
			\qquad 
\int_{\R^d} r^2 \abs{\bm{B}}^2 \abs{u}^2 
\leq c^2\int_{\R^d} \abs{\nabla_{\!A} u}^2 ,
			\end{equation}
			and
			\begin{gather}
			\label{Pauli-cond-V^2}
		\int_{\R^d} r^2 (\Re V^{(2)})_-^2 \abs{u}^2 
\leq b^2 \int_{\R^d} \abs{\nabla_{\!A} u}^2 ,
		\qquad 
		\int_{\R^d} r^2\abs{\Im V^{(2)}}^2 \abs{u}^2 
\leq \beta^2 \int_{\R^d} \abs{\nabla_{\!A} u}^2 ,
		\\
		\label{Pauli-cond-radV^2}
		\int_{\R^d} [\partial_r (r \Re V^{(2)})]_+ \abs{u}^2
\leq \mathfrak{b}^2 \int_{\R^d} \abs{\nabla_{\!A} u}^2 ,
		\end{gather}
		hold true.		
	If, in addition, one has
	\begin{equation*}
		A\in W^{1,d}_\textup{loc}(\R^d)\quad \text{and}\quad \Re V^{(2)}\in W^{1,d/2}_\textup{loc}(\R^d),
	\end{equation*}		
	then $H_\textup{P}(A,\bm{V})$ has no eigenvalues, 
\emph{i.e.}\ $\sigma_\textup{p}(H_\textup{P}(A,\bm{V}))=\varnothing.$
\end{theorem}  

\begin{remark}[Even parity]
	Observe that in the even dimensional case we assume also the component $\bm{V^{(1)}}$ to be diagonal. This is needed in order not to spoil the diagonal form in the definition~\eqref{Pauli-even} of the free Pauli operator, which will represent a crucial point in the strategy underlying the proof (we refer to Section~\ref{Sec:Pauli-Schro} for more details).
\end{remark}

The case of low dimensions $d=1, 2$ is intentionally 
not present in Theorem~\ref{Thm.Pauli} for the following reasons.

\begin{remark}[Dimension one]
As discussed in Section~\ref{rmk:1d}, 
the one-dimensional Pauli operator coincides 
with the scalar potential-free Schr\"odinger operator $-\nabla^2$
(\emph{i.e.}~the one-dimensional Laplacian),
hence the absence of the point spectrum is trivial in this case.
Formally, it is already guaranteed by Theorem~\ref{thm:Schro} with $d=n=1$
(see also Section~\ref{Sec.criticality}). 
\end{remark}	
\begin{remark}[Dimension two]\label{Rem.2D}
The two dimensional case is rather special 
because of the \emph{paramagnetism} of the Pauli operator. 
As a matter of fact, the total absence of the point spectrum 
is no longer guaranteed even in the purely magnetic case 
(\emph{i.e.}\ $\bm{V}=\bm{0}$). 
In this case the Pauli operator has the form 
(see Section~\ref{rmk:Pauli-2d})
	\begin{equation}\label{eq:Pauli2d}
		H_\textup{P}(A, \bm{0})=
		\begin{pmatrix}
			-\nabla_{\!A}^2 + B_{12} & 0\\
			0 & -\nabla_{\!A}^2- B_{12}
		\end{pmatrix}.
	\end{equation} 
For smooth vector potentials, 
the supersymmetry says that the operators 
$-\nabla_{\!A}^2 \pm B_{12}$ have the same spectrum
except perhaps at zero (see \cite[Thm.~6.4]{CFKS}).
Hence the absence of the point spectrum for the two-dimensional 
Pauli operator is in principle governed by 
our Theorem~\ref{thm:Schro} with $d=2$ and $n=1$
(or Theorem~\ref{thm:alternative2d}) 
or its self-adjoint counterpart Theorem~\ref{thm:self-adjoint}
for the special choice $\bm{V} = B_{12} \bm{I}_{\C^2}$.
Unfortunately, we do not see how to derive any non-trivial
condition on~$B_{12}$ to guarantee the \emph{total} absence of eigenvalues 
(\emph{cf}. Remark~\ref{rmk:A-C}).
Physically, it does not come as a big surprise because of 
the celebrated \emph{Aharonov--Casher effect},
which states that the number of zero-eigenstates 
is equal to the integer part of the total magnetic flux
(see \cite[Sec.~6.4]{CFKS}). 
On the one hand, the absence of \emph{negative} eigenvalues
does follow as an immediate consequence of the standard lower bound 
	\begin{equation}\label{B12LB}
		\int_{\R^2} \abs{\nabla_{\!A} u}^2  
\geq \pm \int_{\R^2} B_{12}\abs{u}^2 , 
		\qquad \forall u\in C^\infty_0(\R^2),
	\end{equation} 
which holds with either of the sign $\pm$ 
(see, \emph{e.g.}, \cite[Sec.~2.4]{B_L_S}).

Notice that when an attractive potential is added to 
the two-dimensional Pauli operator, 
it has been proved \cite{Weidl,FMV} that the perturbed Hamiltonian presents always (\emph{i.e.}\ no matter how small is chosen the coupling constant) negative eigenvalues (not only due to the Aharonov--Casher zero modes turning into negative ones, but it is also the essential part of the spectrum that contributes to their appearance). This fact can be seen as a quantification of the aforementioned paramagnetic effect of the Pauli operators in contrast to the diamagnetic effect which holds true for magnetic Schr\"odinger operators.
\end{remark}	

\subsection{The Dirac equation}\label{ss:Dirac} 
Finally, we state our results for the purely magnetic
Dirac operator~\eqref{std_representation}.

\begin{theorem}\label{thm:magnetic_Dirac}
	Let $d\geq 3$ and let $n(d)$ be as in~\eqref{eq:n(d)}. Let $A\in L^2_\textup{loc}(\R^d;\R^d)$ be such that $\bm{B}\in L^2_\textup{loc}(\R^d; \R^{d\times d}).$ Assume that there exists a number $c\in [0,1)$ satisfying
	\begin{equation}
	\label{cond_number_Dirac-gen}
	\frac{2d}{d-2} c<1
	\qquad\mbox{and}\qquad
	2c +\frac{d(2d-3)}{d-2}c  + \frac{\sqrt{2}}{\sqrt{d-2}}\Big(\frac{d}{2}c \Big)^{3/2}<1
	\end{equation}
	such that, for all $n(d)$-vector $u$ with components in $C^\infty_0(\R^d),$ the inequality
	\begin{equation}\label{cond_Dirac_B-higher}
		\int_{\R^d} r^2\abs{\bm{B}}^2\abs{u}^2\leq c^2 \int_{\R^d} \abs{\nabla_{\!A} u}^2
	\end{equation}
	holds true.
If in addition $A\in W^{1,d}_\textup{loc}(\R^d),$
then $H_\textup{D}(A)$ has no eigenvalues, 
\emph{i.e.}\ $\sigma_\textup{p}(H_\textup{D}(A))=\varnothing.$
\end{theorem}

As discussed in Section~\ref{rmk:1d},
the square of the one-dimensional Dirac operator 
is just the one-dimensional Laplacian shifted by a constant
(\emph{cf}.~\eqref{Dirac.1D}), 
hence the absence of the point spectrum follows at once in this case.
On the other hand, the two-dimensional analogue of 
Theorem~\ref{thm:magnetic_Dirac} is unavailable, 
because of the absence of a two-dimensional variant
of Theorem~\ref{Thm.Pauli} in the Pauli case,
\emph{cf}.~Remark~\ref{Rem.2D}.


\section{Scalar electromagnetic Schr\"odinger operators revisited}
\label{Sec:scalar-Schro:1}
%
In this section, we leave aside the operators acting on spinor
Hilbert spaces and focus on \emph{scalar} 
electromagnetic Schr\"odinger operators~\eqref{Schrodinger}. 
This will be useful later on when, in the following sections,
we reduce our analysis to the level of components.
We provide a careful and deep analysis of the method of multipliers, 
stressing on the major outcomes that the technique provides in this context.
Our goal is to represent a reader-friendly overview 
of the original ideas and main outcomes of~\cite{FKV,FKV2} 
to tackle the issue of the total absence of eigenvalues of 
scalar Schr\"odinger operators.
Furthermore, we go through the more technical parts 
by rigorously establishing some results 
that were just sketched in the previous works. 
      
\subsection{Definition of the operators}\label{preliminary_facts}
For the sake of completeness, 
we start with recalling some basic facts 
on the rigorous definition of 
the scalar electromagnetic Schr\"odinger operators. 

Let $d\geq 1$ be any natural number.
Let $A\in L^2_{\textup{loc}}(\R^d;\R^d)$ 
and $V\in L^1_{\textup{loc}}(\R^d; \C)$ 
be respectively a vector potential 
and a scalar potential (the latter possibly complex-valued).
The quantum Hamiltonian apt to describe the motion of a non-relativistic particle interacting with the electric field $-\nabla V$ 
and the magnetic field $\bm{B}:=(\nabla A)- (\nabla A)^t$  
is represented by the scalar electromagnetic Schr\"odinger operator
\begin{equation}\label{eq:ele-magnetic_Hamiltonian}
	H_{A,V}:= -\nabla_{\!A}^2 + V 
\qquad \text{in}\qquad L^2(\R^d).
\end{equation}
Observe that the magnetic field is absent in $\R^1$ 
and~$A$ can be chosen to be equal to zero without loss of generality.
Therefore the two-dimensional framework is the lowest 
in which the introduction of a magnetic field is non-trivial.

As usual, the sum in~\eqref{eq:ele-magnetic_Hamiltonian}
should be understood in the sense of forms
after assuming that~$V$ is relatively form-bounded
with respect to the magnetic Laplacian~$-\nabla_{\!A}^2$
with the relative bound less than one. 
We shall often proceed more restrictively by assuming 
the form-subordination condition
\begin{equation}\label{eq:smallness_assumption}
	\int_{\R^d} \abs{V}\abs{u}^2 
\leq a^2 \int_{\R^d} \abs{\nabla_{\!A} u}^2 , 
 \qquad \forall\, u\in \mathcal{D}_A:=\{u\in L^2(\R^d)\colon \nabla_{\!A} u\in L^2(\R^d)\},
\end{equation}
where $a\in [0,1)$ is a constant independent of~$u$.
Assumption~\eqref{eq:smallness_assumption} in particular implies that the quadratic form
\begin{equation*} 
	h_V[u]:=\int_{\R^d} V \abs{u}^2 , \qquad u\in \mathcal{D}(h_V):=\Big\{u\in L^2(\R^d)\colon \int_{\R^d} \abs{V}\abs{u}^2 <\infty \Big\}
\end{equation*}
is relatively bounded with respect to the quadratic form
\begin{equation*}
	h_A[u]:= \int_{\R^d}\abs{\nabla_{\!A} u}^2 ,
\qquad u\in \mathcal{D}(h_A)=\mathcal{D}_A ,
\end{equation*}
with the relative bound less than one. 
Consequently, the sum $h_{A,V}:=h_A + h_V$ with domain 
$\mathcal{D}(h_{A,V}):=\mathcal{D}_A$
is a closed and sectorial form. 
Therefore $H_{A,V}$ as defined in~\eqref{eq:ele-magnetic_Hamiltonian} 
makes sense as the m-sectorial operator associated to $h_{A,V}$ via the representation theorem (\emph{cf}.~\cite[Thm.~VI.2.1]{Kato}).	


With the aim of including also potentials which are not necessarily subordinated in the spirit of~\eqref{eq:smallness_assumption}, now we present an alternative way to give a meaning to the operator $H_{A,V}$ assuming different conditions on the electric potential $V.$
We introduce the form 
\begin{equation*}
	h_{A,V}^{(1)}[u]:=\int_{\R^d} \abs{\nabla_{\!A} u}^2  + \int_{\R^d} (\Re V)_+ \abs{u}^2 ,\qquad u\in \mathcal{D}(h_{A,V}^{(1)}):=\overline{C^\infty_0(\R^d)}^\normeq{\cdot},
\end{equation*}
with
\begin{equation*}
	\normeq{u}^2:=\int_{\R^d} \abs{\nabla_{\!A} u}^2  + \int_{\R^d} (\Re V)_+ \abs{u}^2  + \int_{\R^d} \abs{u}^2 .
\end{equation*}
The form $h_{A,V}^{(1)}$ is closed by definition. Now instead of assuming  the smallness condition~\eqref{eq:smallness_assumption} for the whole $V$, 
we take the advantage of the splitting in real (positive and negative part) and imaginary part of the potential to require the following more natural subordination: There exist $b, \beta\in [0,1)$ with
\begin{equation}\label{bc}
	b^2+\beta^2<1
\end{equation}
 such that, for any $u\in \mathcal{D}(h_{A,V}^{(1)}),$ 
\begin{equation}\label{Re_-Im}
	\int_{\R^d} (\Re V)_-\abs{u}^2 
\leq b^2 \int_{\R^d} \abs{\nabla_{\!A} u}^2 ,
	\qquad
	\int_{\R^d} \abs{\Im V} \abs{u}^2 
\leq \beta^2 \int_{\R^d} \abs{\nabla_{\!A} u}^2.
\end{equation}
In other words, we require the subordination 
just for the parts $(\Re V)_-$ and $\Im V$ of the potential~$V$.
Hence, defining 
\begin{equation*}
	h_{A,V}^{(2)}[u]:= 
-\int_{\R^d} (\Re V)_- \abs{u}^2  + i \int_{\R^d} \Im V \abs{u}^2 ,
\end{equation*}
the form $h_{A,V}^{(2)}$ is relatively bounded
with respect to  $h_{A,V}^{(1)}$, 
with the relative bound less than one (see~\eqref{bc}).
Consequently, as above, the sum $h_{A,V}=h_{A,V}^{(1)} + h_{A,V}^{(2)}$ is a closed and sectorial form and $\mathcal{D}(h_{A,V})=\mathcal{D}(h_{A,V}^{(1)}).$ Therefore, also in this more general setting, $H_{A,V}$ is the m-sectorial operator associated with $h_{A,V}.$
 
In order to consider simultaneously both these two possible configurations,
we introduce the decomposition $V=V^{(1)} + V^{(2)}$ 
and assume that there exist $a, b, \beta\in [0,1)$ satisfying
\begin{equation}\label{eq:a-b-c} 
  a^2+ b^2+\beta^2<1
\end{equation} 
such that, for any $u\in \mathcal{D}_A,$
\begin{equation}
	\label{eq:smallness_assumption_V1}
	\int_{\R^d} \abs{V^{(1)}} \abs{u}^2  
\leq a^2\int_{\R^d}\abs{\nabla_{\!A} u}^2 ,
\end{equation}
and 
\begin{equation}
	\label{eq:assumption_Re-_Im}
	\int_{\R^d} (\Re V^{(2)})_-\abs{u}^2 
\leq b^2 \int_{\R^d} \abs{\nabla_{\!A} u}^2 ,
	\qquad
	\int_{\R^d} \abs{\Im V^{(2)}} \abs{u}^2  
\leq \beta^2 \int_{\R^d} \abs{\nabla_{\!A} u}^2.
\end{equation}
Let us define
$ h_{A,V}^{(1)}[u]:= \int_{\R^d} \abs{\nabla_{\!A} u}^2 
+ \int_{\R^d} (\Re V^{(2)})_+ \abs{u}^2$  
with 
$\mathcal{D}(h_{A,V}^{(1)}):= \overline{C^\infty_0(\R^d)}^{\normeq{\cdot}},$
where
\begin{equation*}
	\normeq{u}^2:=\int_{\R^d} \abs{\nabla_{\!A} u}^2 
+ \int_{\R^d} (\Re V^{(2)})_+ \abs{u}^2  + \int_{\R^d} \abs{u}^2 ,
\end{equation*}	
and $h_{A,V}^{(2)}[u]:= \int_{\R^d} V^{(1)}\abs{u}^2 
  -\int_{\R^d} (\Re V^{(2)})_-\abs{u}^2  
+ i \int_{\R^d} \Im V^{(2)}\abs{u}^2 $
with $\Dom(h_{A,V}^{(2)}):=\Dom(h_{A,V}^{(1)})$.
By the same reasoning as above, 
one has that $H_{A,V}$ is the m-sectorial operator associated 
with the closed and sectorial form 
$h_{A,V} := h_{A,V}^{(1)}+h_{A,V}^{(2)}$
with $\mathcal{D}(h_{A,V}):=\mathcal{D}(h_{A,V}^{(1)}).$ 
	In order to drop the dependance on the form $h$ in the notation of the domain that will not be used explicitly any more, from now on we will denote
	\begin{equation*}
		\mathcal{D}_{A,V}:= \mathcal{D}(h_{A,V}).
	\end{equation*}

\subsection{Further hypotheses on the potentials}\label{Further hypotheses A,V}
As we shall see below,
in order to justify rigorously the algebraic manipulations that the method of multipliers introduces, we need to assume more regularity 
on the magnetic potential~$A$ and on the electric potential 
$V=V^{(1)} + V^{(2)}$ than the ones required to give a meaning to the electromagnetic Hamiltonian~\eqref{eq:ele-magnetic_Hamiltonian}.

\subsubsection{Further hypotheses on the magnetic potential}
 We assume
\begin{equation}\label{eq:hypotheses_A}
	A\in W^{1,p}_\textup{loc}(\R^d; \R^d)
	\quad \text{where}\quad
	\begin{system}
	&p=2  &\text{if}\quad d=1,\\
	&p>2  &\text{if}\quad d=2,\\
	&p=d  &\text{if}\quad d\geq 3.
	\end{system}
\end{equation}
In particular, these assumptions ensure that for any $u\in \mathcal{D}_A$ then 
\begin{equation}\label{Au_L2}
Au\in L_\textup{loc}^2(\R^d; \R^d)
\end{equation}
and the same can be said for $\partial_l A u,$ with $l=1,2,\dots, d.$
 Indeed, from the H\"older inequality, one has that for any $k=1,2,\dots, d$
\begin{equation}\label{Au-L^2}
	\norm{A_k u}_{L_\textup{loc}^2(\R^d)}^2\leq \norm{A_k}_{L_\textup{loc}^p(\R^d)}\norm{u}_{L_\textup{loc}^q(\R^d)} \quad \text{with} \quad 1/p+1/q=1/2.
\end{equation}
Observe that the diamagnetic inequality~\eqref{eq:diamagnetic}
and $u\in \mathcal{D}_A$ guarantee $\abs{u}\in H^1(\R^d).$  
By the Sobolev embeddings 
\begin{equation}\label{Sobolev_embeddings}
H^1(\R^d)\hookrightarrow L^q(\R^d)
\quad \text{where}\quad
	\begin{system}
		& q=\infty \quad &\text{if}\quad d=1,\\
		& 2\leq q<\infty \quad &\text{if}\quad d=2,\\
		& q=2^\ast:=2d/(d-2) \quad &\text{if}\quad d\geq 3.
	\end{system}
\end{equation}
Consequently, if one chooses $q$ as in~\eqref{Sobolev_embeddings}, 
then $\norm{u}_{L^q(\R^d)}$ is finite. 
If, moreover, the H\"older conjugated exponent~$p$ 
is as in our assumption~\eqref{eq:hypotheses_A}, 
then $\norm{A_k}_{L_\textup{loc}^p(\R^d)}$ is finite 
and therefore, from~\eqref{Au-L^2}, $A_k u\in L_\textup{loc}^2(\R^d).$

Notice that, given any function $u\in \mathcal{D}_A$ as soon as $Au\in L^2(\R^d),$ then $\nabla u\in L^2(\R^d)$ and therefore $u\in H^1(\R^d).$ In other words
\begin{equation}\label{inclusion}
	\{u\in \mathcal{D}_A\ \&\ Au\in L^2(\R^d)\}\subseteq H^1(\R^d).
\end{equation}  

\subsubsection{Further hypotheses on the electric potential}
Recalling the decomposition $V = V^{(1)}+V^{(2)}$,
we assume the following condition on 
the real part of the second component:
\begin{equation}\label{eq:hypotheses_V_alt}
	\Re V^{(2)}\in W^{1,p}_\textup{loc}(\R^d; \R)
	\quad \text{where}\quad
	\begin{system}
	&p=1  &\text{if}\quad d=1,\\
	&p>1  &\text{if}\quad d=2,\\
	&p=d/2  &\text{if}\quad d\geq 3
	\end{system}
\end{equation} 
By the same reasoning as done above for the magnetic potential, one can observe that assumption~\eqref{eq:hypotheses_V_alt} ensures that for any $u\in H^1_A(\R^d),$ then
\begin{equation*}
	\Re V^{(2)}\abs{u}^2\in L_\textup{loc}^1(\R^d),
	\end{equation*}
and the same can be said for $\partial_k \Re V^{(2)},$ with $k=1,2,\dots, d.$

\subsection{The method of multipliers: main ingredients}

The purpose of this subsection is to provide, in a unified and rigorous way, the proof of the common crucial starting point of the series of works~\cite{FKV,FKV2,Cossetti_2017,Cossetti-Krejcirik_2018} for proving the absence of the point spectrum of the electromagnetic Hamiltonians $H_{A,V}$ in various settings.

Since this section is intended as a review of already known results on \emph{scalar} Schr\"odinger Hamiltonians, here we will be concerned almost exclusively with the most interesting and more troublesome case of the spectral parameter $\lambda\in \C$ within the sector of the complex plane given by
\begin{equation}\label{sector}
	\{\lambda \in \C\colon \Re \lambda \geq \abs{\Im \lambda}\}.
\end{equation}
On the other hand, how to deal with the complementary sector, \emph{i.e.}, $\{\lambda \in \C\colon \Re \lambda < \abs{\Im \lambda}\}$ can be seen explicitly in the proof of our original results 
(see Sections~\ref{Sec:Schro} and~\ref{sec:Pauli-Dirac}). 

The proof of the absence of eigenvalues within the sector defined in~\eqref{sector} is based on the following crucial result obtained by means of 
the method of multipliers. 
It basically provides an integral identity for 
weak solutions~$u$ to the resolvent equation $(H_{A,V}-\lambda) u = f$,
where $f\colon \R^d \to \C$ is a suitable function.
More specifically, $u\in \mathcal{D}_{A,V}$ is such that the identity
\begin{equation}\label{eq:weak_resolvent}
	\int_{\R^d} \nabla_{\!A} u \cdot \overline{\nabla_{\!A} v}  + \int_{\R^d} V u \bar{v}  = \lambda \int_{\R^d}u \bar{v}  + \int_{\R^d} f \bar{v} 
\end{equation}
holds for any $v\in \mathcal{D}_{A,V},$ 
where~$f$ is any suitable function for which the last integral in~\eqref{eq:weak_resolvent} is finite. 
The crucial result reads as follows.

\begin{lemma}\label{lemma:crucial_tool}
	Let $d\geq 1$, let $A\in L^2_\textup{loc}(\R^d;\R^d)$ be such that $\bm{B}\in L^2_\textup{loc}(\R^d; \R^{d\times d})$ and~\eqref{eq:hypotheses_A} holds. Suppose that $V\in L^1_\textup{loc}(\R^d;\C)$ admits the decomposition $V=V^{(1)}+ V^{(2)}$ with $\Re V^{(2)}$ satisfying~\eqref{eq:hypotheses_V_alt}.
		Let $u\in \mathcal{D}_{A,V}$ be a solution to~\eqref{eq:weak_resolvent}, with $\abs{\Im \lambda}\leq \Re \lambda$ and $r f\in L^2(\R^d),$ satisfying 
		\begin{equation*}
			\left( 
			r^2 \abs{V^{(1)}}^2
			+ r^2 (\Re V^{(2)})_-^2
			+[\partial_r(r\Re V^{(2)})]_+ 
			+ r^2 \abs{\Im V^{2}}^2
			+r^2 \abs{\bm{B}}^2
			\right)\abs{u}^2\in L^1(\R^d).
		\end{equation*}
		Then also $r\abs{\nabla_{\!A} u^-}^2 + r^{-1}\abs{u}^2 + [\partial_r(r \Re V^{(2)})]_-\abs{u}^2 + r[\Re V^{(2)}]_+\abs{u}^2 \in L^1(\R^d)$ and the identity
\begin{multline}\label{eq:crucial_identity}
	\int_{\R^d} \abs{\nabla_{\!A} u^-}^2\, dx
	+(\Re \lambda)^{-1/2} \abs{\Im \lambda} \int_{\R^d} \abs{x}\abs{\nabla_{\!A} u^-}^2\, dx
	-\frac{(d-1)}{2}(\Re \lambda)^{-1/2}\abs{\Im \lambda} \int_{\R^d} \frac{\abs{u}^2}{\abs{x}}\,dx\\ 
	+2\Im \int_{\R^d} x \cdot \bm{B} \cdot u^- \overline{\nabla_{\!A} u^-}\, dx\\
	+(d-1)\int_{\R^d} \Re V^{(1)}\abs{u}^2\, dx
	+2\Re \int_{\R^d} x\cdot V^{(1)} u^- \overline{\nabla_{\!A} u^-}\, dx
	+(\Re \lambda)^{-1/2} \abs{\Im \lambda} \int_{\R^d}\abs{x} \Re V^{(1)}\abs{u}^2\, dx\\
	-\int_{\R^d} \partial_r (\abs{x} \Re V^{(2)})\abs{u}^2\, dx
	-2\Im \int_{\R^d} x \Im V^{(2)} u^- \overline{\nabla_{\!A} u^-}\, dx
	+(\Re \lambda)^{-1/2} \abs{\Im \lambda}\int_{\R^d}\abs{x} \Re V^{(2)} \abs{u}^2\, dx\\
	=(d-1) \Re \int_{\R^d} f \bar{u}\, dx
	+ 2 \Re \int_{\R^d} x \cdot f^-  \overline{\nabla_{\!A} u^-}\, dx
	+(\Re \lambda)^{-1/2} \abs{\Im \lambda} \Re \int_{\R^d} \abs{x} f \bar{u}\, dx
	\end{multline}
holds true with 
\begin{equation}\label{eq:u^-}
	u^-(x):= e^{-i(\Re \lambda)^{1/2}\sgn(\Im \lambda)\abs{x}} u(x)
\end{equation}
and $f^-$ defined in the analogous way.
\end{lemma}
\begin{remark}[Dimension one]\label{rmk:dim1}
	Since the addition of a magnetic potential is trivial in $\R^1,$ the corresponding identity~\eqref{eq:crucial_identity} with $d=1$ comes with the classical gradient $\nabla$ as a replacement of the magnetic gradient $\nabla_{\!A},$ moreover the term involving $\bm{B}$ is not present.
\end{remark}

The proof of Lemma~\ref{lemma:crucial_tool}
can be found in Subsection~\ref{Proof of Lemma_identity}, 
here we just provide its main steps:
\begin{itemize}
	\item \textbf{Step one:} Approximation of $u$ with a sequence of compactly supported functions $u_R$ (see definition~\eqref{eq:approx_def} below) which satisfy a related problem with  small (in a suitable topology) corrections. This first step is necessary in order to justify rigorously the algebraic manipulations that the method of multipliers introduces when the test function $v$ is chosen to be possibly \emph{unbounded} (so that it is not even \emph{a priori} clear if this specific choice $v$ belongs to $L^2(\R^d)$). 
	\item \textbf{Step two:} Development of the method of multipliers for $u_R$ (main core of the proof) in order to produce the analogue of identity~\eqref{eq:crucial_identity} for the approximating sequence. This step will require a further approximation procedure which will ensure that the chosen multiplier $v$ (see~\eqref{eq:multiplier_explicit} below) is in $\mathcal{D}_{A,V}$ and therefore allowed to be taken as a test function.
	\item \textbf{Step three:} Proof of~\eqref{eq:crucial_identity} by taking the limit as $R \to \infty$
in the previous identity and using the smallness of the corrections which is quantified in Lemma~\ref{lemma:f} below.
\end{itemize}

As a byproduct of the crucial identity of Lemma~\ref{lemma:crucial_tool}, 
we get the following inequality.
For the sake of completeness, we provide it with a proof. 

\begin{lemma}\label{lemma:byproduct}
	Under the hypotheses of Lemma~\ref{lemma:crucial_tool} the following estimate
	\begin{multline}\label{est_Schr}
		\norm{\nabla_{\!A} u^-}_{L^2(\R^d)}^2 
		+ (\Re \lambda)^{-1/2} \abs{\Im \lambda} 
		\Bigg[
		\int_{\R^d} \abs{x} \abs{\nabla_{\!A} u^-}^2\, dx
		- \frac{(d-1)}{2} \int_{\R^d} \frac{~\abs{u^-}^2}{\abs{x}}\, dx
		+ \int_{\R^d} \abs{x} (\Re V^{(2)})_+ \abs{u^-}^2\, dx 
		\Bigg]
		\\
		\leq 
		2\Big( 
		 \norm{\abs{x}\abs{\bm{B}} u^-}_{L^2(\R^d)} 
		 +\norm{\abs{x} \Im V^{(2)} u^-}_{L^2(\R^d)}
		 +\norm{\abs{x} f}_{L^2(\R^d)}
		\Big) \norm{\nabla_{\!A} u^-}_{L^2(\R^d)}
		\\
		+(d-1) \norm{\abs{f}^{1/2} \abs{u^-}^{1/2}}_{L^2(\R^d)}^2
		+ \norm{[\partial_r(\abs{x}\Re V^{(2)})]_+^{1/2} u^-}_{L^2(\R^d)}^2
		\\
		+ \Big(
		\norm{\abs{x} (\Re V^{(2)})_- u^-}_{L^2(\R^d)}
		+\norm{\abs{x} f}_{L^2(\R^d)}
		\Big)
		\Big(
		\norm{\abs{\Im V^{(2)}}^{1/2} u^-}_{L^2(\R^d)}
		+\norm{\abs{f}^{1/2} \abs{u^-}^{1/2}}_{L^2(\R^d)}
		\Big)
	\end{multline}
	holds true.
\end{lemma}
 
\begin{proof}[Proof of Lemma~\ref{lemma:byproduct}]
 	Let us consider identity~\eqref{eq:crucial_identity} with $V^{(1)}=0.$ In passing, notice that requiring $V^{(1)}=0$ do not entails any loss of generality. Indeed since, according to our notations, $V^{(1)}$  represents the component of the electric potential $V$ which is \emph{fully subordinated} to the magnetic Dirichlet form (in the sense given by~\eqref{eq:smallness_assumption_V1}), it can be treated at the same level of the forcing term $f.$
	
				After splitting $\Re V^{(2)}$ in its positive and negative parts, namely using $\Re V^{(2)}= (\Re V^{(2)})_+ - (\Re V^{(2)})_-,$ identity~\eqref{eq:crucial_identity} with $V^{(1)}=0$ reads as follows
	\begin{multline}\label{simplified}
		\int_{\R^d} \abs{\nabla_{\!A} u^-}^2\, dx
	+(\Re \lambda)^{-1/2} \abs{\Im \lambda} 
	\Bigg[
	\int_{\R^d} \abs{x}\abs{\nabla_{\!A} u^-}^2\, dx
	-\frac{(d-1)}{2} \int_{\R^d} \frac{~\abs{u}^2}{\abs{x}}\,dx
	+\int_{\R^d}\abs{x} (\Re V^{(2)})_+ \abs{u}^2\, dx
	\Bigg]
	\\ 
	=-2\Im \int_{\R^d} x \cdot \bm{B} \cdot u^- \overline{\nabla_{\!A} u^-}\, dx
	\\
	+\int_{\R^d} \partial_r (\abs{x} \Re V^{(2)})\abs{u}^2\, dx
	+2\Im \int_{\R^d} x \Im V^{(2)} u^- \overline{\nabla_{\!A} u^-}\, dx
	+(\Re \lambda)^{-1/2} \abs{\Im \lambda}\int_{\R^d}\abs{x} (\Re V^{(2)})_- \abs{u}^2\, dx\\
	+(d-1) \Re \int_{\R^d} f \bar{u}\, dx
	 +2 \Re \int_{\R^d} x f^-  \overline{\nabla_{\!A} u^-}\, dx
	+(\Re \lambda)^{-1/2} \abs{\Im \lambda} \Re \int_{\R^d} \abs{x} f \bar{u}\, dx.
	\end{multline}
	

We consider first
\begin{equation*}
	I:=-2\Im \int_{\R^d} x\cdot \bm{B} \cdot u^-\overline{\nabla_{\!A} u^-}\,dx.
\end{equation*}
By the Cauchy--Schwarz inequality, it immediately follows that
\begin{equation}\label{control_B}
	\abs{I}\leq 2\norm{\abs{x}\abs{\bm{B}} u^-}_{L^2(\R^d)} \norm{\nabla_{\!A} u^-}_{L^2(\R^d)}.
\end{equation}
Now we consider the terms in~\eqref{simplified} involving $V^{(2)},$ that is 
\begin{equation*}
	\begin{split}
	I\!I := & \
	\int_{\R^d} \partial_r (\abs{x} \Re V^{(2)})\abs{u}^2\, dx
	+2\Im \int_{\R^d} x \Im V^{(2)} u^- \overline{\nabla_{\!A} u^-}\, dx
	+(\Re \lambda)^{-1/2} \abs{\Im \lambda}\int_{\R^d}\abs{x} (\Re V^{(2)})_- \abs{u}^2\, dx\\
		= & \ I\!I_1 + I\!I_2+ I\!I_3.
	\end{split}
\end{equation*}
Using that $\abs{u}=\abs{u^-},$ the term $I\!I_1$ can be easily estimated in this way:
\begin{equation}\label{II_1}
	I\!I_1\leq \int_{\R^d} [\partial_r (\abs{x} \Re V^{(2)})]_+ \abs{u}^2\,dx= \norm{[\partial_r (\abs{x} \Re V^{(2)})]_+^{1/2} u^-}_{L^2(\R^d)}^2.
\end{equation}
By the Cauchy--Schwarz inequality one has
\begin{equation}\label{II_2}
	I\!I_2\leq\abs{I\!I_2}\leq 2 \norm{\abs{x} \Im V^{(2)} u^-}_{L^2(\R^d)} \norm{\nabla_{\!A} u^-}_{L^2(\R^d)}.
\end{equation}
Finally, if $\Im \lambda\neq 0,$ we also need to estimate $I\!I_3.$ First notice that choosing $v=\frac{\Im \lambda}{\abs{\Im \lambda}} u$ in~\eqref{eq:weak_resolvent} (with $V^{(1)}=0$) and taking the imaginary part of the resulting identity, gives the following $L^2$- bound
\begin{equation}\label{L^2-bound}
	\norm{u}_{L^2(\R^d)}\leq \abs{\Im \lambda}^{-1/2}
	\Big( 
	\norm{\abs{\Im V^{(2)}}^{1/2} u}_{L^2(\R^d)}
	+\norm{\abs{f}^{1/2} \abs{u}^{1/2}}_{L^2(\R^d)}
	\Big).
\end{equation}
Using the Cauchy--Schwarz inequality, the $L^2$-bound~\eqref{L^2-bound}, the fact that we are working in the sector $\abs{\Im\lambda}\leq \Re \lambda,$ and again using that $\abs{u}=\abs{u^-},$ we have
\begin{equation}\label{II_3}
	\begin{split}
	I\!I_3
	&\leq (\Re \lambda)^{-1/2} \abs{\Im \lambda}\norm{\abs{x}[\Re V^{(2)}]_- u}_{L^2(\R^d)}\norm{u}_{L^2(\R^d)}\\
	&\leq \norm{\abs{x}[\Re V^{(2)}]_- u^-}_{L^2(\R^d)} 
	\Big( 
	\norm{\abs{\Im V^{(2)}}^{1/2} u^-}_{L^2(\R^d)}
	+\norm{\abs{f}^{1/2} \abs{u^-}^{1/2}}_{L^2(\R^d)}
	\Big).
	\end{split}
\end{equation}
Now we estimate the terms in~\eqref{simplified} involving $f,$ namely
\begin{equation*}
	\begin{split}
	I\!I\!I :=& \ (d-1) \Re \int_{\R^d} f \bar{u}\, dx
	 +2 \Re \int_{\R^d} x f^-  \overline{\nabla_{\!A} u^-}\, dx
	+(\Re \lambda)^{-1/2} \abs{\Im \lambda} \Re \int_{\R^d} \abs{x} f \bar{u}\, dx\\
	=& \ I\!I\!I_1+I\!I\!I_2+I\!I\!I_3.
	\end{split}
\end{equation*}
In a similar way as done to estimate $I\!I_1, I\!I_2$ and $I\!I_3,$ one gets
\begin{equation}\label{III_1/2}
	I\!I\!I_1\leq (d-1) \norm{\abs{f}^{1/2} \abs{u^-}^{1/2}}_{L^2(\R^d)}^2,
		\qquad
	I\!I\!I_2\leq 2 \norm{\abs{x} f}_{L^2(\R^d)}\norm{\nabla_{\!A} u^-}_{L^2(\R^d)}
\end{equation}
and
\begin{equation}\label{III_3}
	\begin{split}
		I\!I\!I_3&\leq (\Re \lambda)^{-1/2} \abs{\Im \lambda} \norm{\abs{x}f}_{L^2(\R^d)} \norm{u}_{L^2(\R^d)}\\
		&
	\leq  \norm{\abs{x} f}_{L^2(\R^d)}
	\Big( 
	\norm{\abs{\Im V^{(2)}}^{1/2} u^-}_{L^2(\R^d)}
	+\norm{\abs{f}^{1/2} \abs{u^-}^{1/2}}_{L^2(\R^d)}
	\Big).
	\end{split}
\end{equation}

Applying estimates~\eqref{control_B},~\eqref{II_1},\eqref{II_2} and~\eqref{II_3} together with~\eqref{III_1/2} and~\eqref{III_3} in~\eqref{simplified}, we obtain the thesis.
\end{proof}

Now we are in a position to prove Lemma~\ref{lemma:crucial_tool} on the basis of the three steps presented above.

\subsubsection{Proof of Lemma~\ref{lemma:crucial_tool}}\label{Proof of Lemma_identity}
\begin{itemize}
	\item \textbf{Step one.}
		The desired approximation by compactly supported functions is achieved by a usual ``horizontal cut-off.'' Let $\mu\colon [0,\infty)\to [0,1]$ be a smooth function such that
		\begin{equation*}
			\mu(r)=
			\begin{system}
				&1 &\text{if}\quad 0\leq r\leq 1,\\
				&0  			&\text{if}\quad r\geq 2.
			\end{system}
		\end{equation*}
		Given a positive number~$R$,
we set $\mu_R(x):=\mu(\abs{x}R^{-1}).$ 
Then $\mu_R\colon \R^d \to [0,1]$ is such that
\begin{equation}\label{def:muR}
  \mu_R = 1 \quad \text{in}\quad B_R(0), \qquad 
  \mu_R= 0 \quad \text{in}\quad \R^d\setminus B_{2R}(0), \qquad 
  \abs{\nabla \mu_R}\leq cR^{-1}, \qquad 
  \abs{\Delta \mu_R}\leq c R^{-2},
\end{equation}
where $B_R(0)$ stands for the open ball centered at the origin and with radius $R>0$ and $c> 1$ is a suitable constant independent of~$R$.
For any function $h\colon \R^d \to \C$ 
we then define the compactly supported approximating family of functions 
by setting
\begin{equation}\label{eq:approx_def}
	h_R:=\mu_R h.
\end{equation}
If $u\in \mathcal{D}_{A,V}$ is a weak solution to $-\nabla_{\!A}^2 u + Vu=\lambda u + f$, it is not difficult to show that the compactly supported 
function~$u_R$ belongs to $\mathcal{D}_{A,V}$ and solves in a weak sense the following related problem 
\begin{equation}\label{eq:resolvent}
	-\nabla_{\!A}^2u_R + Vu_R=\lambda u_R + f_R + \textup{err}(R) \quad \text{in}\quad \R^d,
\end{equation}		
where 
\begin{equation}\label{eq:error_terms}
	\textup{err}(R):= -2\nabla_{\!A} u \cdot \nabla \mu_R - u \Delta \mu_R.
\end{equation}		
The next easy result shows that the extra terms~\eqref{eq:error_terms}, 
which originate from the introduction of the horizontal cut-off $\mu_R$,
become negligible as $R$ increases.
\begin{lemma}\label{lemma:f}
Given $u\in \mathcal{D}_{A,V}$,
let $\textup{err}(R)$ be as in~\eqref{eq:error_terms}. 
Then the following limits
\begin{align*} 
  &\norm{\textup{err}(R)}_{L^2(\R^d)}\xrightarrow{R\to \infty} 0, 
  & 
  \norm{\abs{x} \textup{err}(R)}_{L^2(\R^d)}\xrightarrow{R\to \infty} 0
\end{align*}
hold true.
\end{lemma}
\begin{proof}
By~\eqref{def:muR} we have 
\begin{align*}
\norm{\textup{err}(R)}_{L^2(\R^d)}
&\leq 2 \left(\int_{\R^d} \abs{\nabla_{\!A} u}^2 \abs{\nabla \mu_R}^2\right)^{1/2} 
+ \left(\int_{\R^d} \abs{u}^2 \abs{\Delta \mu_R}^2\right)^{1/2}
\\
&\leq \frac{2c}{R} \left(\int_{\{R<\abs{x}<2R\}} \abs{\nabla_{\!A} u}^2\right)^{1/2} 
+ \frac{c}{R^2} \left(\int_{\{R<\abs{x}<2R\}} \abs{u}^2\right)^{1/2}.  
\end{align*}
Since $u\in L^2(\R^d)$ and $\nabla_{\!A} u\in \big[L^2(\R^d)\big]^d$, 
the right-hand side tends to zero as~$R$ goes to infinity.

Similarly,
\begin{align*}
\norm{\abs{x} \textup{err}(R)}_{L^2(\R^d)}
&\leq 2\left (\int_{\R^d} \abs{x}^2 \abs{\nabla_{\!A} u}^2 \abs{\nabla \mu_R}^2\right)^{1/2} 
+ \left(\int_{\R^d} \abs{x}^2 \abs{u}^2 \abs{\Delta \mu_R}^2\right)^{1/2}
\\
&\leq  4c \left(\int_{\{R< \abs{x}<2R\}} \abs{\nabla_{\!A} u}^2\right)^{1/2} 
+ \frac{2c}{R}\left(\int_{\{R<\abs{x}<2R\}} \abs{u}^2\right)^{1/2},
\end{align*}
and again the right-hand side goes to zero as $R$ approaches infinity.
	\end{proof}

	\item \textbf{Step two.}
		This second step represents the main body of the section, it is here that the method of multipliers is fully developed. Informally speaking the method of multipliers is based on producing integral identities
by choosing different test functions~$v$ in~\eqref{eq:weak_resolvent} (see Lemma~\ref{lemma:crucial_identities} below)
and later combining them in a refined way to get, for instance in our case, the analogous to~\eqref{eq:crucial_identity}.
		By virtue of the previous step, we shall develop the method for \emph{compactly supported} solutions $u\in \mathcal{D}_{A,V}$ to~\eqref{eq:weak_resolvent}, it will be in the next \textbf{Step three} that we will get the result also for not necessarily compactly supported solutions.
		
As a starting point we state the aforementioned identities, 
these are collected in the following lemma.
Notice that the lemma is stated for any $\lambda\in \C$ and not necessarily just for $\lambda$ in the sector~\eqref{sector}.

\begin{lemma}\label{lemma:crucial_identities}
	Let $d\geq 1,$ let $A\in L^2_\textup{loc}(\R^d;\R^d)$ be such that $\bm{B}\in L^2_\textup{loc}(\R^d; \R^{d\times d})$ and assume also~\eqref{eq:hypotheses_A}, Suppose that $V\in L^1_\textup{loc}(\R^d;\C)$ admits the decomposition $V=V^{(1)}+ V^{(2)}$ with $\Re V^{(2)}$ satisfying~\eqref{eq:hypotheses_V_alt}. Let $u\in \mathcal{D}_{A,V}$ be any compactly supported solution of~\eqref{eq:weak_resolvent}, with $\lambda$ any complex constant and $\abs{x}f\in L^2_\textup{loc}(\R^d),$ satisfying
	\begin{equation}\label{x_hypothesis}
		\left(
		\abs{x}^2 \abs{V^{(1)}}^2 
		+ \abs{x}^2 \abs{\Im V^{(2)}}^2
		\right) 
		\abs{u}^2\in L^1_\textup{loc}(\R^d). 
	\end{equation}
	Then $\abs{x}^{-1}\abs{u}^2\in L^1_\textup{loc}(\R^d)$ and the following identities 
	\begin{equation}\label{eq:first_const}
		\int_{\R^d} \abs{\nabla_{\!A} u}^2\, dx + \int_{\R^d} \Re V \abs{u}^2\, dx= \Re \lambda \int_{\R^d} \abs{u}^2\, dx + \Re \int_{\R^d} f \bar{u}\, dx. 
	\end{equation}
	\begin{equation}\label{eq:first_x}
		-\frac{d-1}{2} \int_{\R^d} \frac{~\abs{u}^2}{\abs{x}}\, dx + \int_{\R^d} \abs{x} \abs{\nabla_{\!A} u}^2\, dx + \int_{\R^d} \Re V \abs{x} \abs{u}^2\, dx= \Re \lambda \int_{\R^d} \abs{x} \abs{u}^2\, dx + \Re \int_{\R^d} f \abs{x} \bar{u}\, dx. 
	\end{equation}
	\begin{equation}\label{eq:second_const}
		 \int_{\R^d} \Im V \abs{u}^2\, dx= \Im \lambda \int_{\R^d} \abs{u}^2\, dx + \Im \int_{\R^d} f \bar{u}\, dx.
	\end{equation}
	\begin{equation}\label{eq:second_x}
		\Im \int_{\R^d} \frac{x}{\abs{x}}\cdot \bar{u}\nabla_{\!A} u\, dx + \int_{\R^d} \Im V \abs{x} \abs{u}^2\, dx= \Im \lambda \int_{\R^d} \abs{x} \abs{u}^2\, dx + \Im \int_{\R^d} f \abs{x} \bar{u}\, dx.
	\end{equation}
	\begin{multline}\label{eq:third}
		2\int_{\R^d} \abs{\nabla_{\!A} u}^2\, dx 
		+ 2 \Im \int_{\R^d} x\cdot \bm{B} \cdot u \overline{\nabla_{\!A} u}\, dx 
		+ d \int_{\R^d} \Re V^{(1)} \abs{u}^2\, dx 
		+ 2 \Re \int_{\R^d} x\cdot V^{(1)} u \overline{\nabla_{\!A} u}\,dx\\
		-\int_{\R^d} x\cdot \nabla \Re V^{(2)} \abs{u}^2\, dx
		-2 \Im \int_{\R^d} x\cdot \Im V^{(2)} u \overline{\nabla_{\!A} u}\, dx\\
		=-2\Im \lambda \Im \int_{\R^d} x \cdot u \overline{\nabla_{\!A} u}\, dx + d\Re \int_{\R^d} f \bar{u}\, dx + 2 \Re \int_{\R^d} f x\cdot \overline{\nabla_{\!A} u}\, dx.
	\end{multline}
	hold true.
	\end{lemma}
 
	Now we show how to use these identities to prove the analogous of identity~\eqref{eq:crucial_identity} for \emph{compactly supported} solutions of~\eqref{eq:weak_resolvent}. For the sake of clarity, the technical proof of Lemma~\ref{lemma:crucial_identities} is postponed to Subsection~\ref{subsec:proof_identities}.
	
	Let us start our algebraic manipulation of identities~\eqref{eq:first_const}--\eqref{eq:third} by taking the sum
	\begin{equation*}
		-~\eqref{eq:first_const} -2(\Re \lambda)^{1/2}\sgn(\Im \lambda)~\eqref{eq:second_x} +~\eqref{eq:third}.
	\end{equation*}
	
	This gives
	\begin{multline}\label{eq:intermediate}
	\hspace{2cm}
	\int_{\R^d} \abs{\nabla_{\!A} u}^2\, dx 
	-2 (\Re\lambda)^{1/2} \sgn(\Im\lambda) \Im \int_{\R^d} \frac{x}{\abs{x}} \cdot \bar{u}\,\nabla_{\!A} u\, dx + \Re \lambda \int_{\R^d} \abs{u}^2\, dx\\
	+2(\Re \lambda)^{1/2} \abs{\Im \lambda} \int_{\R^d} \abs{x} \abs{u}^2\, dx + 2 \Im\lambda \Im \int_{\R^d} x \cdot u\, \overline{\nabla_{\!A} u}\, dx
	+2\Im \int_{\R^d} x\cdot \bm{B}\cdot u \, \overline{\nabla_{\!A} u}\, dx\\
	+(d-1)\int_{\R^d} \Re V^{(1)} \abs{u}^2\, dx 
	+2\Re \int_{\R^d} x\cdot V^{(1)} u \, \overline{\nabla_{\!A} u}\,dx
	-2(\Re \lambda)^{1/2} \sgn(\Im \lambda) \int_{\R^d} \abs{x} \Im V^{(1)} \abs{u}^2\, dx\\
	-\int_{\R^d} \Re V^{(2)} \abs{u}^2\,dx 
	-\int_{\R^d} x\cdot \nabla \Re V^{(2)} \abs{u}^2\, dx\\
	-2(\Re \lambda)^{1/2} \sgn(\Im \lambda) \int_{\R^d} \abs{x} \Im V^{(2)} \abs{u}^2\, dx
		-2\Im \int_{\R^d} x\cdot \Im V^{(2)} u \overline{\nabla_{\!A} u}\, dx
	\\
	= (d-1)\Re \int_{\R^d} f \bar{u}\, dx 
	+ 2 \Re \int_{\R^d} x\cdot f  \overline{\nabla_{\!A} u}\, dx 
	- 2(\Re \lambda)^{1/2} \sgn(\Im\lambda) \Im \int_{\R^d} \abs{x} f \bar{u}\, dx. 
\end{multline}
		
Recalling definition~\eqref{eq:u^-} of $u^-$, one observes that
\begin{equation}\label{nabla_u^-}
	\nabla_{\!A} u^-(x)= e^{-i (\Re \lambda)^{1/2} \sgn(\Im \lambda) \abs{x}} \left(\nabla_{\!A} u - i (\Re \lambda)^{1/2} \sgn(\Im \lambda) \frac{x}{\abs{x}} u(x)\right),
\end{equation}
and therefore
\begin{equation}\label{eq:square_nabla_u^-}
	\abs{\nabla_{\!A} u^-}^2= \abs{\nabla_{\!A} u}^2 + \Re \lambda \abs{u}^2 - 2 (\Re \lambda)^{1/2} \sgn(\Im\lambda) \frac{x}{\abs{x}} \cdot \Im (\bar{u} \nabla_{\!A} u). 
\end{equation}
Moreover one has
\begin{equation}\label{prop_B}
	x \cdot \bm{B} \cdot u \overline{\nabla_{\!A} u}= x\cdot \bm{B} \cdot u^-\, \overline{\nabla_{\!A} u^-},
\end{equation}
where the previous follows from the fact that being $\bm{B}$ anti-symmetric, then $x\cdot \bm{B} \cdot x=0.$
	
Reintegrating~\eqref{eq:square_nabla_u^-} over~$\R^d$,
we obtain
\begin{equation}\label{eq:nabla_u^-}
	\int_{\R^d} \abs{\nabla_{\!A} u}^2\, dx 
	-2 (\Re\lambda)^{1/2} \sgn(\Im\lambda) \Im \int_{\R^d} \frac{x}{\abs{x}} \cdot \bar{u}\,\nabla_{\!A} u\, dx 
	+ \Re \lambda \int_{\R^d} \abs{u}^2\, dx
	=\int_{\R^d} \abs{\nabla_{\!A} u^-}^2\, dx.
\end{equation}
Adding equation~\eqref{eq:first_x} multiplied by $(\Re \lambda)^{-1/2}\abs{\Im \lambda}$ to~\eqref{eq:intermediate}, plugging~\eqref{eq:nabla_u^-}, using again~\eqref{eq:square_nabla_u^-} and~\eqref{prop_B}, 
we get
\begin{multline}\label{eq:not_yet_minus}
	\int_{\R^d} \abs{\nabla_{\!A} u^-}^2\, dx
	+(\Re \lambda)^{-1/2} \abs{\Im \lambda} \int_{\R^d} \abs{x}\abs{\nabla_{\!A} u^-}^2\, dx
	-\frac{(d-1)}{2}(\Re \lambda)^{-1/2}\abs{\Im \lambda} \int_{\R^d} \frac{\abs{u}^2}{\abs{x}}\,dx\\ 
	+2\Im \int_{\R^d} x \cdot \bm{B} \cdot u^- \overline{\nabla_{\!A} u^-}\, dx\\
	+(d-1)\int_{\R^d} \Re V^{(1)}\abs{u}^2\, dx
	+(\Re \lambda)^{-1/2} \abs{\Im \lambda} \int_{\R^d}\abs{x} \Re V^{(1)}\abs{u}^2\, dx\\
	+2\Re \int_{\R^d} x\cdot V^{(1)} u \left(\overline{\nabla_{\!A} u} + i(\Re \lambda)^{1/2} \sgn(\Im \lambda) \frac{x}{\abs{x}} \bar{u}\right)\, dx\\
	-\int_{\R^d} \partial_r (\abs{x} \Re V^{(2)})\abs{u}^2\, dx
	+(\Re \lambda)^{-1/2} \abs{\Im \lambda}\int_{\R^d}\abs{x} \Re V^{(2)} \abs{u}^2\, dx\\
	-2\Im \int_{\R^d} x \Im V^{(2)} u \left(\overline{\nabla_{\!A} u} + i (\Re \lambda)^{1/2} \sgn(\Im \lambda) \frac{x}{\abs{x}} \bar{u} \right)\, dx\\
	\begin{aligned}
	=&(d-1) \Re \int_{\R^d} f \bar{u}\, dx
	+(\Re \lambda)^{-1/2} \abs{\Im \lambda} \Re \int_{\R^d} \abs{x} f \bar{u}\, dx\\
	&+ 2 \Re \int_{\R^d} x \cdot f  \left( \overline{\nabla_{\!A} u} + i (\Re \lambda)^{1/2} \sgn(\Im \lambda) \frac{x}{\abs{x}} \bar{u} \right)\, dx. \hspace{2.5cm}
	\end{aligned}
\end{multline}	
Then, using~\eqref{nabla_u^-} in the fourth, last but two and last line of the previous identity, we obtain
\begin{multline}\label{eq:last_but_one}
	\int_{\R^d} \abs{\nabla_{\!A} u^-}^2\, dx
	+(\Re \lambda)^{-1/2} \abs{\Im \lambda} \int_{\R^d} \abs{x}\abs{\nabla_{\!A} u^-}^2\, dx
	-\frac{(d-1)}{2}(\Re \lambda)^{-1/2}\abs{\Im \lambda} \int_{\R^d} \frac{\abs{u}^2}{\abs{x}}\,dx\\ 
	+2\Im \int_{\R^d} x \cdot \bm{B} \cdot u^- \overline{\nabla_{\!A} u^-}\, dx\\
	+(d-1)\int_{\R^d} \Re V^{(1)}\abs{u}^2\, dx
	+2\Re \int_{\R^d} x\cdot V^{(1)} u^- \overline{\nabla_{\!A} u^-}\, dx
	+(\Re \lambda)^{-1/2} \abs{\Im \lambda} \int_{\R^d}\abs{x} \Re V^{(1)}\abs{u}^2\, dx\\
	-\int_{\R^d} \partial_r (\abs{x} \Re V^{(2)})\abs{u}^2\, dx
	-2\Im \int_{\R^d} x \Im V^{(2)} u^- \overline{\nabla_{\!A} u^-}\, dx
	+(\Re \lambda)^{-1/2} \abs{\Im \lambda}\int_{\R^d}\abs{x} \Re V^{(2)} \abs{u}^2\, dx\\
	=(d-1) \Re \int_{\R^d} f \bar{u}\, dx
	+ 2 \Re \int_{\R^d} x \cdot f^-  \overline{\nabla_{\!A} u^-}\, dx
	+(\Re \lambda)^{-1/2} \abs{\Im \lambda} \Re \int_{\R^d} \abs{x} f \bar{u}\, dx,
	\end{multline}
	where $f^-(x):=e^{-i(\Re\lambda)^{1/2}\sgn(\Im\lambda)\abs{x}} f(x).$

\item \textbf{Step three.}
Now we want to come back to our approximating sequence $u_R.$ Recalling that $u_R$ is a weak solution to~\eqref{eq:resolvent}, identity~\eqref{eq:last_but_one}, rewritten in terms of $u_R,$ $f_R$ and $\textup{err(R)}$ gives
\begin{multline}\label{eq:last_u_R}
	\int_{\R^d} \abs{\nabla_{\!A} u_R^-}^2\, dx
	+(\Re \lambda)^{-1/2} \abs{\Im \lambda} \int_{\R^d} \abs{x}\abs{\nabla_{\!A} u_R^-}^2\, dx
	-\frac{(d-1)}{2}(\Re \lambda)^{-1/2}\abs{\Im \lambda} \int_{\R^d} \frac{\abs{u_R}^2}{\abs{x}}\,dx\\ 
	+2\Im \int_{\R^d} x \cdot \bm{B} \cdot u_R^- \overline{\nabla_{\!A} u_R^-}\, dx\\
	+(d-1)\int_{\R^d} \Re V^{(1)}\abs{u_R}^2\, dx
	+2\Re \int_{\R^d} x\cdot V^{(1)} u_R^- \overline{\nabla_{\!A} u_R^-}\, dx
	+(\Re \lambda)^{-1/2} \abs{\Im \lambda} \int_{\R^d}\abs{x} \Re V^{(1)}\abs{u_R}^2\, dx\\
	-\int_{\R^d} \partial_r (\abs{x} \Re V^{(2)})\abs{u_R}^2\, dx
	-2\Im \int_{\R^d} x \Im V^{(2)} u_R^- \overline{\nabla_{\!A} u_R^-}\, dx
	+(\Re \lambda)^{-1/2} \abs{\Im \lambda}\int_{\R^d}\abs{x} \Re V^{(2)} \abs{u_R}^2\, dx\\
	=(d-1) \Re \int_{\R^d} f_R \overline{u_R}\, dx
	+ 2 \Re \int_{\R^d} x \cdot f_R^-  \overline{\nabla_{\!A} u_R^-}\, dx
	+(\Re \lambda)^{-1/2} \abs{\Im \lambda} \Re \int_{\R^d} \abs{x} f_R \overline{u_R}\, dx\\
	+(d-1) \Re \int_{\R^d} \textup{err}(R) \overline{u_R}\, dx
	+ 2 \Re \int_{\R^d} x \cdot \textup{err}(R)^-  \overline{\nabla_{\!A} u_R^-}\, dx
	+(\Re \lambda)^{-1/2} \abs{\Im \lambda} \Re \int_{\R^d} \abs{x} \textup{err}(R) \overline{u_R}\, dx.
	\end{multline} 

Letting $R$ go to infinity, the thesis follows from dominated and monotone convergence theorems and Lemma~\ref{lemma:f}.  
\end{itemize}

	\subsection{The method of multipliers: proof of the crucial Lemma~\ref{lemma:crucial_identities}}\label{subsec:proof_identities}
	
	This part is entirely devoted to the rigorous proof of the crucial identities contained in Lemma~\ref{lemma:crucial_identities}.
Let us start proving~\eqref{eq:first_const} and~\eqref{eq:first_x}. 
Choosing in~\eqref{eq:weak_resolvent} $v:= \varphi u,$ with $\varphi\colon \R^d \to \R$ being a radial function such that $v\in \mathcal{D}_{A,V}$ (since the support of $u$ is compact, any locally bounded $\varphi$ together with locally bounded partial derivatives of first order is admissible). Using the generalised Leibniz rule for the magnetic gradient, namely
	\begin{equation}\label{eq:chain_rule_magnetic}
		\nabla_{\!A}(g h)= (\nabla_{\!A} g) h + g \nabla h
	\end{equation}
	valid for any $g, h\colon \R^d \to \C,$
we get 
\begin{equation*}
	\int_{\R^d} \varphi \abs{\nabla_{\!A} u}^2
	 + \int_{\R^d} \bar{u}\nabla_{\!A} u \cdot \nabla \varphi
	  + \int_{\R^d} V \varphi \abs{u}^2
	  = \lambda \int_{\R^d} \varphi \abs{u}^2 + \int_{\R^d} f \varphi \bar{u}.
\end{equation*}
Taking the real part of the obtained identity,
using that being~$A$ a real-valued vector field one has
one has
\begin{equation}\label{eq:real_magnetic_A}
	\Re(\bar{u}\nabla_{\!A} u)= \Re(\bar{u}\nabla u) 
\end{equation}
 and performing an integration by parts give
\begin{equation*}
	-\frac{1}{2} \int_{\R^d} \Delta \varphi \abs{u}^2 + \int_{\R^d} \varphi \abs{\nabla_{\!A} u}^2 + \int_{\R^d} \Re V \varphi \abs{u}^2= \Re \lambda \int_{\R^d} \varphi \abs{u}^2 + \Re \int_{\R^d} f \varphi \bar{u}.  
\end{equation*}
Taking $\varphi:=1$ and $\varphi(x):=\abs{x},$ we get~\eqref{eq:first_const} and~\eqref{eq:first_x}.
Equation~\eqref{eq:second_const} and~\eqref{eq:second_x} are obtained as in the previous case choosing in~\eqref{eq:weak_resolvent} $v:=\psi u,$ with $\psi\colon \R^d\to \R$ being a radial function such that $v\in \mathcal{D}_{A,V}$ and taking the imaginary part of the resulting identity. Finally, one chooses $\psi:=1$ and $\psi(x):=\abs{x},$ respectively.

The remaining identity~\eqref{eq:third} is formally obtained by plugging into~\eqref{eq:weak_resolvent} the multiplier
\begin{equation}\label{Morawetz_multiplier}
	v:=[\nabla_{\!A}^2,\phi]u=\Delta \phi u + 2 \nabla \phi\cdot \nabla_{\!A} u \qquad \text{with}\qquad \phi(x):=\abs{x}^2,
\end{equation}
taking the real part and integrating by parts. However, such $v$ does not need to belong to $\mathcal{D}_A$ (and therefore neither to $\mathcal{D}_{A,V}$).
Indeed, though on the one hand the unboundedness of the radial function~$\phi$ does not pose any problems because the support of $u$ is assumed to be compact at this step, on the other hand $\nabla_{\!A} u$ does not necessarily belong to $\mathcal{D}_A.$ Following the strategy developed in~\cite{Cossetti-Krejcirik_2018}, we replace~\eqref{Morawetz_multiplier} by its regularised version
\begin{equation}\label{Morawetz_multiplier_reg}
	v:= \Delta \phi u + \nabla \phi \cdot [\nabla_{\!A}^{\delta, N} u + \nabla_{\!A}^{-\delta,N}u]=\Delta \phi u + \partial_k \phi \,[\partial_{k,A}^{\delta,N} u + \partial_{k,A}^{-\delta,N}u]
	\qquad \text{with}\qquad 
	\phi(x):=\abs{x}^2,
\end{equation}
where 
\begin{equation}\label{eq:regularised_magnetic_derivative}
	\nabla_{\!A}^{\delta,N} u
	:=(\partial_{1,A}^{\delta,N}, \dots, \partial_{d,A}^{\delta,N})u
	\qquad  \text{with}\qquad 
	\partial_{k,A}^{\delta,N}u:=\partial_k^\delta u + i T_N(A_k)u, 
	\qquad k=1,2,\dots,d,
\end{equation} 
and where
\begin{equation}\label{eq:difference_quotient}
	\partial_k^\delta u(x)
  :=\frac{\tau_k^\delta u(x)-u(x)}{\delta}
  \qquad \mbox{with} \qquad
  \tau_k^\delta u(x) := u(x+ \delta e_k),
	\qquad k=1,2,\dots,d,
\end{equation}
with $\delta\in \R \setminus \{0\}$ is the standard \emph{difference quotient} of $u$ (we refer to \cite[Sec.~5.8.2]{Evans} or~\cite[Sec.~10.5]{Leoni}
for basic facts about the difference quotients) and the Lipschitz continuous function
\begin{equation*}
	T_N(s):=\max\{-N, \min\{s,N\}\}
\end{equation*}
with $N>0$ is the usual \emph{truncation function}.
After the second equality of~\eqref{Morawetz_multiplier_reg} and in the sequel, we use the Einstein summation convention.

We start showing that $v$ defined as in~\eqref{Morawetz_multiplier_reg} belongs to $\mathcal{D}_{A,V},$ which is saying $v\in L^2(\R^d),$ $\partial_{l,A}v:=(\partial_l+iA_l)v\in L^2(\R^d)$ for any $l=1, \dots, d$ and $\sqrt{(\Re V^{(2)})_+}v\in L^2(\R^d).$ To see that, let us rewrite explicitly~\eqref{Morawetz_multiplier_reg} with the choice $\phi(x):=\abs{x}^2,$ that is
\begin{equation}\label{eq:multiplier_explicit}
	v:=2d u + 2x_k[\partial_{k,A}^{\delta, N} u + \partial_{k,A}^{-\delta,N} u].
\end{equation}
Clearly, being $u\in \mathcal{D}_{A,V},$ the first term in $v$ belongs to $\mathcal{D}_{A,V}$ and therefore we need to comment further just on the second term of the sum, namely $x_k\,\partial_{k,A}^{\delta,N} u$ (the part involving $\partial_{k,A}^{-\delta,N}u$ is analogous). One can check that $x_k \partial_{k,A}^{\delta, N}u:=x_k(\partial_k^\delta + i T_N(A_k))u\in L^2(\R^d)$; this is a consequence of $u\in L^2(\R^d)$ being compactly supported and of the boundedness of $T_N(A_k).$ It is less trivial to prove that for any $l=1,2,\dots,d,$ one has $\partial_{l,A}[x_k \partial_{k,A}^{\delta, N}u]\in L^2(\R^d).$

To begin with, it is easy to check that the following commutation relation between the magnetic gradient $\partial_{l,A}$ and its \emph{regularised} version $\partial_{k,A}^{\delta, N}$ holds true
\begin{equation}\label{eq:commutation_relation}
	\big[\partial_{l,A}, \partial_{k,A}^{\delta, N}\big]:=i [(\partial_l A_k) \chi_{\{\abs{A_k}\leq N\}} - (\partial_{k}^\delta A_l) \tau_k^\delta], \qquad k,l=1,2\dots,d.
\end{equation}
Here $[\cdot, \cdot]$ denotes the usual commutator operator, for any given subset $S\subseteq \R^d,$ the function $\chi_S$ is the characteristic function of the set $S$ and $\tau_k^\delta$ is the translation operator as defined in~\eqref{eq:difference_quotient}.

Using~\eqref{eq:chain_rule_magnetic}, the fact that, by definition of the commutator operator, $\partial_{l,A} \partial_{k,A}^{\delta, N}= \partial_{k,A}^{\delta, N} \partial_{l,A} + [\partial_{l,A} \partial_{k,A}^{\delta, N}]$ and eventually using~\eqref{eq:commutation_relation} one has
\begin{equation*}
	\begin{split}
		\partial_{l,A}[x_k \partial_{k,A}^{\delta, N}u]&=\delta_{l,k} \partial_{k,A}^{\delta,N}u + x_k \partial_{l,A}\partial_{k,A}^{\delta,N}u\\
		&=\delta_{l,k} \partial_{k,A}^{\delta,N}u + x_k \partial_{k,A}^{\delta,N}\partial_{l,A}u + x_k [\partial_{l,A}, \partial_{k,A}^{\delta, N}]u\\
 &=v_1+v_2+v_3,
	\end{split}
\end{equation*}
where
\begin{equation*}
	v_1:=	\delta_{l,k} \partial_{k,A}^{\delta,N}u,
	\qquad
	v_2:=
	x_k \partial_{k,A}^{\delta,N}\partial_{l,A}u,
	\qquad
	v_3:=
	x_ki [(\partial_l A_k) \chi_{\{\abs{A_k}\leq N\}} - (\partial_{k}^\delta A_l) \tau_k^\delta]u	.
\end{equation*}
Here and hence $\delta_{l,k}$ for every $k,l=1,2,\dots, d$ denotes the Kronecker symbol.

Now, being $u\in \mathcal{D}_{A,V}$ (thus in particular $u\in L^2(\R^d)$) and since $T_N(A_k)\in L^\infty(\R^d),$ then
\begin{equation*}
	v_1= \delta_{l,k} \partial_{k,A}^{\delta, N}u:=\delta_{l,k} (\partial_k^\delta + i T_N(A_k))u
\end{equation*}
is clearly in $L^2(\R^d).$
Moreover, since $u\in \mathcal{D}_{A,V}$ (thus in particular $\partial_{l,A} u\in L^2(\R^d)$) is compactly supported, one can conclude the same for $v_2.$ With respect to $v_3,$ since $A_k\in W^{1,p}_\textup{loc}(\R^d)$ with $p$ as in~\eqref{eq:hypotheses_A}, then $(\partial_l A_k) u\in L^2(\R^d)$ (see~\eqref{Au_L2}). Similar reasoning allows us to conclude that also $(\partial_k^\delta A_l) \tau_k^\delta u\in L^2(\R^d).$ Therefore $v_3\in L^2(\R^d).$

Now we are left to show just that $\sqrt{(\Re V^{(2)})_+} [x_k \partial_{k,A}^{\delta, N} u]\in L^2(\R^d).$
First let us write
\begin{equation*}
	\sqrt{(\Re V^{(2)})_+} [x_k \partial_{k,A}^{\delta, N} u]
	=v_4+v_5,
\end{equation*}
where
\begin{equation*}
	v_4:=
	x_k \sqrt{(\Re V^{(2)})_+} \partial_k^\delta u,
	\qquad
	v_5:=  i x_k \sqrt{(\Re V^{(2)})_+} T_N(A_k) u.
\end{equation*}
Observe that being $u\in \mathcal{D}_{A,V}$ (thus in particular $\sqrt{(\Re V^{(2)})_+} u \in L^2(\R^d)$) and compactly supported and since $T_N(A_k)\in L^\infty(\R^d),$  one has that $v_5\in L^2(\R^d).$ Making explicit the difference quotient $\partial_k^\delta u,$ one can also see that $v_4\in L^2(\R^d)$ by using that $(\Re V^{(2)})_+\in L^p_\textup{loc}(\R^d)$ with $p$ as in~\eqref{eq:hypotheses_V_alt} and the fact that $\abs{u}\in H^1(\R^d).$ 
  
Gathering these facts together, we guaranteed that our multiplier $v$ as defined in~\eqref{eq:multiplier_explicit} belongs to $\mathcal{D}_{A,V}$ and hence we have justified its choice as a test function in the weak formulation~\eqref{eq:weak_resolvent}.

Now we are in a position to prove identity~\eqref{eq:third}. For a moment, we proceed in a greater generality by considering $\phi$ in~\eqref{Morawetz_multiplier_reg} to be an arbitrary smooth function $\phi\colon \R^d\to \R.$ We plug~\eqref{Morawetz_multiplier_reg} in~\eqref{eq:weak_resolvent} and take the real part. Below, for the sake of clarity, we consider each integral of the resulting identity separately.

\subsubsection*{$\bullet$ Kinetic term}
Let us start with the ``kinetic'' part of~\eqref{eq:weak_resolvent}:
\begin{equation}\label{eq:kinetic_part}
	K:=\Re \int_{\R^d} \nabla_{\!A}u \cdot \overline{\nabla_{\!A} v}.
\end{equation} 
Using
\begin{equation*}
	\overline{\partial_{l,A}v}=(\partial_l \Delta \phi) \bar{u} + \Delta \phi \overline{\partial_{l,A}u} + \partial_{lk}\phi \,[\overline{\partial_{k,A}^{\delta, N} u} + \overline{\partial_{k,A}^{-\delta, N} u}] + \partial_k \phi \,[ \overline{\partial_{l,A} \partial_{k,A}^{\delta, N} u} + \overline{\partial_{l,A} \partial_{k,A}^{-\delta, N} u}],
\end{equation*}
we write $K=K_1+K_2+K_3+K_4$ with
\begin{equation}\label{eq:K's}
\begin{aligned}
K_1&:=\Re \int_{\R^d} \partial_{l,A} u (\partial_l \Delta \phi) \bar{u}, 
&
K_2&:=\int_{\R^d} \abs{\nabla_{\!A} u}^2 \Delta \phi,
\\
K_3&:=\Re \int_{\R^d} \partial_{lk}\phi\, \partial_{l,A} u\, [\overline{\partial_{k,A}^{\delta, N} u} + \overline{\partial_{k,A}^{-\delta, N} u}], 
& 
K_4&:=\Re \int_{\R^d} \partial_k \phi \, \partial_{l,A} u\, [\overline{\partial_{l,A} \partial_{k,A}^{\delta, N} u} + \overline{\partial_{l,A} \partial_{k,A}^{-\delta, N} u}].
\end{aligned}
\end{equation}

Using~\eqref{eq:real_magnetic_A} and integrating by parts in $K_1$  give
\begin{equation*}
	K_1=-\frac{1}{2}\int_{\R^d} \Delta^2 \phi \abs{u}^2.
\end{equation*}
Now we consider $K_4.$ Using simply the definition of the commutator operator,
 we write
\begin{equation*}
	K_4=
	K_{4,1} + K_{4,2},
\end{equation*}
where 
\begin{equation*}
K_{4,1}:=\Re \int_{\R^d} \partial_k \phi \, \partial_{l,A} u\, \big\{\overline{\partial_{k,A}^{\delta, N} \partial_{l,A} u} + \overline{\partial_{k,A}^{-\delta, N} \partial_{l,A} u}\big\},
\qquad
K_{4,2}:=
\Re \int_{\R^d} \partial_k \phi \partial_{l,A} u \big\{ \overline{[\partial_{l,A}, \partial_{k,A}^{\delta,N}] u} + \overline{[\partial_{l,A}, \partial_{k,A}^{-\delta,N}] u}\big\}.
\end{equation*}
We start considering $K_{4,1}.$ Using an analogous version to~\eqref{eq:real_magnetic_A} for the regularised magnetic gradient, namely
\begin{equation}\label{eq:real_magnetic_A_reg}
	\Re (\bar{u}\, \partial_{k,A}^{\delta, N} u)= \Re(\bar{u}\, \partial_k^\delta u),\qquad k=1,2,\dots, d
\end{equation} 
and the identity
\begin{equation}\label{eq:real}
  2\Re(\bar{\psi} \partial_k^\delta \psi)
  =\partial_k^\delta \abs{\psi}^2 -\delta \abs{\partial_k^\delta \psi}^2
\end{equation}
valid for every $\psi:\R^d \to \C$,
we write $K_{4,1}=K_{4,1,1} + K_{4,1,2}$ with
\begin{equation*}
	K_{4,1,1}:=\frac{1}{2}\int_{\R^d} \partial_k \phi \{\partial_k^\delta \abs{\partial_{l,A}u}^2 + \partial_{k}^{-\delta} \abs{\partial_{l,A} u}^2\},
	\quad \text{and}\quad
	K_{4,1,2}:=-\frac{\delta}{2} \int_{\R^d} \partial_k \phi \{ \abs{\partial_k^\delta \partial_{l,A} u}^2 - \abs{\partial_k^{-\delta} \partial_{l,A} u}^2 \}.	
\end{equation*}
Making use of the integration-by-parts formula for difference quotients (see \cite[Sec.~5.8.2]{Evans})
\begin{equation}\label{int_by_parts}
  \int_{\R^d} \varphi \ \partial_k^{\delta}\psi
  = - \int_{\R^d} (\partial_k^{-\delta} \varphi) \ \psi
\end{equation}
which holds true for every $\varphi,\psi \in L^2(\R^d),$ one gets
\begin{equation*}
	K_{4,1,1}=-\frac{1}{2} \int_{\R^d} \big \{ \partial_k^{-\delta} \partial_k \phi + \partial_k^\delta \partial_k \phi \big\} \abs{\nabla_{\!A} u}^2.
\end{equation*}
At the same time, making explicit the difference quotient and changing variable in $K_{4,1,2}$ give (summation both over $k$ and $l$)
\begin{equation*}
	K_{4,1,2}=-\frac{\delta}{2} \int_{\R^d} \{\partial_k \phi - (\tau_{k}^\delta \partial_k \phi)\} \abs{\partial_k^{\delta} \partial_{l,A} u}^2. 
\end{equation*}
Now we choose the multiplier $\phi(x):=\abs{x}^2$ and observe that
\begin{equation}\label{eq:phi_derivatives}
	\partial_k \phi=2x_k, \qquad \partial_{lk} \phi=2\delta_{k,l}, \qquad \partial_k^{\pm\delta}\partial_k\phi=2,\qquad \nabla \Delta \phi=0, \qquad \Delta^2\phi=0.
\end{equation}
Consequently,
\begin{equation*}
	K_1=0,\qquad
	K_2=2d \int_{\R^d} \abs{\nabla_{\!A} u}^2\, dx,\qquad
	K_3=2 \Re \int_{\R^d} \partial_{l,A} u\,[\overline{\partial_{l,A}^{\delta,N} u} + \overline{\partial_{l,A}^{-\delta, N} u}]\, dx,
\end{equation*}
and
\begin{multline*}
K_{4}=-2d\int_{\R^d} \abs{\nabla_{\!A} u}^2\, dx
+\int_{\R^d} \abs{\tau_k^\delta \nabla_{\!A} u -\nabla_{\!A} u}^2\, dx\\
+ 2 \Im \int_{\R^d} x_k \partial_{l,A}u \big[(\partial_l A_k)\chi_{\{\abs{A_k}\leq N\}}\bar{u} - (\partial_k^\delta A_l) \tau_k^\delta \bar{u} \big]\, dx + 2\Im \int_{\R^d} x_k \partial_{l,A}u \big[(\partial_l A_k)\chi_{\{\abs{A_k}\leq N\}}\bar{u} - (\partial_k^{-\delta} A_l) \tau_k^{-\delta} \bar{u} \big]\, dx.
\end{multline*}
In summary,
\begin{multline*}
	K=2\Re \int_{\R^d} \partial_{l,A}u\,[\overline{\partial_{l,A}^{\delta,N} u} + \overline{\partial_{l,A}^{-\delta, N} u}]\, dx
		+\int_{\R^d} \abs{\tau_k^\delta \nabla_{\!A} u -\nabla_{\!A} u}^2\, dx\\
+ 2 \Im \int_{\R^d} x_k \partial_{l,A}u \big[(\partial_l A_k)\chi_{\{\abs{A_k}\leq N\}}\bar{u} - (\partial_k^\delta A_l) \tau_k^\delta \bar{u} \big]\, dx + 2\Im \int_{\R^d} x_k \partial_{l,A}u \big[(\partial_l A_k)\chi_{\{\abs{A_k}\leq N\}}\bar{u} - (\partial_k^{-\delta} A_l) \tau_k^{-\delta} \bar{u} \big]\, dx.
\end{multline*}

Now we want to see what happens when $\delta$ goes to zero and $N$ goes to infinity. To do so, we need the following lemma.
\begin{lemma}\label{lemma:auxiliary}
Under the hypotheses of Lemma~\ref{lemma:crucial_identities},
the following limits hold true:
\begin{equation}\label{first_conv}
	\partial_{l,A}^{\delta, N}u \xrightarrow{\substack{\delta\to 0\\N\to \infty}} \partial_{l,A}u \qquad \text{in}\quad L^2(\R^d)
\end{equation}
and
\begin{equation}\label{second_conv}
	\big[(\partial_l A_k) \chi_{\{\abs{A_k}\leq N\}} - (\partial_k^\delta A_l)\tau_k^\delta\big]u \xrightarrow{\substack{\delta\to 0\\N\to \infty}} [\partial_l A_k - \partial_k A_l]u \qquad \text{in}\quad L^2(\R^d).
\end{equation}
\end{lemma}
\begin{proof}
Let us start with~\eqref{first_conv}. Using the explicit expression~\eqref{eq:regularised_magnetic_derivative} for $\partial_{l,A}^{\delta, N} u,$ one easily has
\begin{equation*}
	\int_{\R^d} \abs{\partial_{l,A}^{\delta, N} u - \partial_{l,A} u}^2\, dx
	\leq 
	2 \int_{\R^d} \abs{\partial_{l}^{\delta} u - \partial_{l} u}^2\, dx + 2 \int_{\R^d} \abs{T_N(A_l) u - A_l u}^2\, dx.
\end{equation*}
Now, as a consequence of the $L^2$-strong convergence of the difference quotients (which can be used here because $u\in H^1(\R^d)$ (see~\eqref{inclusion})), the first integral converges to zero as $\delta$ goes to zero. As regards with the second integral we use that, by definition, $T_N(s)$ converges to $s$ as $N$ tends to infinity, the bound $\abs{T_N(s)}\leq \abs{s}$ and the fact that by virtue of~\eqref{eq:hypotheses_A} the function $A_l u\in L^2(\R^d),$ these allow us to conclude that the integral goes to zero as $N$ goes to infinity via the dominated convergence theorem. This concludes the proof of~\eqref{first_conv}.

Now we prove~\eqref{second_conv}. Observe that~\eqref{second_conv} follows as soon as one proves that the limits
\begin{equation*}
	(\partial_l A_k) \chi_{\{\abs{A_k}\leq N\}} u \xrightarrow{N\to \infty} \partial_l A_k u\qquad \text{in}\quad L^2(\R^d)
\end{equation*}
and 
\begin{equation*}
	(\partial_k^\delta A_l) \tau_k^\delta u \xrightarrow{\delta \to 0} \partial_k A_l u\qquad \text{in}\quad L^2(\R^d)
\end{equation*}
hold true.
As hypothesis~\eqref{eq:hypotheses_A} implies that $\partial_l A_k u \in L^2(\R^d),$ the first limit is an immediate consequence of the dominated convergence theorem. 
With respect to the second one, one has
\begin{equation*}
	\int_{\R^d} \abs{(\partial_k^\delta A_l) \tau_k^\delta u - \partial_k A_l u}^2\, dx
	\leq 
	2\int_{\R^d} \abs{\partial_k^\delta A_l}^2 \abs{\tau_k^\delta u-u}^2\, dx
	+2\int_{\R^d} \abs{\partial_k^\delta A_l - \partial_k A_l}^2\abs{u}^2\, dx
\end{equation*}
and the two integrals tend to zero as $\delta$ goes to zero as a consequence of the $L^q$-continuity of the translations with $1\leq q<\infty$ 
and the strong $L^p$-convergence of the difference quotients with $1\leq p<\infty$ together with assumption~\eqref{eq:hypotheses_A}.
\end{proof}

With Lemma~\ref{lemma:auxiliary} at hand, it follows as a mere consequence of the Cauchy--Schwarz inequality that 
\begin{equation*}
	K\xrightarrow{\substack{\delta \to 0\\ N\to \infty}} 4\int_{\R^d} \abs{\partial_{l,A} u}^2\, dx + 4\Im \int_{\R^d} x_k \partial_{l,A} u\, [\partial_l A_k - \partial_k A_l] \bar{u}\, dx.
\end{equation*} 

\subsubsection*{$\bullet$ Source term}
Let us now consider simultaneously
the ``source'' and ``eigenvalue''
parts of~\eqref{eq:weak_resolvent}, that is,
\begin{equation}\label{eq:source_part}
	F:=\Re \left( \lambda \int_{\R^d} u \bar{v}  + \int_{\R^d} f \bar{v} \right).
\end{equation}
This can be written as $F=F_1+F_2+F_3+F_4$ with
\begin{equation}\label{F's}
	\begin{aligned}
	F_1&:=\Re\lambda \int_{\R^d} \Delta \phi \abs{u}^2,
&
	F_2&:=\Re \lambda \Re \int_{\R^d} \partial_k \phi \, u [\overline{\partial_{k, A}^{\delta,N} u} + \overline{\partial_{k,A}^{-\delta,N} u}],
\\
	F_3&:=-\Im \lambda \Im \int_{\R^d} \partial_k \phi \, u [\overline{\partial_{k, A}^{\delta,N} u} + \overline{\partial_{k,A}^{-\delta,N} u}],
&
	F_4&:= \Re \int_{\R^d} f\{\Delta \phi \bar{u} + \partial_k \phi \, [\overline{\partial_{k, A}^{\delta,N} u} + \overline{\partial_{k,A}^{-\delta,N} u}]\}.
	\end{aligned}
\end{equation}
Applying~\eqref{eq:real_magnetic_A_reg} and~\eqref{eq:real}, we further split $F_2=F_{2,1} + F_{2,2},$ where
\begin{equation*}
	F_{2,1}:=\frac{1}{2} \Re \lambda \int_{\R^d} \partial_k \phi\, \{\partial_k^\delta \abs{u}^2 + \partial_k^{-\delta} \abs{u}^2 \}
	\quad \text{and}\quad
	F_{2,2}:=-\frac{\delta}{2} \Re \lambda \int_{\R^d} \partial_k \phi\, \{ \abs{\partial_k^\delta u}^2 - \abs{\partial_k^{-\delta} u}^2 \}.
\end{equation*}
Using the integration-by-parts formula~\eqref{int_by_parts}, we get
\begin{equation*}
	F_{2,1}=-\frac{1}{2}\Re\lambda \int_{\R^d} \{\partial_k^{-\delta} \partial_k \phi + \partial_k^{\delta} \partial_k \phi\}\abs{u}^2.
\end{equation*}
Choosing $\phi(x):=\abs{x}^2$ 
in the previous identities and using~\eqref{eq:phi_derivatives} gives
\begin{equation*}
	\begin{aligned}
	& F_1=2d \Re\lambda \int_{\R^d} \abs{u}^2\, dx,
	&
	& F_2=-2d \Re \lambda \int_{\R^d} \abs{u}^2\, dx 
	       -\delta \Re \lambda \int_{\R^d} x_k \{\abs{\partial_k^\delta u}^2 - \abs{\partial_k^{-\delta} u}^2\}\, dx,\\ 
	&F_3=-2\Im \lambda \Im \int_{\R^d} x_k u\, [\overline{\partial_{k,A}^{\delta,N} u} + \overline{\partial_{k,A}^{-\delta,N} u}]\, dx,
	&
	&F_4= \Re \int_{\R^d} f \{2d \bar{u} +2 x_k [\overline{\partial_{k,A}^{\delta,N} u} + \overline{\partial_{k,A}^{-\delta,N} u}]\}\, dx.
\end{aligned}
\end{equation*}
Using limit~\eqref{first_conv} in Lemma~\ref{lemma:auxiliary}, one gets from the Cauchy--Schwarz inequality that
\begin{equation*}
	F\xrightarrow{\substack{\delta \to 0\\ N\to \infty}} -4\Im \lambda \Im \int_{\R^d} x_k u\, \overline{\partial_{k,A} u}\, dx
		+ \Re \int_{\R^d} f \{2d\bar{u} + 4 x_k \, \overline{\partial_{k,A} u}\}\, dx.
\end{equation*}

\subsubsection*{$\bullet$ Electric potential term}
Let us now consider the contribution 
of the ``potential'' part of~\eqref{eq:weak_resolvent}, that is, 
\begin{equation}\label{eq:potential_part}
	J:=\Re \int_{\R^d} V u \bar{v}.
\end{equation}
Using the decomposition $V=V^{(1)} + V^{(2)},$ it can be written as $J=J_1+J_2$ with
\begin{equation*}
	J_1:=\Re \int_{\R^d} V^{(1)} u \bar{v}
	\qquad \text{and}\qquad
	J_2:=\Re \int_{\R^d} V^{(2)} u \bar{v}
\end{equation*}
First of all,
\begin{equation*}
	J_1= \int_{\R^d} \Re V^{(1)} \Delta \phi \abs{u}^2 + \Re \int_{\R^d} \partial_k \phi V^{(1)} \,  u\, [\overline{\partial_{k,A}^{\delta,N} u} + \overline{\partial_{k,A}^{-\delta,N}u}]. 
\end{equation*}
Let us consider now the part involving $V^{(2)}.$ We can write
\begin{equation*}
	J_2
	=J_{2,1}+ J_{2,2}+ J_{2,3},
\end{equation*}
where
\begin{gather*}
	J_{2,1}:= \int_{\R^d} \Re V^{(2)} \Delta \phi \abs{u}^2,
	\qquad
	J_{2,2}:= 
	\int_{\R^d} \Re V^{(2)} \partial_k\phi\, \Re\{u\,[\overline{\partial_{k,A}^{\delta,N} u} + \overline{\partial_{k,A}^{-\delta,N}u}]\}\\
	J_{2,3}:=
	- \Im \int_{\R^d} \Im V^{(2)} \partial_k\phi\, u\,[\overline{\partial_{k,A}^{\delta,N} u} + \overline{\partial_{k,A}^{-\delta,N}u}]
\end{gather*}
Let us consider $J_{2,2}.$ Using~\eqref{eq:real_magnetic_A_reg},~\eqref{eq:real} and integrating by parts we get
\begin{equation*}
	\begin{split}
	J_{2,2}
	&=-\frac{1}{2} \int_{\R^d} \{\partial _k^{-\delta}[\partial_k \phi \Re V^{(2)}] + \partial_k^\delta [\partial_k \phi \Re V^{(2)}]\} \abs{u}^2 
	- \frac{\delta}{2} \int_{\R^d} \Re V^{(2)} \partial_k \phi \{\abs{\partial_k^\delta u}^2 -\abs{\partial_k^{-\delta} u}^2 \}\\
	&=
	-\frac{1}{2} \int_{\R^d} \{\partial _k^{-\delta}[\partial_k \phi \Re V^{(2)}] + \partial_k^\delta [\partial_k \phi \Re V^{(2)}]\} \abs{u}^2
	+ \frac{1}{2} \int_{\R^d} \partial_k^\delta [\partial_k \phi \Re V^{(2)}] \abs{\tau_k^\delta u -u}^2.
	\end{split}
\end{equation*}
Choosing $\phi(x):=\abs{x}^2$ in the previous identities and using~\eqref{eq:phi_derivatives} we can write
\begin{equation*}
	J_1
	=J_{1,1} +J_{1,2},		
\end{equation*}
where	
\begin{equation*}
	J_{1,1}:=2d\int_{\R^d} \Re V^{(1)} \abs{u}^2\, dx
	\quad \text{and} \quad
	J_{1,2}:=2 \Re \int_{\R^d} x_k V^{(1)} \,  u\, [\overline{\partial_{k,A}^{\delta,N} u} + \overline{\partial_{k,A}^{-\delta,N}u}]\, dx.
\end{equation*}
Moreover
\begin{equation*}
	J_{2,1}= 2d \int_{\R^d} \Re V^{(2)} \abs{u}^2\, dx,
	\qquad 
	J_{2,3}=-2\Im \int_{\R^d} x_k \Im V^{(2)} u [\overline{\partial_{k,A}^{\delta, N}u} + \overline{\partial_{k,A}^{-\delta, N}u}]\, dx,
\end{equation*}
\begin{equation*}
	J_{2,2}=
	- \int_{\R^d} \{\partial _k^{-\delta}[x_k \Re V^{(2)}] + \partial_k^\delta [x_k \Re V^{(2)}]\} \abs{u}^2\, dx
	+ \int_{\R^d} \partial_k^\delta [x_k \Re V^{(2)}] \abs{\tau_k^\delta u -u}^2\, dx.	
\end{equation*}
By virtue of hypothesis~\eqref{x_hypothesis}, $\abs{x}\abs{V^{(1)}} \abs{u}\in L^2_\textup{loc}(\R^d)$ and then, using the Cauchy--Schwarz inequality and limit~\eqref{first_conv} in Lemma~\ref{lemma:auxiliary}, one has
\begin{equation}\label{eq:J_2_limit}
	J_{1,2}\xrightarrow{\substack{\delta \to 0\\ N\to \infty}} 4\Re \int_{\R^d} x_k V^{(1)} u\, \overline{\partial_{k,A} u}\, dx.
\end{equation}

Similarly, using that $\abs{x}\abs{\Im V^{(2)}}\abs{u}\in L^2_\textup{loc}(\R^d)$ (see~\eqref{x_hypothesis}) and again~\eqref{first_conv}, via the Cauchy--Schwarz inequality one also has
\begin{equation*}
J_{2,3}\xrightarrow{\substack{\delta \to 0\\ N\to \infty}} -4\Im \int_{\R^d} x_k \Im V^{(2)} u\, \overline{\partial_{k,A} u}\, dx.
\end{equation*}
Since $x_k \Re V^{(2)}\in W_\textup{loc}^{1,p}(\R^d)$ with $p$ as in~\eqref{eq:hypotheses_V_alt}, using the strong $L^p$-convergence of the difference quotients with $1\leq p<\infty$ and via the H\"older inequality, it is not difficult to see that 
\begin{equation*}
	\begin{split}
		J_{2,2}\xrightarrow{\delta \to 0} & -2 \int_{\R^d} \partial_k [x_k \Re V^{(2)}]\abs{u}^2\, dx\\
		& =-2d \int_{\R^d} \Re V^{(2)} \abs{u}^2\, dx -2\int_{\R^d} x_k \partial_k \Re V^{(2)} \abs{u}^2\, dx,
	\end{split}
\end{equation*}
where the last identity follows from the Leibniz rule applied to $\partial_k(x_k \Re V^{(2)}).$

In summary, 
gathering the previous limits altogether, one gets
\begin{equation*}
	J_1\xrightarrow{\substack{\delta \to 0\\ N\to \infty}} 2d\int_{\R^d} \Re V^{(1)} \abs{u}^2\, dx+4 \Re \int_{\R^d} x_k V^{(1)} u\, \overline{\partial_{k,A} u}\, dx.
\end{equation*}
and 
\begin{equation*}
	J_2\xrightarrow{\substack{\delta \to 0\\ N\to \infty}} 
	-2 \int_{\R^d} x_k \partial_k \Re V^{(2)} \abs{u}^2\, dx 
	- 4 \Im \int_{\R^d} x_k \Im V^{(2)} u \overline{\partial_{k,A}u}\, dx.
\end{equation*}

Passing to the limit $\delta \to 0$ and $N\to \infty$ in~\eqref{eq:weak_resolvent} and multiplying the resulting identity by $1/2,$ one obtains~\eqref{eq:third}.

\qedhere

\subsection{Potentials with just one singularity: 
alternative proof of the crucial Lemma~\ref{lemma:crucial_identities}}\label{Sec:one_singularity}

In this section we consider the case of potentials (both electric and magnetic) with capacity zero set of singularities, in fact with just one singularity at the origin. This will allow us to remove 
the unpleasant hypotheses~\eqref{eq:hypotheses_A} 
and~\eqref{eq:hypotheses_V_alt}.  
Since the point has a positive capacity in dimension one,
here we exclusively consider $d\geq 2$.
(As a matter of fact, if $d=1$,
hypothesis~\eqref{eq:hypotheses_V_alt} is rather natural,
while~\eqref{eq:hypotheses_A} is automatically satisfied 
because of the absence of magnetic fields on the real line.)

To be more specific, in the sequel we consider the following setup.
Let $A\in L^{2}_\textup{loc}(\R^d\setminus\{0\};\R^d)$ and $V\in L^1_\textup{loc}(\R^d\setminus \{0\};\C)$ and assume 
\begin{equation}\label{eq:con_AB}
	\Re V\in L^\infty_\textup{loc}(\R^d\setminus \{0\}) 
\qquad \text{and}\qquad A\in W^{1,\infty}_\textup{loc}(\R^d\setminus \{0\}).
\end{equation}
Notice that assumption~\eqref{eq:con_AB} is satisfied by a large class of potentials, namely $V(x)=a/\abs{x}^\alpha$ with $a>0$ and $\alpha>0$ and the Aharonov--Bohm vector field~\eqref{eq:A-B}.

Observe that since it is no more necessarily true that $V\in L^1_\textup{loc}(\R^d;\C)$ and $A\in L^2_\textup{loc}(\R^d;\R^d),$ the procedure developed in Subsection~\ref{preliminary_facts} in order to rigorously introduced the Hamiltonian $H_{A,V}$ formally defined in~\eqref{eq:ele-magnetic_Hamiltonian} must be adapted. The modification of the procedure consists merely in taking the Friedrichs extension of the operator initially defined on $C^{\infty}_0(\R^d\setminus \{0\})$
instead of $C^{\infty}_0(\R^d)$. To be more specific, 
we first introduce the closed quadratic form
\begin{equation}\label{eq:quadratic_AB}
	h_{A,V}^{(1)}[u]:=\int_{\R^d} \abs{\nabla_{\!A} u}^2\, dx + \int_{\R^d} (\Re V)_+ \abs{u}^2\, dx, \qquad u \in\mathcal{D}(h_{A,V}^{(1)}):= 
\overline{C^{\infty}_0(\R^d\setminus \{0\})}^{\normeq{\cdot}},
\end{equation}
where 
\begin{equation*}
	\normeq{u}^2:= 
h_{A,V}^{(1)}[u] + \norm{u}_{L^2(\R^d)}^2.
\end{equation*}
Assume that there exist $b,\beta \in [0,1)$ with
\begin{equation*}
	b^2 + \beta^2<1,
\end{equation*} 
such that, for any $u\in \mathcal{D}(h_{A,V}^{(1)}),$
\begin{equation}\label{eq:subordination}
	\int_{\R^d} (\Re V)_-\abs{u}^2\, dx\leq b^2 \int_{\R^d} \abs{\nabla_{\!A} u}^2\, dx,
	\qquad
	\int_{\R^d} \abs{\Im V} \abs{u}^2\, dx \leq \beta^2 \int_{\R^d} \abs{\nabla_{\!A} u}^2\, dx.
\end{equation}
Then, defining 
\begin{equation*}
	h_{A,V}^{(2)}[u]:= -\int_{\R^d} (\Re V)_- \abs{u}^2\, dx + i \int_{\R^d} \Im V \abs{u}^2\, dx,
\qquad
u \in \mathcal{D}(h_{A,V}^{(1)}),
\end{equation*}
the form $h_{A,V}^{(2)}$ is relatively bounded 
with respect to $h_{A,V}^{(1)}$, with the relative bound less than one.
Consequently, the sum $h_{A,V}:=h_{A,V}^{(1)} + h_{A,V}^{(2)}$ 
with domain $\mathcal{D}(h_{A,V}):=\mathcal{D}(h_{A,V}^{(1)})$
is a closed and sectorial form and $H_{A,V}$ is understood as 
the m-sectorial operator associated with $h_{A,V}$
via the representation theorem. 
Again, we abbreviate
\begin{equation*}
	\mathcal{D}_{A,V}:=\mathcal{D}(h_{A,V}).
\end{equation*}

\subsubsection{Proof of identity~\texorpdfstring{\eqref{eq:third}}{third}}
This subsection is concerned with the proof of Lemma~\ref{lemma:crucial_identities} in the present alternative framework. More specifically we will provide the proof of identity~\eqref{eq:third} only, which is the one whose changes are significant.
For the sake of clarity, 
we restate it with the alternative hypotheses assumed in this section.
(Without loss of generality, 
we consider just the situation in which $V^{(1)}=0$; 
indeed, the assumption \eqref{eq:hypotheses_V_alt} that we remove now
concerned the component $V^{(2)}$ only.)

\begin{lemma}\label{lemma:id-AB}
	Let $d\geq 2$. 
Let $A\in L^2_\textup{loc}(\R^d\setminus \{0\})$ be such that $\bm{B}\in L^2_\textup{loc}(\R^d\setminus \{0\})$ and let $V\in W^{1,1}_\textup{loc}(\R^d\setminus \{0\})$
be potentials satisfying~\eqref{eq:con_AB}. 
Let $u\in \mathcal{D}_{A,V}$ be any compactly supported solution of~\eqref{eq:weak_resolvent}, with $\lambda$ being any complex constant and $\abs{x}f\in L^2_\textup{loc}(\R^d),$ satisfying
	\begin{equation*}
		\left(
		\abs{x}^2\abs{\bm{B}}^2 + \abs{x}^2\abs{\Im V}^2 + [x\cdot \nabla \Re V]_+
		\right)\abs{u}^2 \in L^1_\textup{loc}(\R^d).
	\end{equation*}
	Then $[x\cdot \nabla \Re V]_-\abs{u}^2\in L^1_\textup{loc}(\R^d)$ and the following identity
	\begin{multline}
		\tag{\ref*{eq:third}'}
		\label{eq:third-alternative}
		2\int_{\R^d} \abs{\nabla_{\!A} u}^2\, dx 
		+ 2 \Im \int_{\R^d} x\cdot \bm{B} \cdot u \overline{\nabla_{\!A} u}\, dx 
		-\int_{\R^d} x\cdot \nabla \Re V^{(2)} \abs{u}^2\, dx
		-2 \Im \int_{\R^d} x\cdot \Im V^{(2)} u \overline{\nabla_{\!A} u}\, dx\\
		=-2\Im \lambda \Im \int_{\R^d} x \cdot u \overline{\nabla_{\!A} u}\, dx + d\Re \int_{\R^d} f \bar{u}\, dx + 2 \Re \int_{\R^d} f x\cdot \overline{\nabla_{\!A} u}\, dx
	\end{multline}
	holds true.
\end{lemma}  
\begin{proof}
For $d\geq 3$ we define $\xi\colon[0,\infty) \to [0,1]$ to be a smooth function such that
\begin{equation*}
	\xi(r):=
	\begin{system}
		 &0 &\text{if}\quad r\leq 1,\\
		&1 &\text{if}\quad r\geq 2,
	\end{system} 
\end{equation*}
and set $\xi_\varepsilon(x):=\xi(\abs{x}/\varepsilon).$
For $d=2,$ let $\xi \in C^\infty([0,1])$ such that $\xi=0$ in a right neighborhood of $0$ and $\xi=1$ in a left neighborhood of $1$; then we define the smooth function
\begin{equation*}
	\xi_\varepsilon(x):=
	\begin{system}
		&0 &\text{if}\quad \abs{x}\leq \varepsilon , \\
		&\xi(\log_2(\abs{x}/\varepsilon)) &\text{if}\quad \varepsilon\leq \abs{x}\leq 2\varepsilon, \\
		&1 &\text{if}\quad \abs{x}\geq 2\varepsilon.
	\end{system}
\end{equation*}
It comes from a straightforward computation to check that in both cases,
there exists a constant $\widetilde{c}>0$ such that the following control on the first derivatives  
\begin{equation}\label{der_xi_eps}
	\abs{\nabla \xi_\varepsilon}\leq \widetilde{c}/\varepsilon 
\end{equation} 
holds true.
%

%

We take as the test function in~\eqref{eq:weak_resolvent} a slight modification of the multiplier~\eqref{Morawetz_multiplier_reg} chosen above, namely 
\begin{equation}\label{eq:multiplier_AB_gen}
	v:=\Delta \phi u + \xi_\varepsilon \partial_k \phi [\partial_{k,A}^{\delta} u + \partial_{k,A}^{-\delta} u] \qquad \text{with} \qquad
	\phi(x):=\abs{x}^2,
\end{equation}
where 
\begin{equation*}
	\partial_{k,A}^{\delta}u:= \partial_k^\delta u + i A_k u, \qquad k=1,2,\dots, d, 
\end{equation*}	
	with $\partial_k^\delta$ defined as in~\eqref{eq:difference_quotient}.
More specifically,
\begin{equation}\label{eq:multiplier_AB}
	v = 
2du + 2\xi_\varepsilon x_k [\partial_{k,A}^{\delta} u + \partial_{k,A}^{-\delta} u].
\end{equation} 
 Observe that in this framework we do not need the truncation of the magnetic potential. 

Mimicking the arguments of Section~\ref{subsec:proof_identities},
one can show that $v$ defined as in~\eqref{eq:multiplier_AB} 
belongs to $\mathcal{D}_{A,V}.$
In fact, one has
$v\in L^2(\R^d),$ 
$\partial_{l,A}v:=(\partial_l+iA_l)v\in L^2(\R^d)$ for any $l=1, \dots, d$ 
and $\sqrt{(\Re V)_+}v\in L^2(\R^d)$.
We comment just on $\xi_\varepsilon x_k \partial_{k,A}^\delta u$ in~\eqref{eq:multiplier_AB}. Being $\xi_\varepsilon$ supported off the origin, 
$A_k\in L^\infty(\supp\xi_\varepsilon)$, 
therefore $\xi_\varepsilon x_k \partial_{k,A}^\delta u:=\xi_\varepsilon x_k (\partial_k^\delta + i A_k)u\in L^2(\R^d).$ Now we want to show that $\partial_{l,A}[\xi_\varepsilon x_k \partial_{k,A}^\delta u]\in L^2(\R^d).$
First observe that using the chain rule for magnetic derivatives~\eqref{eq:chain_rule_magnetic}, one can write
\begin{equation*}
	\partial_{l,A}[\xi_\varepsilon x_k \partial_{k,A}^\delta u]
	=v_1 + v_2,
\end{equation*}
where
\begin{equation*}
	v_1:=\xi_\varepsilon \partial_{l,A}[x_k \partial_{k,A}^\delta u],
	\quad \text{and} \quad v_2:=\partial_l \xi_\varepsilon [x_k \partial_{k,A}^\delta u].
\end{equation*}%
Clearly, exactly as above, $v_2\in L^2(\R^d).$  Using again that $\norm{A_k}_{L^\infty(\supp \xi_\varepsilon)}<\infty$ and the fact that $x_k \partial_{k,A}^\delta u= x_k \partial_{k,A}^{\delta, N} u$ with $N= \norm{A_k}_{L^\infty(\supp \xi_\varepsilon)},$ where $\partial_{k,A}^{\delta, N}$ are defined as in~\eqref{eq:regularised_magnetic_derivative}, 
one can reason as in Section~\ref{subsec:proof_identities} to conclude that $v_1\in L^2(\R^d)$ as well
(observe that here it comes into play the assumption $\partial_l A_k\in L^\infty(\R^d\setminus \{0\}),$ as in the previous section it came into play the assumption $\partial_l A_k\in L^p_\textup{loc}(\R^d)$ with $p$ as in~\eqref{eq:hypotheses_A}).   
It remains just to prove that $\sqrt{(\Re V)_+} [\xi_\varepsilon x_k \partial_{k,A}^\delta u]\in L^2(\R^d)$, but this follows immediately observing that, on the support of $\xi_\varepsilon,$  $(\Re V)_+$ is bounded. 

Now we are in position to prove identity~\eqref{eq:third-alternative}. Also in this section we proceed in a greater generality by considering $\phi$ in~\eqref{eq:multiplier_AB_gen} to be an arbitrary smooth function $\phi\colon \R^d \to \R.$ After we will plug in our choice $\phi(x)=\abs{x}^2.$ We consider identity~\eqref{eq:weak_resolvent} with the test function $v$ as in~\eqref{eq:multiplier_AB} and we take the real part. Each resulting integrals are treated separately.

\subsubsection*{$\bullet$ Kinetic term}
	Let us start with the ``kinetic'' part of~\eqref{eq:weak_resolvent}, 
\emph{i.e.}~\eqref{eq:kinetic_part}.
	Using 
	\begin{equation*}
		\begin{split}
		\overline{\partial_{l,A} v}=& 
		(\partial_l \Delta \phi)\bar{u} 
		+ \Delta \phi \overline{\partial_{l,A}u}\\
		&	+ \xi_\varepsilon \partial_{lk} \phi [\overline{\partial_{k,A}^\delta u} + \overline{\partial_{k,A}^{-\delta} u}] 
		+ \xi_\varepsilon \partial_k \phi \,[ \overline{\partial_{l,A} \partial_{k,A}^{\delta} u} + \overline{\partial_{l,A} \partial_{k,A}^{-\delta} u}] 
		+ \partial_l \xi_\varepsilon \partial_k \phi [\overline{\partial_{k,A}^\delta u} + \overline{\partial_{k,A}^{-\delta}u}],
		\end{split}
	\end{equation*}
	we write $K=K_0^\varepsilon + K_1 + K_2 + K_3^\varepsilon + K_4^\varepsilon$ with $K_1$ and $K_2$ as in~\eqref{eq:K's} and 
	\begin{equation*}
		K_0^\varepsilon:= \Re \int_{\R^d} \partial_l\xi_\varepsilon \partial_k \phi \partial_{l,A} u\, [\overline{\partial_{k,A}^\delta u} + \overline{\partial_{k,A}^{-\delta}u}],
		\end{equation*}
		\begin{equation*}
		K_3^\varepsilon:=\Re \int_{\R^d} \xi_\varepsilon \partial_{lk}\phi\, \partial_{l,A} u\, [\overline{\partial_{k,A}^{\delta} u} + \overline{\partial_{k,A}^{-\delta} u}], 
\qquad
K_4^\varepsilon:=\Re \int_{\R^d} \xi_\varepsilon\partial_k \phi \, \partial_{l,A} u\, [\overline{\partial_{l,A} \partial_{k,A}^{\delta} u} + \overline{\partial_{l,A} \partial_{k,A}^{-\delta} u}].
	\end{equation*}
	As regards with $K_4^\varepsilon,$ proceeding in the same way 
as done in Section~\ref{subsec:proof_identities}
to treat the term $K_4$, we end up with
	\begin{equation*}
		K_4^\varepsilon=K_{4,1,1}^\varepsilon + K_{4,1,2}^\varepsilon + K_{4,2}^\varepsilon, 
	\end{equation*}
	where
	\begin{equation*}
		K_{4,1,1}^\varepsilon= -\frac{1}{2} \int_{\R^d} \{\partial_k^{-\delta}(\xi_\varepsilon \partial_k \phi) +\partial_k^{\delta}(\xi_\varepsilon \partial_k \phi)\} \abs{\nabla_{\!A} u}^2, 
		\qquad
		K_{4,1,2}^\varepsilon=-\frac{\delta}{2} \int_{\R^d} \{\xi_\varepsilon \partial_k \phi - \tau_{k}^\delta (\xi_\varepsilon \partial_k \phi)\} \abs{\partial_k^{\delta} \partial_{l,A} u}^2
	\end{equation*}
	and
	\begin{equation*}
		K_{4,2}^\varepsilon= \Im \int_{\R^d} \xi_\varepsilon \partial_k \phi \partial_{l,A}u \big[\partial_l A_k\bar{u} - (\partial_k^\delta A_l) \tau_k^\delta \bar{u} \big]
+ \Im \int_{\R^d} \xi_\varepsilon \partial_k \phi \partial_{l,A}u \big[\partial_l A_k\bar{u} - (\partial_k^{-\delta} A_l) \tau_k^{-\delta} \bar{u} \big].
	\end{equation*}
	Now we choose $\phi(x):= \abs{x}^2.$ Using~\eqref{eq:phi_derivatives} we get
	\begin{equation*}
		K_0^\varepsilon=2\Re \int_{\R^d} \partial_l \xi_\varepsilon x_k \partial_{l,A} u\,  [\overline{\partial_{k,A}^\delta u} + \overline{\partial_{k,A}^{-\delta} u} ]\, dx,
	\end{equation*}
	\begin{equation*}	
		K_1=0,\qquad
		K_2=2d \int_{\R^d} \abs{\nabla_{\!A} u}^2\, dx,\qquad
	  K_3^\varepsilon=2 \Re \int_{\R^d} \xi_\varepsilon \partial_{l,A} u \,[\overline{\partial_{l,A}^{\delta} u} + \overline{\partial_{l,A}^{-\delta} u}]\, dx,
	\end{equation*}
and
\begin{equation*}
	K_{4,1,1}^\varepsilon=-\int_{\R^d} \{ \partial_k^{-\delta}(\xi_\varepsilon x_k) + \partial_k^{\delta}(\xi_\varepsilon x_k)\} \abs{\nabla_{\!A} u}^2\, dx,
	\qquad
	K_{4,1,2}^\varepsilon= \int_{\R^d} \partial_k^\delta (\xi_\varepsilon x_k) \abs{\tau_k^\delta \partial_{l,A} u - \partial_{l,A} u}^2\, dx,
\end{equation*}
\begin{equation*}
	K_{4,2}^\varepsilon= 2 \Im \int_{\R^d} \xi_\varepsilon x_k \partial_{l,A}u \big[\partial_l A_k\bar{u} - (\partial_k^\delta A_l) \tau_k^\delta \bar{u} \big]\, dx 
+ 2\Im \int_{\R^d} \xi_\varepsilon x_k \partial_{l,A}u \big[\partial_l A_k\bar{u} - (\partial_k^{-\delta} A_l) \tau_k^{-\delta} \bar{u} \big]\, dx.
\end{equation*}
Now we need the following analogous version to Lemma~\ref{lemma:auxiliary}.
\begin{lemma}\label{lemma:auxiliary-bis}
	Under the hypotheses of Lemma~\ref{lemma:id-AB}, 
the limits 
	\begin{equation*}
		\partial_{l,A}^\delta u \xrightarrow{\delta \to 0} \partial_{l,A} u\qquad \text{in}\quad L^2_\textup{loc}(\R^d\setminus \{0\})
	\end{equation*}
	and
	\begin{equation*}
		(\partial_k^\delta A_l)\tau_k^\delta u \xrightarrow{\delta \to 0} \partial_k A_l u\qquad \text{in}\quad L^2_\textup{loc}(\R^d\setminus \{0\})
	\end{equation*}
hold true.
\end{lemma}

Using Lemma~\ref{lemma:auxiliary-bis} and letting $\delta$ go to zero, 
it is easy to see that
\begin{equation}\label{K_eps}
	\begin{split}
	&K_0^\varepsilon \xrightarrow{\delta \to 0} 4\Re \int_{\R^d} \partial_l \xi_\varepsilon x_k \partial_{l,A}u \overline{\partial_{k,A} u}\, dx,\\
	&K_3^\varepsilon \xrightarrow{\delta \to 0} 4\int_{\R^d} \xi_\varepsilon \abs{\partial_{l,A} u}^2\, dx,\\
	&\begin{aligned}
		K_{4,1,1}^\varepsilon \xrightarrow{\delta \to 0} 
		&-2\int_{\R^d} \partial_k(\xi_\varepsilon x_k)\abs{\nabla_{\!A} u}^2\, dx\\
		&=-2\int_{\R^d} \partial_k \xi_\varepsilon x_k \abs{\nabla_{\!A} u}^2\, dx -2d \int_{\R^d}\xi_\varepsilon \abs{\nabla_{\!A} u}^2\, dx,
	\end{aligned}
	\\
	&K_{4,1,2}^\varepsilon \xrightarrow{\delta \to 0} 0,\\
	&K_{4,2}^\varepsilon \xrightarrow{\delta \to 0} 4 \Im \int_{\R^d} \xi_\varepsilon x_k \partial_{l,A}u [\partial_l A_k - \partial_k A_l]\bar{u}\, dx.
	\end{split}
\end{equation}
Now we want to see what happens in the limit of~$\varepsilon$ approaching zero. 
In order to do that we will use the following lemma.
\begin{lemma}\label{limit}
	Let $g\in L^1(\R^d)$ and let $\xi_\varepsilon$ be defined as above. Then
	\begin{equation}\label{limit1_2}
	\int_{\R^d} \xi_\varepsilon g\, dx \xrightarrow{\varepsilon \to 0} \int_{\R^d} g\, dx
\qquad \text{and}\qquad
		\int_{\R^d} \partial_l \xi_\varepsilon x_k g\, dx \xrightarrow{\varepsilon \to 0} 0
		\quad k,l=1,2\dots,d.
	\end{equation}
	\begin{proof}
	The first limit in~\eqref{limit1_2} immediately follows from the definition of $\xi_\varepsilon$ via the dominated convergence theorem.
	On the other hand, using~\eqref{der_xi_eps}, one has 
	\begin{equation*}
		\int_{\R^d} \abs{\partial_l \xi_\varepsilon} \abs{x_k} \abs{g}\, dx\leq 2\widetilde{c} \int_{\varepsilon<\abs{x}<2\varepsilon} \abs{g}\, dx\xrightarrow{\varepsilon \to 0} 0,
	\end{equation*}
	which yields the second limit in~\eqref{limit1_2},
 again from the dominated convergence theorem.
	\end{proof}
\end{lemma}
	
	Using Lemma~\ref{limit} and passing to the limit in~\eqref{K_eps}, one easily gets
	\begin{equation*}
		K \xrightarrow{\substack{\delta \to 0\\ \varepsilon \to 0}} 4\int_{\R^d} \abs{\partial_{l,A}u}^2\, dx + 4 \Im \int_{\R^d} x_k \partial_{l,A}u [\partial_l A_k - \partial_k A_l]\bar{u}\, dx.
	\end{equation*}
		Notice that here we have used that, by hypothesis, $\abs{x}^2 \abs{\bm{B}}^2\abs{u}^2\in L^1_\textup{loc}(\R^d).$

\subsubsection*{$\bullet$ Source term}
Now consider simultaneously
the ``source'' and ``eigenvalue''
parts of~\eqref{eq:weak_resolvent}, \emph{i.e.}~\eqref{eq:source_part}. 
Plugging in~\eqref{eq:source_part} our chosen test function $v$ defined in~\eqref{eq:multiplier_AB_gen}, we can write $F=F_1+F_2^\varepsilon+F_3^\varepsilon+F_4^\varepsilon$ with $F_1$ as in~\eqref{F's} and 
\begin{align*}
	F_2^\varepsilon&:=\Re \lambda \Re \int_{\R^d} \xi_\varepsilon \partial_k \phi \, u [\overline{\partial_{k, A}^{\delta} u} + \overline{\partial_{k,A}^{-\delta} u}],
\\
	F_3^\varepsilon&:=-\Im \lambda \Im \int_{\R^d} \xi_\varepsilon \partial_k \phi \, u [\overline{\partial_{k, A}^{\delta} u} + \overline{\partial_{k,A}^{-\delta} u}] ,
&
	F_4^\varepsilon&:= \Re \int_{\R^d} f\{\Delta \phi \bar{u} +\xi_\varepsilon \partial_k \phi \, [\overline{\partial_{k, A}^{\delta} u} + \overline{\partial_{k,A}^{-\delta} u}]\}.
\end{align*}
As regards with $F_2^\varepsilon,$ 
proceeding as in Section~\ref{subsec:proof_identities} when we treated $F_2,$ we end up with
\begin{equation*}
	F_2^\varepsilon= F_{2,1}^\varepsilon + F_{2,2}^\varepsilon
\end{equation*}
with
\begin{equation*}
	F_{2,1}^\varepsilon=-\frac{1}{2}\Re\lambda \int_{\R^d} \{\partial_k^{-\delta} (\xi_\varepsilon \partial_k \phi) + \partial_k^{\delta} (\xi_\varepsilon \partial_k \phi)\}\abs{u}^2
	\quad \text{and}\quad
	F_{2,2}^\varepsilon:=-\frac{\delta}{2} \Re \lambda \int_{\R^d} \xi_\varepsilon \partial_k \phi\, \{ \abs{\partial_k^\delta u}^2 - \abs{\partial_k^{-\delta} u}^2 \}.
\end{equation*}
Choosing $\phi(x):=\abs{x}^2$ 
in the previous identities and using~\eqref{eq:phi_derivatives} give
\begin{align*}
	&F_1=2d \Re\lambda \int_{\R^d} \abs{u}^2\, dx,\\
	&F_{2,1}^\varepsilon=-\Re\lambda \int_{\R^d} \{\partial_k^{-\delta} (\xi_\varepsilon x_k) + \partial_k^{\delta} (\xi_\varepsilon x_k)\}\abs{u}^2\, dx,
	\qquad
	F_{2,2}^\varepsilon= -\delta \Re \lambda \int_{\R^d} \xi_\varepsilon x_k \, \{ \abs{\partial_k^\delta u}^2 - \abs{\partial_k^{-\delta} u}^2 \}\, dx , \\
 &F_3^\varepsilon=-2\Im \lambda \Im \int_{\R^d} \xi_\varepsilon x_k u\, [\overline{\partial_{k,A}^{\delta} u} + \overline{\partial_{k,A}^{-\delta} u}]\, dx,
	\qquad
	F_4^\varepsilon= \Re \int_{\R^d} f \{2d \bar{u} +2 \xi_\varepsilon x_k [\overline{\partial_{k,A}^{\delta} u} + \overline{\partial_{k,A}^{-\delta} u}]\}\, dx.
\end{align*}
Reasoning as above, one gets
\begin{equation*}
	\begin{split}
		&F_{2,1}^\varepsilon\xrightarrow{\delta \to 0} -2\Re\lambda \int_{\R^d} \partial_k \xi_\varepsilon x_k \abs{u}^2\, dx 
	-2d\Re \lambda \int_{\R^d} \xi_\varepsilon \abs{u}^2\, dx,\\
	&F_{2,2}^\varepsilon\xrightarrow{\delta \to 0} 0,\\
	&F_3^\varepsilon\xrightarrow{\delta \to 0}-4\Im \lambda \Im \int_{\R^d} \xi_\varepsilon x_k u \overline{\partial_{k,A}u}\, dx,\\
	&F_4^\varepsilon\xrightarrow{\delta \to 0}\Re \int_{\R^d} f\{2d \bar{u} + 4\xi_\varepsilon x_k \overline{\partial_{k,A} u}\}\,dx. 
	\end{split}
\end{equation*}
Using Lemma~\ref{limit}, we conclude that
\begin{equation*}
	F\xrightarrow{\substack{\delta\to 0\\ \varepsilon \to 0}} 
	-4\Im \lambda \Im \int_{\R^d} x_k u\overline{\partial_{k,A}u}\, dx
  +\Re \int_{\R^d} f\{2d \bar{u} + 4 x_k \overline{\partial_{k,A} u}\}\,dx.
\end{equation*}

\subsubsection*{$\bullet$ Electric potential term}
Let us now consider the contribution 
of the ``potential'' part of~\eqref{eq:weak_resolvent}, 
\emph{i.e.}~\eqref{eq:potential_part}. 
Plugging $v$ defined as in~\eqref{eq:multiplier_AB_gen} 
into~\eqref{eq:potential_part}, we write $J=J_1 + J_2^\varepsilon$ with 
\begin{equation*}
	J_1:=\int_{\R^d} \Re V \Delta \phi \abs{u}^2
	\quad \text{and} \quad
	J_2^\varepsilon:= \Re \int_{\R^d} V \xi_\varepsilon \partial_k \phi u [\overline{\partial_{k,A}^{\delta}u} + \overline{\partial_{k,A}^{-\delta}u}].
\end{equation*}
Choosing $\phi(x):=\abs{x}^2$ in the previous identities and using~\eqref{eq:phi_derivatives}, we obtain
\begin{equation*}
	J_1= 2d \int_{\R^d} \Re V \abs{u}^2\, dx
	\quad \text{and} \quad
	J_2^\varepsilon= 2 \Re \int_{\R^d} \xi_\varepsilon x_k V  u [\overline{\partial_{k,A}^{\delta}u} + \overline{\partial_{k,A}^{-\delta}u}]\, dx. 
\end{equation*}
Now we write
\begin{equation*}
	J_2^\varepsilon=
	J_{2,1}^\varepsilon + J_{2,2}^\varepsilon,
\end{equation*}
where
\begin{equation*}
	J_{2,1}^\varepsilon:= 
	2\Re \int_{\R^d} \xi_\varepsilon x_k \Re V u [\overline{\partial_{k,A}^{\delta}u} + \overline{\partial_{k,A}^{-\delta}u}]\, dx
	\quad \text{and} \quad
	J_{2,2}^\varepsilon:=
	-2\Im \int_{\R^d} \xi_\varepsilon x_k \Im V u [\overline{\partial_{k,A}^{\delta}u} + \overline{\partial_{k,A}^{-\delta}u}]\, dx.
\end{equation*}
Using that $\Re V$ is bounded on $\supp \xi_\varepsilon,$ taking the limit as $\delta$ goes to zero, it follows from Lemma~\ref{lemma:auxiliary-bis} 
\begin{equation*}
	\begin{split}
	J_{2,1}^\varepsilon \xrightarrow{\delta \to 0} 
	& \ 4\Re \int_{\R^d} \xi_\varepsilon x_k \Re V  u \overline{\partial_{k,A}u}\, dx\\
	  & \ =-2\int_{\R^d} \partial_k\xi_\varepsilon x_k \Re V \abs{u}^2\, dx
	  -2d\int_{\R^d} \xi_\varepsilon \Re V \abs{u}^2\, dx
		-2\int_{\R^d} \xi_\varepsilon x_k \partial_k \Re V  \abs{u}^2\, dx,
		\end{split}
\end{equation*}
where in the last identity we have just integrated by parts.
Moreover, using that by hypothesis $\abs{x}^2\abs{\Im V}^2\abs{u}^2\in L^1_\textup{loc}(\R^d)$, we have
\begin{equation*}
	J_{2,2}^\varepsilon \xrightarrow{\delta \to 0} 
	-4\Im \int_{\R^d} \xi_\varepsilon x_k \Im V u \overline{\partial_{k,A}u}\, dx.
\end{equation*}
	
Finally, 
using that $\Re V \abs{u}^2$ and $[x_k\partial_k \Re V]_+ \abs{u}^2\in L^1(\R^d)$ and again $\abs{x}^2\abs{\Im V}^2\abs{u}^2\in L^1_\textup{loc}(\R^d),$ then Lemma~\ref{limit} gives
\begin{equation*}
		J\xrightarrow{\substack{\delta \to 0\\ \varepsilon \to 0}} 
		 -2\int_{\R^d} [x_k \partial_k \Re V]_+ \abs{u}^2\, dx
		 +2\int_{\R^d} [x_k \partial_k \Re V]_- \abs{u}^2\, dx
		 -4\Im \int_{\R^d} x_k \Im V u \overline{\partial_{k,A}u}\, dx.
\end{equation*}  
Observe that in order to pass to the limit in the integral involving $[x_k \partial_k \Re V]_-,$ we have used the monotone convergence theorem being $\xi_\varepsilon\nearrow 1$ as $\varepsilon$ tends to zero.

\medskip

In summary,
passing to the limit $\delta \to 0$ and $\varepsilon \to 0$ in~\eqref{eq:weak_resolvent} and multiplying the resulting identity by $1/2,$ one obtains~\eqref{eq:third-alternative}.
This concludes the proof of Lemma~\ref{lemma:id-AB}.
\end{proof}

\section{Absence of eigenvalues of matrix Schr\"odinger operators}
\label{Sec:Schro}
We start our investigation on Schr\"odinger operators by considering first the most delicate case represented by the non self-adjoint results Theorem~\ref{thm:Schro} (and its particular case Theorem~\ref{thm:Schro_simplified}) and the alternatives in $d=2$ given by Theorem~\ref{thm:alternative2d} and Theorem~\ref{thm:A-B}. The self-adjoint situation is treated afterward (Subsection~\ref{Ssec:self-adj}).
 
\subsection{Non self-adjoint case}

\begin{proof}[Proof of Theorem~\ref{thm:Schro}]
Let $u$ be any weak solution to the eigenvalue equation
\begin{equation}\label{eq:evs_Schro}
	H_\textup{S}(A, \bm{V}) u= \lambda u
\end{equation}
with $H_\textup{S}(A, \bm{V})$ being defined as in~\eqref{eq:Schro_int-unit} 
and $\lambda$ being any complex constant. More precisely, $u$ satisfies
\begin{equation}\label{eq:weak_resolvent_Schro}
	\int_{\R^d} \nabla_{\!A} u_j \cdot \overline{\nabla_{\!A} v_j}\, dx  + \int_{\R^d} V^{(2)} u_j \overline{v_j}\, dx= \lambda \int_{\R^d} u_j \overline{v_j}\, dx + \int_{\R^d} f_j \overline{v_j}\, dx 
\end{equation} 
for $j=1,2\dots,n$ and for any $v_j\in \mathcal{D}_{A,V}.$

Here, since we want to use directly the estimate in Lemma~\ref{lemma:byproduct},
we have defined $f:=-V^{(1)} u.$
In passing, observe that by virtue of our hypothesis~\eqref{hyp:V^1}, it is not difficult to check that $f,$ so defined, satisfies
\begin{equation}\label{Schro_condition:f}
	\sum_{j=1}^n \norm{\abs{f_j}^{1/2}\abs{u_j}^{1/2}}_{L^2(\R^d)}^2\leq a_1^2 \norm{\nabla_{\!A} u^-}_{[L^2(\R^d)]^n}^2
\qquad \text{and}\qquad
	\norm{\abs{x}f}_{[L^2(\R^d)]^n}\leq a_2 \norm{\nabla_{\!A} u^-}_{[L^2(\R^d)]^n},
\end{equation}
 with $a_1$ and $a_2$ as in~\eqref{hyp:V^1} and $u^-$ as in~\eqref{eq:u^-}. Notice that here we have used that $\abs{u}=\abs{u^-}.$

  The strategy of our proof is to show that, under the hypotheses of Theorem~\ref{thm:Schro}, $u$ is identically zero. In order to do that, as customary, we split the proof into two cases: $\abs{\Im\lambda}\leq \Re\lambda$ and $\abs{\Im\lambda}>\Re \lambda.$

\subsubsection*{$\bullet$ Case $\abs{\Im \lambda}\leq \Re \lambda.$}

Since $u_j$, for $j=1,2,\dots, n$, is a solution to~\eqref{eq:weak_resolvent_Schro}, 
we can use directly Lemma~\ref{lemma:byproduct} to get the estimate 
\begin{multline*}
		\norm{\nabla_{\!A} u_j^-}_{L^2(\R^d)}^2 
		+ (\Re \lambda)^{-1/2} \abs{\Im \lambda} 
		\Bigg[
		\int_{\R^d} \abs{x} \abs{\nabla_{\!A} u_j^-}^2\, dx
		- \frac{(d-1)}{2} \int_{\R^d} \frac{~\abs{u_j^-}^2}{\abs{x}}\, dx
		+ \int_{\R^d} \abs{x} (\Re V^{(2)})_+ \abs{u_j^-}^2\, dx 
		\Bigg]
		\\
		\leq 
		2\Big( 
		 \norm{\abs{x}\abs{\bm{B}} u_j^-}_{L^2(\R^d)} 
		 +\norm{\abs{x} \Im V^{(2)} u_j^-}_{L^2(\R^d)}
		 +\norm{\abs{x} f_j}_{L^2(\R^d)}
		\Big) \norm{\nabla_{\!A} u_j^-}_{L^2(\R^d)}
		\\
		+(d-1) \norm{\abs{f_j}^{1/2} \abs{u_j^-}^{1/2}}_{L^2(\R^d)}^2
		+ \norm{[\partial_r(\abs{x}\Re V^{(2)})]_+^{1/2} u_j^-}_{L^2(\R^d)}^2
		\\
		+ \Big(
		\norm{\abs{x} (\Re V^{(2)})_- u_j^-}_{L^2(\R^d)}
		+\norm{\abs{x} f_j}_{L^2(\R^d)}
		\Big)
		\Big(
		\norm{\abs{\Im V^{(2)}}^{1/2} u_j^-}_{L^2(\R^d)}
		+\norm{\abs{f_j}^{1/2} \abs{u_j^-}^{1/2}}_{L^2(\R^d)}
		\Big).
	\end{multline*}
Summing over $j=1,2, \dots, n$ and using the Cauchy--Schwarz inequality for discrete measures, we easily obtain
\begin{multline*}
		\norm{\nabla_{\!A} u^-}_{[L^2(\R^d)]^n}^2 
		+ (\Re \lambda)^{-1/2} \abs{\Im \lambda} 
		\Bigg[
		\int_{\R^d} \abs{x} \abs{\nabla_{\!A} u^-}^2\, dx
		- \frac{(d-1)}{2} \int_{\R^d} \frac{~\abs{u^-}^2}{\abs{x}}\, dx
		+ \int_{\R^d} \abs{x} (\Re V^{(2)})_+ \abs{u^-}^2\, dx 
		\Bigg]
		\\
		\leq 
		2\Big( 
		 \norm{\abs{x}\abs{\bm{B}} u^-}_{[L^2(\R^d)]^n} 
		 +\norm{\abs{x} \Im V^{(2)} u^-}_{[L^2(\R^d)]^n}
		 +\norm{\abs{x} f}_{[L^2(\R^d)]^n}
		\Big) \norm{\nabla_{\!A} u^-}_{[L^2(\R^d)]^n}
		\\
		+(d-1) \sum_{j=1}^n \norm{\abs{f_j}^{1/2} \abs{u_j^-}^{1/2}}_{L^2(\R^d)}^2
		+ \norm{[\partial_r(\abs{x}\Re V^{(2)})]_+^{1/2} u^-}_{[L^2(\R^d)]^n}^2
		\\
		+ \Big(
		\norm{\abs{x} (\Re V^{(2)})_- u^-}_{[L^2(\R^d)]^n}
		+\norm{\abs{x} f}_{[L^2(\R^d)]^n}
		\Big)
		\Big(
		\norm{\abs{\Im V^{(2)}}^{1/2} u^-}_{[L^2(\R^d)]^n}
		+ \Big(\sum_{j=1}^n\norm{\abs{f_j}^{1/2} \abs{u_j^-}^{1/2}}_{L^2(\R^d)}^2 \Big)^{1/2}
		\Big).
		\end{multline*}
Using assumptions \eqref{hyp:ReV^2-}--\eqref{hyp:B} 
together with~\eqref{Schro_condition:f}, one has
\begin{multline}\label{eq:ineq-I}
	\Big(1- \big(2c + 2\beta_2 + 2a_2 + (d-1)a_1^2 + \mathfrak{b}^2 + (b_2 + a_2)(\beta_1 + a_1)\big) \Big)\norm{\nabla_{\!A} u^-}_{[L^2(\R^d)]^n}^2\\ 
		+ (\Re \lambda)^{-1/2} \abs{\Im \lambda} 
		\Bigg[
		\int_{\R^d} \abs{x} \abs{\nabla_{\!A} u^-}^2\, dx
		- \frac{(d-1)}{2} \int_{\R^d} \frac{~\abs{u^-}^2}{\abs{x}}\, dx
		+ \int_{\R^d} \abs{x} (\Re V^{(2)})_+ \abs{u^-}^2\, dx 
		\Bigg]
		\leq 0.
\end{multline}
Now we need to estimate the squared bracket of the latter inequality, namely
\begin{equation}\label{eq:def-I}
	I:=\int_{\R^d} \abs{x} \abs{\nabla_{\!A} u^-}^2\, dx
		- \frac{(d-1)}{2} \int_{\R^d} \frac{~\abs{u^-}^2}{\abs{x}}\, dx
		+ \int_{\R^d} \abs{x} (\Re V^{(2)})_+ \abs{u^-}^2\, dx.
\end{equation}
Notice that, since $I$ appears as a ``coefficient'' of the \emph{positive} spectral quantity $(\Re\lambda)^{-1/2}\abs{\Im \lambda},$ we would like to get a positive contribution out of it to eventually discard this term in the previous estimate. Notice that only the second term in $I$  could spoil such positivity and therefore our aim is to control its magnitude in size by means of the positivity of the other terms in $I.$

To do so, we will proceed distinguishing the cases $d=1, d=2$ and $d\geq 3.$

Let us start with the easiest $d=1.$ In this situation the second term in $I$ cancels out and therefore $I\geq 0.$

We go further considering the case $d\geq 3.$ Here we employ the weighted magnetic Hardy-inequality
\begin{equation}\label{eq:weighted_Hardy_higher}
	\int_{\R^d} \abs{x} \abs{\nabla_{\!A} u}^2\, dx \geq \frac{(d-1)^2}{4} \int_{\R^d} \frac{~\abs{u}^2}{\abs{x}}\, dx.
\end{equation}
More specifically, using~\eqref{eq:weighted_Hardy_higher} we have
\begin{equation}\label{est:higher}
	I\geq \frac{d-3}{d-1}\int_{\R^d} \abs{x} \abs{\nabla_{\!A} u^-}^2\, dx + \int_{\R^d} \abs{x} (\Re V^{(2)})_+ \abs{u^-}^2\, dx,
\end{equation}
 which again is positive because we are considering $d\geq 3.$

Observe that in both cases treated so far, namely $d=1$ and $d\geq 3,$ the positivity of the real part of $V^{(2)},$ namely the term $\int_{\R^d} \abs{x}[\Re V^{(2)}]_+\abs{u}^2\,dx,$ did not really enter the proof of the positivity of $I.$  The situation is different when considering $d=2.$ Indeed, although~\eqref{eq:weighted_Hardy_higher} is valid also for $d=2,$ in this case the right-hand side of estimate~\eqref{est:higher} is not necessarily positive. Thus assumption~\eqref{positivity-d=2} comes into play here. Indeed, thanks to~\eqref{positivity-d=2}, it is immediate that
\begin{equation*}
	I:= \int_{\R^2}\abs{x}\abs{\nabla_{\!A} u^-}^2\,dx -\frac{1}{2}\int_{\R^2}\frac{~\abs{u^-}^2}{\abs{x}}\, dx + \int_{\R^2} \abs{x} (\Re V^{(2)})_+ \abs{u^-}^2\,dx \geq 0.
\end{equation*}
Hence we have proved that in any dimension $d\geq 1$ we have $I\geq 0.$ This yields that 
\begin{equation*}
	\Big(1- \big(2c + 2\beta_2 + 2a_2 + (d-1)a_1^2 + \mathfrak{b}^2 + (b_2 + a_2)(\beta_1 + a_1)\big) \Big)\norm{\nabla_{\!A} u^-}_{[L^2(\R^d)]^n}^2\leq 0,
\end{equation*}
which, by virtue of~\eqref{Schro_cond_numbers}, 
implies that $u^-$ (and therefore $u$) 
is identically equal to zero.

\subsubsection*{$\bullet$ Case $\abs{\Im \lambda}> \Re \lambda.$}

Let $u_j$ for $j=1,2,\dots, n$ be a solution to~\eqref{eq:weak_resolvent_Schro}. Choosing as a test function $v_j:= u_j$ and taking the real part of the resulting identity and adding/subtracting, instead of the real part, the imaginary part of the resulting identity, one gets
\begin{multline*}
	\int_{\R^d} \abs{\nabla_{\!A} u_j}^2\,dx + \int_{\R^d} (\Re V^{(2)})_+ \abs{u_j}^2\, dx - \int_{\R^d} (\Re V^{(2)})_- \abs{u_j}^2\, dx \pm \int_{\R^d} \Im V^{(2)}\abs{u_j}^2\, dx
	\\=(\Re \lambda \pm \Im \lambda) \int_{\R^d} \abs{u_j}^2\, dx
	 + \Re \int_{\R^d} f_j \overline{u_j}\, dx \pm \Im \int_{\R^d} f_j \overline{u_j}\, dx. 
\end{multline*} 
Summing over $j=1,2,\dots, n$ and discarding the positive term on the left-hand side involving $(\Re V^{(2)})_+$, one easily gets
\begin{multline*}
	\norm{\nabla_{\!A} u}_{[L^2(\R^d)]^n}^2 \\ \leq (\Re \lambda \pm \Im \lambda) \int_{\R^d} \abs{u}^2\, dx
	 + \int_{\R^d}(\Re V^{(2)})_- \abs{u}^2\, dx 
	 + \int_{\R^d}\abs{\Im V^{(2)}} \abs{u}^2\, dx
	 + 2 \sum_{j=1}^n \norm{\abs{f_j}^{1/2} \abs{u_j}^{1/2}}_{L^2(\R^d)}^2.
\end{multline*} 
Using the first inequalities 
in~\eqref{hyp:ReV^2-},~\eqref{hyp:ImV^2} and~\eqref{Schro_condition:f}, 
we have
\begin{equation*}
	\Big(1-(b_1^2 + \beta_1^2 + 2 a_1^2) \Big)\norm{\nabla_{\!A} u}_{[L^2(\R^d)]^n}^2\leq (\Re \lambda \pm \Im \lambda) \norm{u}_{[L^2(\R^d)]^n}^2. 
\end{equation*}
Therefore, since by the first inequality in~\eqref{Schro_cond_numbers} 
we have $b_1^2 + \beta_1^2 + 2 a_1^2<1,$ then $\Re \lambda \pm \Im\lambda \geq 0$ unless $u=0.$ But since $\abs{\Im \lambda}>\Re \lambda$ we conclude that $u=0.$

This concludes the proof of Theorem~\ref{thm:Schro}.
\end{proof}


Now we prove the alternative Theorem~\ref{thm:alternative2d} valid in $d=2.$ 

\begin{proof}[Proof of Theorem~\ref{thm:alternative2d}]
Since the proof follows analogously to the one of Theorem~\ref{thm:Schro} presented above, except for the analysis in the sector $\abs{\Im \lambda}\leq\Re \lambda,$ we shall comment just on this situation.

	As in the proof of Theorem~\ref{thm:Schro}, 
we want to estimate the term $I$ defined in~\eqref{eq:def-I}, 
which appears multiplied by the spectral coefficient $(\Re\lambda)^{-1/2}\abs{\Im \lambda}$ in~\eqref{eq:ineq-I}. A first application of the weighted inequality~\eqref{eq:weighted_Hardy_higher} gives
	\begin{equation}\label{eq:bound-I-2d}	
		I\geq -\frac{1}{4} \int_{\R^2} \frac{~\abs{u^-}^2}{\abs{x}}\, dx + \int_{\R^2} \abs{x} (\Re V^{(2)})_+ \abs{u^-}\, dx\\
		\geq -\frac{1}{4} \int_{\R^2} \frac{~\abs{u^-}^2}{\abs{x}}\, dx,
	\end{equation}
	where the last inequality follows by discarding the positive term involving the potential $V^{(2)}.$ Now, we proceed estimating the term $\int_{\R^2} \frac{\abs{u^-}^2}{\abs{x}}\, dx.$ In order to do that we will strongly use the following Hardy--Poincaré-type inequality
\begin{equation}\label{eq:Hardy-Poincare}	
	 \int_{B_R} \abs{\nabla \psi}^2\, dx \geq \frac{1}{4R}\int_{B_R}\frac{\abs{\psi}^2}{\abs{x}}\, dx,	
\end{equation}
	valid for all $\psi\in W^{1,2}_0(B_R),$ where $B_R:=\{x\in \R^2\colon \abs{x}<R\}$ denotes the open disk of radius $R>0$ (see~\cite{FKV2} for an explicit proof of~\eqref{eq:Hardy-Poincare}). 
	
	Following the strategy of~\cite{FKV2}, given two positive numbers $R_1<R_2,$ we introduce the function $\eta\colon [0,\infty)\to [0,1]$ such that $\eta=1$ on $[0, R_1],$ $\eta=0$ on $[R_2, \infty)$  and $\eta(r)=(R_2-r)/(R_2-R_1)$ for $r\in (R_1, R_2).$ We denote by the same symbol $\eta$ the radial function $\eta \circ r\colon \R^2 \to [0,1].$ Now, writing $u^-= \eta u^- +(1-\eta)u^-$ and using~\eqref{eq:Hardy-Poincare}, we have
	\begin{equation*}
		\begin{split}
			\int_{\R^2} \frac{\abs{u^-}^2}{\abs{x}}\,dx&\leq
			 2 \int_{B_{R_2}} \frac{(\eta \abs{u^-})^2}{\abs{x}}\, dx + 2\int_{\R^2} \frac{\big((1-\eta)\abs{u^-}\big)^2}{\abs{x}}\, dx\\
			 &\leq 8 R_2 \int_{B_{R_2}} \abs{\nabla(\eta \abs{u^-})}^2\, dx + \frac{2}{R_1}\int_{\R^2} \abs{u^-}^2\, dx\\
			 &\leq 16 R_2 \int_{\R^2} \abs{\nabla \abs{u^-}}^2\, dx + 16\frac{R_2}{(R_2-R_1)^2} \int_{\R^2} \abs{u^-}^2\, dx + \frac{2}{R_1}\int_{\R^2} \abs{u^-}^2\, dx.
		\end{split}	
	\end{equation*}
	Choosing $R_1=R_2/2$ and using the diamagnetic inequality~\eqref{eq:diamagnetic} give
	\begin{equation*}
		\int_{\R^2} \frac{\abs{u^-}^2}{\abs{x}}\, dx
		\leq 16 R_2 \int_{\R^2} \abs{\nabla_{\!A} u^-}^2\, dx + \frac{68}{R_2} \int_{\R^2} \abs{u^-}^2\, dx.
	\end{equation*}
	Now we fix conveniently $R_2$; 
namely, given any positive number $\epsilon$, 
we set $R_2:= \epsilon (\Re \lambda)^{1/2}/\abs{\Im \lambda}$ 
in the previous inequality. Then multiplying the resulting inequality by $(\Re \lambda)^{-1/2} \abs{\Im \lambda} \frac{1}{4}$, we get
	\begin{equation}\label{eq:and}
		\begin{split}
			 (\Re\lambda)^{-1/2}\abs{\Im \lambda}\frac{1}{4}\int_{\R^2} \frac{\abs{u^-}^2}{\abs{x}}\,dx&\leq  
			4 \epsilon \int_{\R^2} \abs{\nabla_{\!A} u^-}^2\, dx + \frac{17}{\epsilon} \abs{\Im \lambda} \int_{\R^2} \abs{u^-}^2\, dx
			\\
			&\leq 
			4 \epsilon \int_{\R^2} \abs{\nabla_{\!A} u^-}^2\, dx + \frac{17}{\epsilon}  \int_{\R^2} \abs{\Im \bm{V}} \abs{u^-}^2\, dx
			\\		
			&\leq 
			\left[
			4 \epsilon +
			\frac{17}{\epsilon} (a_1^2 + \beta_1^2) 
			\right]
			 \int_{\R^2} \abs{\nabla_{\!A} u^-}^2\, dx,
		\end{split}
	\end{equation}
	where in the first inequality we have used the restriction to the sector $\abs{\Im \lambda}\leq \Re\lambda,$ the second estimate follows from~\eqref{eq:second_const} with $f=0$ and the third inequality from~\eqref{hyp:V^1} and~\eqref{hyp:ImV^2}.

	Using that, from~\eqref{eq:bound-I-2d} and~\eqref{eq:and}, one has
	\begin{equation*}
		\begin{split}
		(\Re\lambda)^{-1/2} \abs{\Im \lambda}\,  I &\geq 
		-(\Re\lambda)^{-1/2} \abs{\Im \lambda}\, \frac{1}{4}
		\int_{\R^2} \frac{\abs{u^-}^2}{\abs{x}}\,dx\\
		&\geq-
		 \left[
			4 \epsilon +
			\frac{17}{\epsilon} (a_1^2 + \beta_1^2) 
		\right]
		 \int_{\R^2} \abs{\nabla_{\!A} u^-}^2\, dx	
		 \end{split}
	\end{equation*}
	and plugging this last bound in~\eqref{eq:ineq-I}, we get
	\begin{equation*}
		\left[1- \Big(2c + 2\beta_2 + 2a_2 + a_1^2 + \mathfrak{b}^2 + (b_2 + a_2)(\beta_1 + a_1) + 4 \epsilon +
			\frac{17}{\epsilon} (a_1^2 + \beta_1^2)  \Big) 
			\right]
			\norm{\nabla_{\!A} u^-}_{[L^2(\R^2)]^n}^2
		\leq 0.
	\end{equation*}
	From hypothesis~\eqref{eq:alt-Schr-cond-2d}, 
we therefore conclude that $u=0$ as above.
\end{proof}

Finally, we prove the two dimensional result in which the magnetic potential is fixed to be the Aharonov--Bohm one.
 
\begin{proof}[Proof of Theorem~\ref{thm:A-B}]
 As in the proof of Theorem~\ref{thm:alternative2d}, 
we need to estimate the term $I$ defined in~\eqref{eq:def-I}, 
which appears in~\eqref{eq:ineq-I}. 
Notice that in this specific case
(due to the triviality of the magnetic field,
everywhere except at the origin, see~\eqref{eq:B_A-B}), 
in~\eqref{eq:ineq-I} there does not appear the constant~$c$ 
related to the smallness condition assumed for~$B$. 
In order to estimate~$I$,
we will use the following weighted Hardy inequality, which is also an improvement upon~\eqref{eq:weighted_Hardy-2d}	, it reads
	\begin{equation}\label{eq:weighted-A-B}
		\int_{\R^2} \abs{x}\abs{\nabla_A \psi}^2\, dx \geq \left( \frac{1}{4} + \gamma^2 \right) \int_{\R^2} \frac{\abs{\psi}^2}{\abs{x}}\, dx, \qquad \forall\, \psi \in C^\infty_0(\R^d\setminus \{0\}),
	\end{equation}
	where $\gamma:=\dist\{\bar{\alpha}, \mathbb{Z}\}$ and $\bar{\alpha}$ is as in~\eqref{eq:mag-flux}  
(see~\cite[Lem.~3]{FKV2} for a proof of~\eqref{eq:weighted-A-B}).
 
	A first application of~\eqref{eq:weighted-A-B} gives
	\begin{equation}\label{eq:A-B-I}
		(\Re\lambda)^{-1/2} \abs{\Im \lambda}\, I\geq 
		- (\Re\lambda)^{-1/2} \abs{\Im \lambda} \,\left(\frac{1}{4} - \gamma^2 \right)\int_{\R^2} \frac{\abs{u^-}^2}{\abs{x}} \, dx,
	\end{equation}
	where we discarded the positive term in $I$ involving the potential $V^{(2)}.$ Notice that since we are assuming $\bar{\alpha}\notin \mathbb{Z},$ then $ \gamma\in (0, 1/2],$ this gives $1/4 - \gamma^2 \geq 0.$

	Now, we proceed estimating the term $\int_{\R^2} \frac{\abs{u^-}^2}{\abs{x}}\, dx.$ Given any positive number $R,$ we write
	\begin{equation*}
			\int_{\R^2} \frac{\abs{u^-}^2}{\abs{x}}\, dx =
			\int_{B_R} \frac{\abs{u^-}^2}{\abs{x}}\, dx + \int_{\R^2\setminus \overline{B_R}} \frac{\abs{u^-}^2}{\abs{x}}\, dx
			\leq R \int_{B_R} \frac{\abs{u^-}^2}{\abs{x}^2}\, dx
			+ \frac{1}{R} \int_{\R^2} \abs{u^-}^2\, dx,
	\end{equation*}
	where, also here, $B_R$ denotes the open disk of radius $R>0.$	
	
	Choosing in the previous inequality $R:=\epsilon \gamma^2 (\Re \lambda)^{1/2}/ \abs{\Im \lambda}$ with any positive constant $\epsilon,$ and multiplying the resulting estimate by the quantity $(\Re \lambda)^{-1/2} \abs{\Im \lambda} \left( \frac{1}{4} - \gamma^2\right),$ we get
	\begin{equation*}
	\begin{split} 
			(\Re \lambda)^{-1/2} \abs{\Im \lambda} \left( \frac{1}{4} - \gamma^2\right) \int_{\R^2}\frac{\abs{u^-}^2}{\abs{x}}\, dx	
			&\leq \left( \frac{1}{4} - \gamma^2\right)
			\left[
			 \epsilon \gamma^2  \int_{\R^2} \frac{\abs{u^-}^2}{\abs{x}^2}\, dx
			+ \frac{1}{\epsilon \gamma^2} \int_{\R^2} \abs{\Im \bm{V}} \abs{u^-}^2\, dx
			\right]			
			\\
			&\leq \left( \frac{1}{4} - \gamma^2\right) \left[ \epsilon + \frac{(a + \beta)}{\epsilon \gamma^3} \right] \int_{\R^2} \abs{\nabla_{\!A} u^-}^2\, dx.
	\end{split}	
	\end{equation*}	
In the first inequality we have used the restriction to the sector $\abs{\Im \lambda}\leq \Re\lambda,$ while in the second inequality we have used first the Hardy inequality~\eqref{eq:magnetic_Hardy} and then the hypotheses on the potential~\eqref{eq:V^1-AB} 
together with the second inequality of~\eqref{eq:V^2-AB}.
	Plugging the last estimate in~\eqref{eq:A-B-I} and the resulting estimate in~\eqref{eq:ineq-I}, and using an analog reasoning as in Remark~\ref{rmk:overabundant}, give
	\begin{equation*}
	\left[
	 1- \Bigg(2\beta + 2a + \frac{a}{\gamma} + \mathfrak{b}^2 + \frac{1}{\sqrt{\gamma}}(b + a)(\sqrt{a} + \sqrt{\beta} )+ \left( \frac{1}{4} - \gamma^2 \right) \left[\epsilon + \frac{(a+\beta)}{\epsilon \gamma^3}\right]\Bigg)
	\right]
	\norm{\nabla_{\!A} u^-}_{[L^2(\R^2)]^n}^2
		\leq 0.
	\end{equation*}
	From hypothesis~\eqref{eq:alt-Schr-cond-AB} we therefore conclude that $u=0$ as above.
\end{proof}

\subsection{Self-adjoint case: Proof of Theorem~\ref{thm:self-adjoint}}
\label{Ssec:self-adj}
Now we prove the much simpler and less involved analogous result to Theorem~\ref{thm:Schro} for self-adjoint Schr\"odinger operators, namely Theorem~\ref{thm:self-adjoint}.

\begin{proof}[Proof of Theorem~\ref{thm:self-adjoint}]
Let $u$ be any weak solution to the eigenvalues equation~\eqref{eq:evs_Schro}, with $\bm{V}$ real-valued. 

The proof of this theorem is based exclusively on the identity~\eqref{eq:third}. More precisely, using that $\bm{V}$ is real-valued, so necessarily $\Im \lambda=0,$ from~\eqref{eq:third} (with $f=0$) we get
\begin{equation}\label{id:self-adj}
	\begin{split}
	2\int_{\R^d} \abs{\nabla_{\!A} u_j}^2\, dx
	=&
	- 2\Im \int_{\R^d} x\cdot \bm{B}\cdot u_j \overline{\nabla_{\!A} u_j}\, dx
	+\int_{\R^d} \abs{x} \partial_r V^{(2)} \abs{u_j}^2\, dx\\
	&
	- d \int_{\R^d} V^{(1)}\abs{u}^2\, dx 
	- 2\Re \int_{\R^d} x \cdot V^{(1)} u_j \overline{\nabla_{\!A} u_j}\, dx.
	\end{split}
\end{equation} 
Observing that 
\begin{equation*}
	\int_{\R^d} \abs{x} \partial_r V^{(2)} \abs{u_j}^2\, dx
	\leq
	\int_{\R^d} [\abs{x} \partial_r V^{(2)}]_+ \abs{u_j}^2\, dx,
\end{equation*}
using the Cauchy--Schwarz inequality and summing over $j=1,2,\dots, n,$ one has
\begin{equation*}
	\begin{split}
	2\norm{\nabla_{\!A} u}_{[L^2(\R^d)]^n}^2
	\leq&  
	2\Big(\int_{\R^d} \abs{x}^2\abs{\bm{B}}^2 \abs{u}^2\, dx\Big)^{1/2}
	\norm{\nabla_{\!A} u}_{[L^2(\R^d)]^n}
	+\int_{\R^d} [\abs{x} \partial_r V^{(2)}]_+ \abs{u}^2\, dx\\
	&
	+ d \int_{\R^d} \abs{V^{(1)}}\abs{u}^2\, dx
	+2\Big(\int_{\R^d} \abs{x}^2\abs{V^{(1)}}^2 \abs{u}^2\, dx\Big)^{1/2} 
	\norm{\nabla_{\!A} u}_{[L^2(\R^d)]^n}. 
	\end{split}
\end{equation*}
Now, using~\eqref{hyp:V^1},~\eqref{hyp:B} and~\eqref{hyp:r-radial_der-V^2}, 
one easily gets
\begin{equation*}
	\Big( 
	2-(2c+ \mathfrak{b}^2 + d a_1^2 + 2a_2)
	\Big)\norm{\nabla_{\!A} u}_{[L^2(\R^d)]^n}^2\leq 0.
\end{equation*}
This immediately gives a contradiction in virtue of~\eqref{Schro_cond_numbers:s-a}. This concludes the proof. 
\end{proof}

In passing, observe that here we did not need to split the proof and proving separately absence of positive and non-positive eigenvalues. 
Indeed, we got the absence of the total point spectrum in just one step.

\begin{remark}[Two-dimensional Pauli operators as a special case]\label{rmk:A-C}
One reason for investigating matrix self-adjoint Schr\"odinger operators in this work, comes from our interest in pointing out a pathological behavior of the two dimensional purely magnetic (and so self-adjoint) Pauli Hamiltonian. From the explicit expression~\eqref{eq:Pauli2d} of the two dimensional Pauli operators, it is evident the relation with the scalar Schr\"odinger operator
	\begin{equation*} 
		-\nabla_{\!A}^2 + V^{(1)}
		\qquad \mbox{with} \qquad
V^{(1)}:= \pm B_{12}. 
	\end{equation*}
In this specific situation identity~\eqref{id:self-adj}, which was the crucial identity to prove absence of point spectrum in the self-adjoint situation, reads (after multiplying by $1/2$)
\begin{equation*}
	\int_{\R^2} \abs{\nabla_{\!A} u}^2\, dx =- \Im \int_{\R^2} x\cdot B_{12} u \overline{\nabla_{\!A} u}^\perp\, dx -\int_{\R^2} B_{12} \abs{u}^2\, dx - \Re \int_{\R^2} x\cdot B_{12} u \overline{\nabla_{\!A} u}\, dx.
\end{equation*} 
We stress that differently to the proof presented above, here the presence of the second term on the right-hand side involving the magnetic field does not allow us to get a contradiction. Indeed, roughly speaking, all the positivity coming from the left-hand side and that is customarily used to get the contradiction under the smallness assumption on the magnetic field is exploited to control the second term on the right-hand side (due to inequality~\eqref{B12LB}), therefore, using~\eqref{hyp:B}, one is left with a term of the type
\begin{equation*}
	-2c \norm{\nabla_A u}_{L^2(\R^2)}^2\leq 0,
\end{equation*}
which leads to no contradiction, however small is chosen the constant $c.$
\end{remark}

\section{Absence of eigenvalues of Pauli and Dirac operators}\label{sec:Pauli-Dirac}
This section is devoted to the proof of emptiness 
of the point spectrum of Pauli and Dirac Hamiltonians.

\subsection{Warm-up in the 3d case}\label{ss:warmup}
Even though the three dimensional setting proposed in the introduction is clearly covered by the more general results Theorem~\ref{Thm.Pauli} and Theorem~\ref{thm:magnetic_Dirac}, we decided to dedicate to the 3d case a separate section. Indeed, due to the physical relevance of this framework, we want to make it easier to spot the conditions which guarantee the absence of the point spectrum in this case, avoiding the interested reader working his/her way through the statements of the theorems in the general setting. 

\subsubsection{Absence of eigenvalues of Pauli operators: proof of Theorem~\protect\ref{thm:main_result}}
Let $u$ be any weak solution to the eigenvalue equation
\begin{equation}\label{eq:evs_Pauli}
	H_\textup{P}(A,\bm{V})u=\lambda u,
\end{equation} 
with $H_\textup{P}(A, \bm{V})$ 
defined as in~\eqref{Pauli_perturbed} and where $\lambda$ is any complex constant. 

Using~\eqref{Pauli_perturbed} 
and the decomposition $\bm{V}=\bm{V^{(1)}} + \bm{V^{(2)}},$ problem~\eqref{eq:evs_Pauli} can be written as an eigenvalue problem for matrix Schr\"odinger operators, namely  
\begin{equation*}
	H_\textup{S}(A, \bm{W})u =\lambda u, 
\end{equation*}
where $H_\textup{S}(A,\bm{W})$ is defined 
in~\eqref{eq:Schro_int-unit} and where $\bm{W}=\bm{W}^{(1)} + \bm{W^{(2)}}$ with
\begin{equation}\label{components}
	\bm{W^{(1)}} := \sigma \cdot B + \bm{V^{(1)}} 
	\qquad \text{and} \qquad
	\bm{W^{(2)}} := \bm{V^{(2)}}.
\end{equation} 
In light of the assumptions in~\eqref{eq:conditions_B_V-3d} about $\bm{V^{(1)}}$ and $B,$ which intrinsically are \emph{both} full-subordination conditions to the magnetic Dirichlet form, it is indeed natural to treat $\bm{V^{(1)}}$ and $B$ in a unified way defining $\bm{W^{(1)}}$ as in~\eqref{components}.

Assuming the hypotheses of Theorem~\ref{thm:main_result}
and using that $\abs{\sigma}=\sqrt{3}$ 
due to the fact that the Pauli matrices have norm one,
one easily verifies the bound
\begin{equation*}
	\int_{\R^3} \abs{x}^2 \abs{\bm{W^{(1)}}}^2 \abs{u}^2\, dx\leq (a+ \sqrt{3}c)^2 \int_{\R^3} \abs{\nabla_{\!A} u}^2\, dx.
\end{equation*}
Hence, hypotheses of Theorem~\ref{thm:Schro_simplified} are satisfied (with $\bm{W}$ instead of $\bm{V}$ and with $a+ \sqrt{3}c$ as a replacement for $a$ in~\eqref{eq:conditions_B_V}). From this we conclude the absence of eigenvalues of $H_\textup{S}(A, \bm{W})$ and, in turn clearly of $H_\textup{P}(A, \bm{V}),$ which is the thesis.
\qed

\subsubsection{Absence of eigenvalues of Dirac operators: proof of Theorem~\protect\ref{thm:magnetic_Dirac3d}}

Now we are in position to prove Theorem~\ref{thm:magnetic_Dirac3d}. As we will see, it follows as a consequence of the corresponding result for Pauli operators, namely Theorem~\ref{thm:main_result}.

Let $u$ be any solution to the eigenvalues equation
\begin{equation}\label{eq:evs_Dirac}
	H_\textup{D}(A) u=k u,
\end{equation}
with $H_\textup{D}(A):=H_\textup{D}(A,\bm{0})$ the three dimensional self-adjoint Dirac operator defined in~\eqref{magnetic_Dirac} and where $k$ is any real constant.   
A second application of the Dirac operator to the eigenvalues problem~\eqref{eq:evs_Dirac} gives that if $u$ is a solution to~\eqref{eq:evs_Dirac}, then it satisfies 
\begin{equation*}
	H_\textup{D}(A)^2 u= k^2 u.
\end{equation*}
More explicitly, using expression~\eqref{favorable_form} and defining $u_{1,2}:=(u_1,u_2)$ and $u_{3,4}:=(u_3,u_4)$ the two-vectors with components respectively the first and the second component of $u=(u_1,u_2,u_3,u_4)$, and the third and the fourth, one gets that $u_{1,2}$ and $u_{3,4}$ satisfy
\begin{equation*}
	\begin{system}
		&H_\textup{P}(A) u_{1,2} + \quarter u_{1,2}= k^2 u_{1,2},\\
		&H_\textup{P}(A) u_{3,4} + \quarter u_{3,4}= k^2 u_{3,4}.
	\end{system}
\end{equation*}
In other words, the two-vectors $u_{1,2}$ and $u_{3,4}$ are solutions to the eigenvalue problems associated to the shifted Pauli operators $H_\textup{P}(A) + \quarter$ with eigenvalues $k^2.$  

Notice that since~\eqref{cond_Dirac_B} holds for any $u=(u_1,u_2,u_3,u_4),$  in particular it holds for the four-vector $(u_1,u_2,0,0)$ and $(0,0,u_3,u_4).$ This fact implies that the second condition in~\eqref{eq:conditions_B_V-3d} of Theorem~\ref{thm:main_result} holds with the same constant $c$ as in~\eqref{cond_Dirac_B}. This means that we are in the hypotheses of Theorem~\ref{thm:main_result} (once we set a purely magnetic framework, namely $\bm{V}=0$), so $H_\textup{P}(A)$ has no eigenvalues. As a consequence, the shifted operator $H_\textup{P}(A) + \quarter\bm{I_{\C^2}}$ has no eigenvalues too. Hence $u_{1,2}$ and $u_{3,4}$ are vanishing and with them $u=(u_{1,2}, u_{3,4})$ itself. 

This concludes the proof of Theorem~\ref{thm:magnetic_Dirac3d}.          
\qed

\subsection{Absence of eigenvalues of Pauli operators in any dimension}
\label{Sec:Pauli-Schro}

Now we are in position to prove the general Theorem~\ref{Thm.Pauli}.

\begin{proof}[Proof of Theorem~\ref{Thm.Pauli}]

We divide the proof depending on the parity of the space dimension.
 
\subsubsection{Odd dimensions}
In odd dimensions,  
the proof follows the same scheme as the one presented 
in the three-dimensional case.

Looking at expression~\eqref{Pauli-odd} 
and using the decomposition of $\bm{V}=\bm{V^{(1)}}+\bm{V^{(2)}},$ 
one defines $\bm{W}=\bm{W^{(1)}} + \bm{W^{(2)}}$ such that
	\begin{equation*}
		\bm{W^{(1)}}= -\frac{i}{2}a\cdot \bm{B}\cdot a + \bm{V^{(1)}}
		\qquad \text{and} \qquad
		\bm{W^{(2)}}= \bm{V^{(2)}}.
	\end{equation*}
	It is easy to see that 
	\begin{equation*}
		\int_{\R^d} \abs{x}^2\abs{\bm{W^{(1)}}}^2 \abs{u}^2\, dx\leq \Big(a + \frac{d}{2} c\Big)^2 \int_{\R^d} \abs{\nabla_{\!A}u}^2\, dx, 
	\end{equation*}
	where we have used the validity of~\eqref{eq:conditions_B_V} and the fact that $\abs{a}=\sqrt{d}$ (see Remark~\ref{rmk:norm=1}).
	
	Thus, the proof follows exactly as the one of Theorem~\ref{thm:main_result} using, this time, the general result for Schr\"odinger operators Theorem~\ref{thm:Schro}.	

\subsubsection{Even dimensions}
	Let $u$ be any solution to the eigenvalue problem
	\begin{equation*}
		H_\textup{P}^\textup{even}(A, \bm{V}) u=\lambda u,
	\end{equation*}
	where $H_\textup{P}^\textup{even}(A, \bm{V})$ is defined in~\eqref{Pauli-even} and $\lambda$ is any complex constant. 
	In passing notice that according to~\eqref{eq:n'(d)}, since $d$ is even, then $n'(d)=n(d).$ 
	
Defining $u_\textup{up}:=(u_1, u_2, \dots, u_{n(d)/2})$ 
and $u_\textup{down}:=(u_{n(d)/2 +1}, u_{n(d)/2 +2}, \dots, u_{n(d)})$, 
the $n(d)/2$-vectors with components respectively the first half and the second half of the components of $u=(u_1, u_2, \dots, u_{n(d)}),$ one gets
	\begin{equation*}
		\begin{system}
			H_\textup{S}(A, \bm{W_\textup{up}})u_\textup{up}&= \lambda u_\textup{up},\\
			H_\textup{S}(A, \bm{W_\textup{down}})u_\textup{down}&= \lambda u_\textup{down},
			\end{system}
	\end{equation*}
	where $\bm{W_\textup{up}}=\bm{W_\textup{up}^{(1)}} + \bm{W_\textup{up}^{(2)}}$ with 
	\begin{equation*}
		\bm{W_\textup{up}^{(1)}}:= -\frac{i}{2} a^\ast \cdot \bm{B} \cdot a + V^{(1)}\bm{I_{\C^{n(d)/2}}}
		\qquad \text{and} \qquad
		\bm{W_\textup{up}^{(2)}}:= V^{(2)}\bm{I_{\C^{n(d)/2}}},
	\end{equation*}
	and where $\bm{W_\textup{down}}=\bm{W_\textup{down}^{(1)}} + \bm{W_\textup{down}^{(2)}}$ with 
	\begin{equation*}
		\bm{W_\textup{down}^{(1)}}:= -\frac{i}{2} a \cdot \bm{B} \cdot a^\ast + V^{(1)}\bm{I_{\C^{n(d)/2}}}
		\qquad \text{and} \qquad
		\bm{W_\textup{down}^{(2)}}:= V^{(2)}\bm{I_{\C^{n(d)/2}}}.
	\end{equation*}
	Notice that here we have also used that the component $\bm{V^{(1)}}$ and $\bm{V^{(2)}}$ of $\bm{V}=\bm{V^{(1)}}+\bm{V^{(2)}}$ are diagonal by the hypothesis.
	
	It is easy to see that 
	\begin{equation*}
		\int_{\R^d} \abs{x}^2\abs{\bm{W_\textup{up}^{(1)}}}^2 \abs{u_\textup{up}}^2\, dx\leq \Big(a + \frac{d}{2} c\Big)^2 \int_{\R^d} \abs{\nabla_{\!A}u_\textup{up}}^2\, dx, 
	\end{equation*}
	and 
	\begin{equation*}
		\int_{\R^d} \abs{x}^2\abs{\bm{W_\textup{down}^{(1)}}}^2 \abs{u_\textup{down}}^2\, dx\leq \Big(a + \frac{d}{2} c\Big)^2 \int_{\R^d} \abs{\nabla_{\!A}u_\textup{down}}^2\, dx, 
	\end{equation*}
	where we have used~\eqref{eq:conditions_B_V} for the vector $(u_\textup{up}, 0)$ and $(0, u_\textup{down}),$ respectively, and the fact that $\abs{a}=\sqrt{d}.$ 
	
	This means that we are in the hypotheses of Theorem~\ref{thm:Schro} (once we replace $\bm{V}$ with $\bm{W_\textup{up}}$ and $\bm{W_\textup{down}}$ and with $a + \frac{d}{2}c$ instead of $a_2$ in~\eqref{hyp:V^1}) and therefore $H_\textup{S}(A, \bm{W_\textup{up}})$ and $H_\textup{S}(A, \bm{W_\textup{down}})$ have no eigenvalues. Hence $u_\textup{up}$ and $u_\textup{down}$ are vanishing and with them $u=(u_\textup{up}, u_\textup{down}).$  

\medskip
This concludes the proof of Theorem~\ref{Thm.Pauli}.
\end{proof}

\subsection{Absence of eigenvalues of Dirac operators in any dimension}
\label{Sec:Dirac}

Now we can conclude our discussion by proving the absence of eigenvalues of Dirac operators in the general case, namely proving Theorem~\ref{thm:magnetic_Dirac}.

Let us start commenting on the odd-dimensional case. Due to expression~\eqref{eq:oddDirac} for the squared Dirac in odd dimensions and due to the analogy with~\eqref{favorable_form} in the three-dimensional case, one can proceed as in the proof of Theorem~\ref{thm:magnetic_Dirac3d} using the validity of the corresponding result Theorem~\ref{Thm.Pauli} for Pauli operators to get the result.

Turning to the even-dimensional situation, one realises from~\eqref{eq:evenDirac} that the squared Dirac operator equals a shifted Pauli operator. 
Therefore Theorem~\ref{thm:magnetic_Dirac} follows as a consequence of Theorem~\ref{Thm.Pauli} for even Pauli operators.

\subsection*{Acknowledgment}
The first author (L.C.) gratefully acknowledges financial support by the Deutsche Forschungsgemeinschaft (DFG) through CRC 1173.
The research of the second author (D.K.) was partially supported 
by the GACR grant No.\ 18-08835S.	
%

	
\providecommand{\bysame}{\leavevmode\hbox to3em{\hrulefill}\thinspace}
\providecommand{\MR}{\relax\ifhmode\unskip\space\fi MR }
\providecommand{\MRhref}[2]{%
  \href{http://www.ams.org/mathscinet-getitem?mr=#1}{#2}
}
\providecommand{\href}[2]{#2}

\end{document}